%% file: sextic0.tex
\documentclass[11pt]{article}

\usepackage{graphicx}

\newtheorem{thm}{Theorem}
\newtheorem{prop}{Proposition}
\newtheorem{conj}{Conjecture}
\newtheorem{fact}{Fact}

\newcommand{\pf}{\emph{Proof.}$\ $}
\newcommand{\fp}{\hfill $\triangle$ \\[5pt]}

\newcommand{\CC}[1]{\ensuremath{\mathbf{C}^{#1}}}
\newcommand{\CP}[1]{\ensuremath{\mathbf{CP}^{#1}}}

\newcommand{\R}[1]{\ensuremath{\mathbf{R}^{#1}}}
\newcommand{\RP}[1]{\ensuremath{\mathbf{RP}^{#1}}}

\newcommand{\RR}[2]{\ensuremath{\mathcal{R}_{\overline{#1} #2}}}

\newcommand{\Cu}[1]{\ensuremath{\mathcal{C}_{#1}}}
\newcommand{\Cb}[1]{\ensuremath{\mathcal{C}_{\overline{#1}}}}
\newcommand{\CCu}[1]{\ensuremath{{C}_{#1}}}
\newcommand{\CCb}[1]{\ensuremath{{C}_{\overline{#1}}}}

\newcommand{\LL}{\ensuremath{\mathcal{L}}}

\newcommand{\Orb}[1]{\ensuremath{\mathcal{O}_{#1}}}

\newcommand{\A}[1]{\ensuremath{\mathcal{A}_{#1}}}

\newcommand{\Sym}[1]{\ensuremath{\mathcal{S}_{#1}}}

\newcommand{\D}[1]{\ensuremath{\mathcal{D}_{#1}}}

\newcommand{\DD}[2]{\ensuremath{\mathcal{D}_{\Bar{#1}#2}}}

\newcommand{\G}{\ensuremath{\mathcal{G}}}

\newcommand{\Tet}{\ensuremath{\mathcal{T}}}

\newcommand{\Oct}{\ensuremath{\mathcal{O}}}

\newcommand{\I}{\ensuremath{\mathcal{I}}}
\newcommand{\Iu}[1]{\ensuremath{\mathcal{I}_{#1}}}
\newcommand{\Ib}[1]{\ensuremath{\mathcal{I}_{\overline{#1}}}}

\newcommand{\V}{\ensuremath{\mathcal{V}}}
\newcommand{\Va}{\ensuremath{\mathcal{V}_{3 \cdot 360}}}
\newcommand{\Vb}{\ensuremath{\mathcal{V}_{6 \cdot 360}}}

\newcommand{\Z}[1]{\ensuremath{\mathbf{Z}/{#1}}}

\newcommand{\St}[1]{\ensuremath{\mathrm{Stab}{#1}}}

\newcommand{\ZZ}[2]{\ensuremath{Z_{\Bar{#1} #2}}}

\newcommand{\PP}[2]{\ensuremath{P_{\overline{#1} #2}}}

\newcommand{\id}{\ensuremath{\mathrm{id}}}

\newcommand{\Bub}{\ensuremath{\mathrm{bub}}}
\newcommand{\bub}[2]{\ensuremath{\mathrm{bub}_{\overline{#1} #2}}}

\newcommand{\Ub}[1]{\ensuremath{U_{\overline{#1}}}}

\newcommand{\pZ}[2]{\ensuremath{p_{\:\!\overline{#1}#2}}}
\newcommand{\LZ}[2]{\ensuremath{\mathcal{L}_{\overline{#1}#2}}}
\newcommand{\LLZ}[2]{\ensuremath{L_{\overline{#1}#2}}}

\newcommand{\pT}[2]{\ensuremath{p_{#1* #2}}}
\newcommand{\LT}[2]{\ensuremath{\mathcal{L}_{#1* #2}}}

\newcommand{\pTbar}[2]{\ensuremath{p_{\:\!\overline{#1* #2}}}}
\newcommand{\LTbar}[2]{\ensuremath{\mathcal{L}_{\overline{#1* #2}}}}

\newcommand{\pQ}[3]{\ensuremath{p_{\:\!\overline{#1}#2_{#3}}}}
\newcommand{\LQ}[3]{\ensuremath{\mathcal{L}_{\overline{#1}#2_{#3}}}}

\newcommand{\pD}[2]{\ensuremath{p_{\:\!\overline{#1} #2}}}
\newcommand{\LD}[2]{\ensuremath{\mathcal{L}_{\overline{#1} #2}}}
\newcommand{\LLD}[2]{\ensuremath{L_{\overline{#1} #2}}}

\newcommand{\pP}[3]{\ensuremath{p_{\:\!\overline{#1} #2_{#3}}}}
\newcommand{\LP}[3]{\ensuremath{\mathcal{L}_{\overline{#1} #2_{#3}}}}

\newcommand{\F}[1]{\ensuremath{F_{#1}}}

\newcommand{\h}[1]{\ensuremath{h_{#1}}}

\newcommand{\Bar}[1]{\ensuremath{\overline{#1}}}

\title{Solving the sextic by iteration:\\ 
A study in complex geometry and dynamics}

\author{
Scott Crass\\
Department of Mathematics\\
Buffalo State College\\
Buffalo, NY  14222\\
\texttt{crasssw@buffalostate.edu}}

\begin{document}

\maketitle

\input{sextic2}
\input{sextic3}

\input{sextic4}

\input{sextic5}

\input{sextic6}
\input{sextic7}
\input{sextic8}

\input{sextic9}

\appendix

\include{appA}
\include{appB}

\include{bib}
\end{document}

%% file: sextic2.tex
\section{Introduction}\footnotetext{For insight and inspiration, I would
like thank Peter Doyle and Curt McMullen, the founders of this project.}
\subsection{Overview}

Recently, \cite{DM} devised an iterative
solution to the fifth degree polynomial.  At the
method's core is a rational mapping $f$ of \CP{1}\ with
the icosahedral symmetry of a general quintic.
Algebraically, this means that $f$ commutes with a group
of M\"{o}bius transformations that is isomorphic to the
alternating group \A{5}. Moreover, this
\A{5}-\emph{equivariant} posseses
\emph{nice} dynamics:  for almost any initial point
$a\in\CP{1}$, the sequence of iterates
$f^k(a)$ converges to one of the periodic cycles that
comprise an icosahedral orbit.\footnote{For a geometric
description see
\cite[p. 163]{DM}.}  This breaking of
\A{5}-symmetry provides for a \emph{reliable} or
\emph{generally-convergent} quintic-solving algorithm: 
with almost any fifth-degree equation, associate a
rational mapping that has nice dynamics and whose
attractor consists of a single orbit from which one
computes a root.  

An algorithm that solves the sixth-degree equation calls
for a dynamical system with \Sym{6} or \A{6} symmetry. 
Since neither \Sym{6} nor \A{6} acts on
\CP{1}, attention turns to higher dimensions.  Acting on
\CP{2}\ is an \A{6}-isomorphic group of projective
transformations found by Valentiner in the late
nineteenth century.  The present work exploits this
2-dimensional \A{6}\ ``soccer ball" in order to discover
a ``Valentiner-symmetric" rational mapping of \CP{2}\ 
whose dynamics
\emph{experimentally appear} to be nice in the above
sense---transferred to the \CP{2}\ setting.  This map
provides the central feature of a conjecturally-reliable
sextic-solving algorithm analogous to that employed in
the quintic case.

\subsection{Solving equations by iteration}

For $n \leq 4$, the symmetric groups \Sym{n}\ act
faithfully on \CP{1}.  Corresponding to each action is a
map whose nice dynamics provides for algorithmic
convergence to roots of given
$n$th-degree equations.  For instance, Newton's method
provides a direct iterative solution to quadratic
polynomials, but, due to a lack of symmetry, not to
higher degree equations.  My interests here are the
geometric and dynamical properties of  complex
projective mappings rather than numerical estimates. 

The search for elegant complex geometry and dynamics
continues into degree five where \A{5}\ is the
appropriate group, since \Sym{5}\ fails to act on the
sphere.  This reduction in the galois group requires the
extraction of the square root of a polynomial's
discriminant.  Such root-taking is itself the result of
a reliable iteration, namely, Newton's method.  In
practical terms, the Doyle-McMullen algorithm solves a
family of fifth-degree \emph{resolvents} the members of
which possess \A{5}\ symmetry.  A map with icosahedral
symmetry and nice dynamics plays the leading role.

Pressing on to the sixth-degree leads to the 2-dimensional
\A{6} action of the Valentiner group
\V.  Here, the problem shifts to one of finding a nice
\V-symmetric mapping of \CP{2} from whose attractor one
calculates a given sextic's root.\footnote{The
solution-procedure follows that of the quintic algorithm. 
See Section~\ref{sec:solve}.}  Providing the overall
framework is the 2-dimensional \A{6}\ analogue of the
icosahedron.

\subsection{Proofs and Computations}

At the moment, many of this work's results have only
computational support.  As such, I call them ``Facts". 
Furthermore, its conjectural nature calls for a deeper
understanding of Valentiner geometry and dynamics.  As the
theory of complex dynamics in several dimensions develops
more sophisticated weaponry, the barricades to
understanding might become assailable.  For now, I hope
that these discoveries provide a stimulus to such
development.  

\section{Valentiner's Group:  The \A{6} Action on
\CP{2}}

\begin{quote}
\small

If \dots any equation $f(x)=0$ is given, we will
investigate what is the smallest number of variables with which we
can construct a group of linear substitutions which is isomorphic
with the Galois group of $f(x)=0$.~\cite[p. 138]{Klein}

\end{quote}
\normalsize

In the wake of the
mid-nineteenth century abstractionist turn in mathematics the theory
of group representations began to emerge.  Part
of the concrete yield from work on symmetric
groups was Valentiner's discovery\cite{Val} of a complex
projective group that is isomorphic to \A{6}---the
alternating group of six things. 
Shortly thereafter \cite{Wiman} explored some
of the geometric and invariant structure
determined by this action on
the complex projective plane.   A
more thorough exposition appeared in \cite{Fricke}.  

Here, I take a new approach to the generation of this
``Valentiner" group and then explore some of its rich
combinatorial geometry.  The core of this work involves
the development of a combinatorially sensitive
description of the basic geometric structures.  In so
doing, I reproduce some of the Wiman and Fricke
results.  The study culminates in an account of the
system of Valentiner-invariant polynomials and,
thereby, lays the algebraic foundation for constructing
Valentiner-symmetric mappings of \CP{2}.

\subsection{Basics of \A{6}}

The following table indicates the non-identity
elements of \A{6}\ by order, cycle-structure, and
cardinality.

$$ \begin{array}{c|l|l}

\mbox{Order}&\mbox{Structure}&\mbox{Number}\\
\hline
&&\\
2&(ab)(cd)&45=\frac{1}{2}\left(\begin{array}{c}
    6\\4 \end{array}\right) \cdot 3!\\[15pt]

 3&(abc)&40=\left(\begin{array}{c} 6\\3
\end{array}\right)\cdot
 2\\[15pt]
 3&(abc)(def)&40=\left(\begin{array}{c} 6\\3
\end{array}\right)\cdot
 2\\[15pt]

 4&(abcd)(ef)&90=\left(\begin{array}{c}
6\\4\end{array}\right)\cdot
 3!\\[15pt]
 5&(abcde)&144=6\cdot 4!

\end{array} $$

Sitting inside \A{6}\ are twelve versions of \A{5}
that decompose into two conjugate systems of six:
\begin{itemize}
\item[1)] the stabilizers \St{\{k\}}\ of one thing
\item[2)] the permutations of the six pairs of antipodal
icosahedral vertices under rotation.
\end{itemize}
Acting by conjugation, \A{6} permutes each of
the two systems individually.  A given \A{5}\ subgroup fixes
itself set-wise and permutes its five conjugates according to the
rotational icosahedral group's action on the five cubes found in the
icosahedron.  Meanwhile, the other system of six
subgroups undergo the permutations of the six pairs
of antipodal vertices.  Consequently, the
intersection of two \A{5}\ subgroups in the \emph{same} system is
isomorphic to \A{4}---the tetrahedral rotations---while two in
\emph{different} systems give a dihedral group \D{5}.

%% file: sextic3.tex
\subsection{Generating the Valentiner group} 
\label{sec:GenVal}

An \A{5}\ subgroup of \A{6}, say \St{\{1\}}, extends to
\A{6}\ by addition of the generator (12)(3456).  Furthermore,
this order-four element generates an \Sym{4}\ over the
\A{4}\ subgroup
$$ \left<(35)(46),(456)\right> \subset \St{\{1\}}. $$ 
This structure suggests a method for producing an
\A{6}-isomorphic group
\V\ in $\mbox{PGL}_{3}(\CC{})$:  
\begin{itemize}
\item[1)]  take a tetrahedral subgroup \Tet of an
icosahedral group \I;
\item[2)]  by addition of an order-4 transformation $Q$,
extend \Tet\ to an octahedral group $\Oct =
\left<\Tet,Q\right>$;
\item[3)]  generate $\V = \left<\I,Q\right>\simeq \A{6}$.
\end{itemize}

\begin{figure}[hp]

\resizebox{\textwidth}{!}{\includegraphics{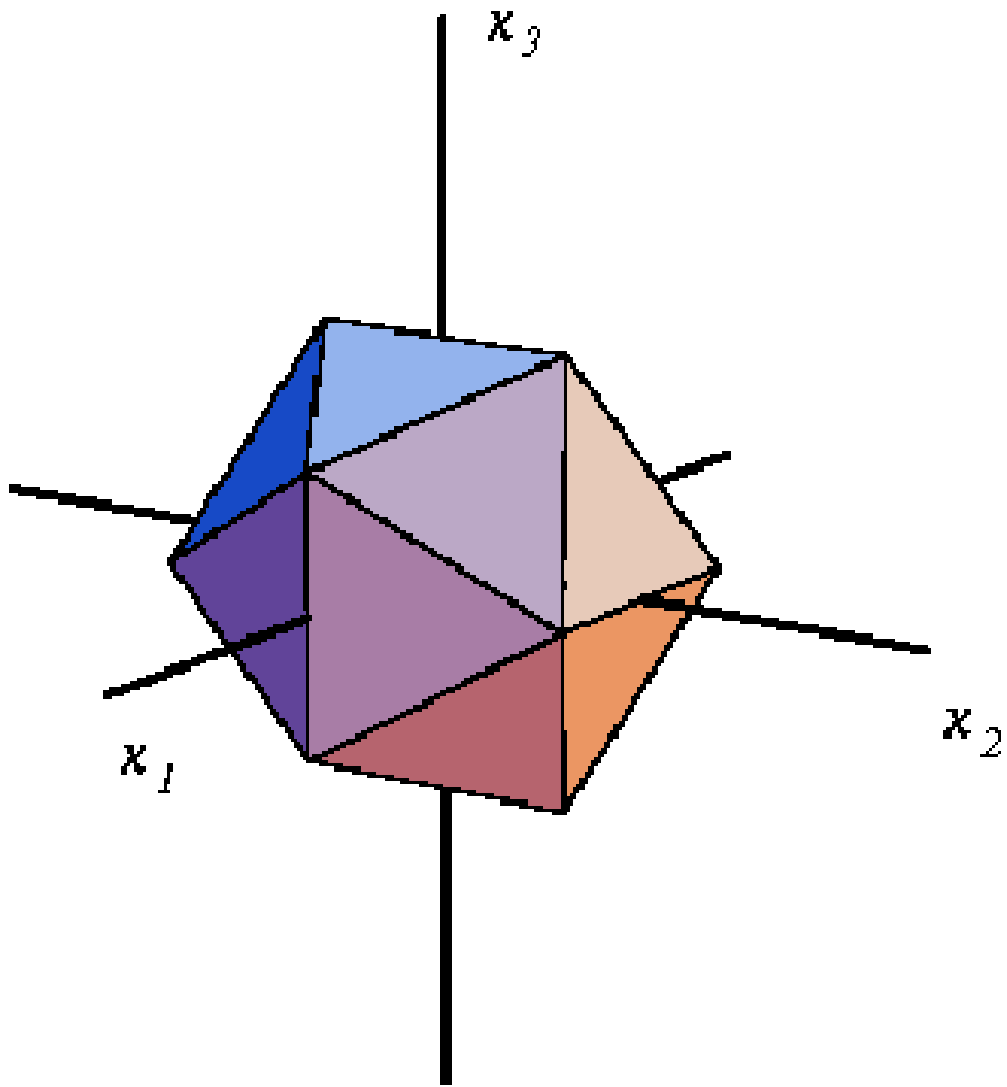}}

\caption{The icosahedron in ``octahedral" coordinates.}
\label{fig:IcosInOct}

\end{figure}

The 15 pairs of antipodal edges of the standard icosahedron
decompose into five triples such that three lines joining
antipodal edge-midpoints are mutually perpendicular. 
Stabilizing each such triple is one of the five tetrahedral
subgroups of the icosahedral group.  Alternatively, the
lines in such a triple correspond to the two-fold rotational
axes of a tetrahedron whose four vertices are face-centers
of the icosahedron.  With such a triple of lines as
coordinate axes\footnote{See Figure~\ref{fig:IcosInOct}.} in
$\{x_{1},x_{2},x_{3}\}$, the points\footnote{Square brackets
indicate points in projective space.}
$$A=\left\{[1,1,1],[-1,-1,1],[1,-1,-1],[-1,1,-1]\right\}$$
constitute a set of tetrahedral vertices.  The corresponding
tetrahedral group 
$$\Tet=\St{(A)}$$ 
consists of the identity and the 11
\emph{orthogonal} transformations:
$$ \begin{array}{c}

\begin{array}{lll}
Z_{1}=\left( \begin{array}{ccc}
1&0&0\\0&-1&0\\0&0&-1 \end{array}
\right)&
Z_{2}= \left( \begin{array}{ccc}
-1&0&0\\0&1&0\\0&0&-1 \end{array}
\right)&
Z_{3}= \left( \begin{array}{ccc}
-1&0&0\\0&-1&0\\0&0&1 \end{array}
\right) 
\end{array}\\[.3in]

\begin{array}{ll}
T_{1}= \left( \begin{array}{ccc}
0&0&1\\1&0&0\\0&1&0 \end{array}
\right)&
T_{1}^2= \left(\begin{array}{ccc}
0&1&0\\0&0&1\\1&0&0 \end{array}
\right)\\[.3in]
T_{2}=\left( \begin{array}{ccc}
0&0&-1\\1&0&0\\0&-1&0 \end{array}
\right)&
T_{2}^2= \left( \begin{array}{ccc}
0&1&0\\0&0&-1\\-1&0&0 \end{array}
\right)\\[.3in]
T_{3}= \left( \begin{array}{ccc}
0&0&-1\\-1&0&0\\0&1&0 \end{array}
\right)&
T_{3}^2= \left( \begin{array}{ccc}
0&-1&0\\0&0&1\\-1&0&0 \end{array}
\right)\\[.3in]
T_{4}= \left( \begin{array}{ccc}
0&0&1\\-1&0&0\\0&-1&0 \end{array}
\right)&
T_{4}^2= \left( \begin{array}{ccc}
0&-1&0\\0&0&-1\\1&0&0 \end{array}
\right).
\end{array}

\end{array} $$
Being orthogonal, \Tet\ preserves the quadratic form 
$$ C(x) = x_{1}^{2} + x_{2}^{2} + x_{3}^{2} $$
and hence, the conic $\Cu{}=\{C=0\}$ in \CP{2}.  This
``tetrahedral" conic contains two sets of
orbits\footnote{This is a manifestation of the pairs of
``antipodal" tetrahedra situated in the icosahedron.} of
size four: two of the points fixed projectively by each of
the three-fold $T_{k}$.  The third fixed point is one of the
elements of the set $A$ above.  With $\rho=e^{2\pi i/3}$ the
respective points are
$$ \begin{array}{ll}
v_{1}=[\rho,\rho^2,1]  \hspace*{1in}
&v_{\overline{1}}=[\rho^2,\rho,1]\\
v_{2}=[-\rho,-\rho^2,1]\hspace*{1in}
&v_{\overline{2}}=[-\rho^2,-\rho,1]\\
v_{3}=[-\rho,\rho^2,1] \hspace*{1in}
&v_{\overline{3}}=[-\rho^2,\rho,1]\\
v_{4}=[\rho,-\rho^2,1] \hspace*{1in}
&v_{\overline{4}}=[\rho^2,-\rho,1].
\end{array} $$
The ``barred' notation $v_{\Bar{a}}$ derives from
Fricke, being suggested by an antiholomorphic relationship
between the two systems of tetrahedra.  Indeed, in the $x$
coordinates chosen above, the conjugation map $x
\rightarrow \overline{x}$ exchanges a tetrahedron
and its ``antipode".

The order-four transformation
$$ \begin{array}{c}
Q = \left( \begin{array}{ccc}
1&0&0\\0&0&\rho^{2}\\0&-\rho&0 \end{array}
\right)
\end{array} $$
cyclically permutes the $v_a$ but not the $v_{\Bar{a}}$
while 
$$  \begin{array}{c}
\Bar{Q} = \left( \begin{array}{ccc}
1&0&0\\0&0&\rho\\0&-\rho^{2}&0 \end{array}
\right)
\end{array} $$
cyclically permutes the $v_{\Bar{a}}$  but not the
$v_a$.  Also, $Q^2=\Bar{Q}^2=Z_1$.  Since
$\mbox{PGL}_{3}(\CC{})$ is four-times transitive, $Q$ and
$R$ are the unique such projective transformations. 
Accordingly, the groups
$\Oct=\left<\Tet,Q\right>$ and
$\Bar{\Oct}=\left<\Tet,\Bar{Q}\right>$ are octahedral, i.e.,
isomorphic to \Sym{4}.

\begin{figure}[hp]

\resizebox{\textwidth}{!}{\includegraphics{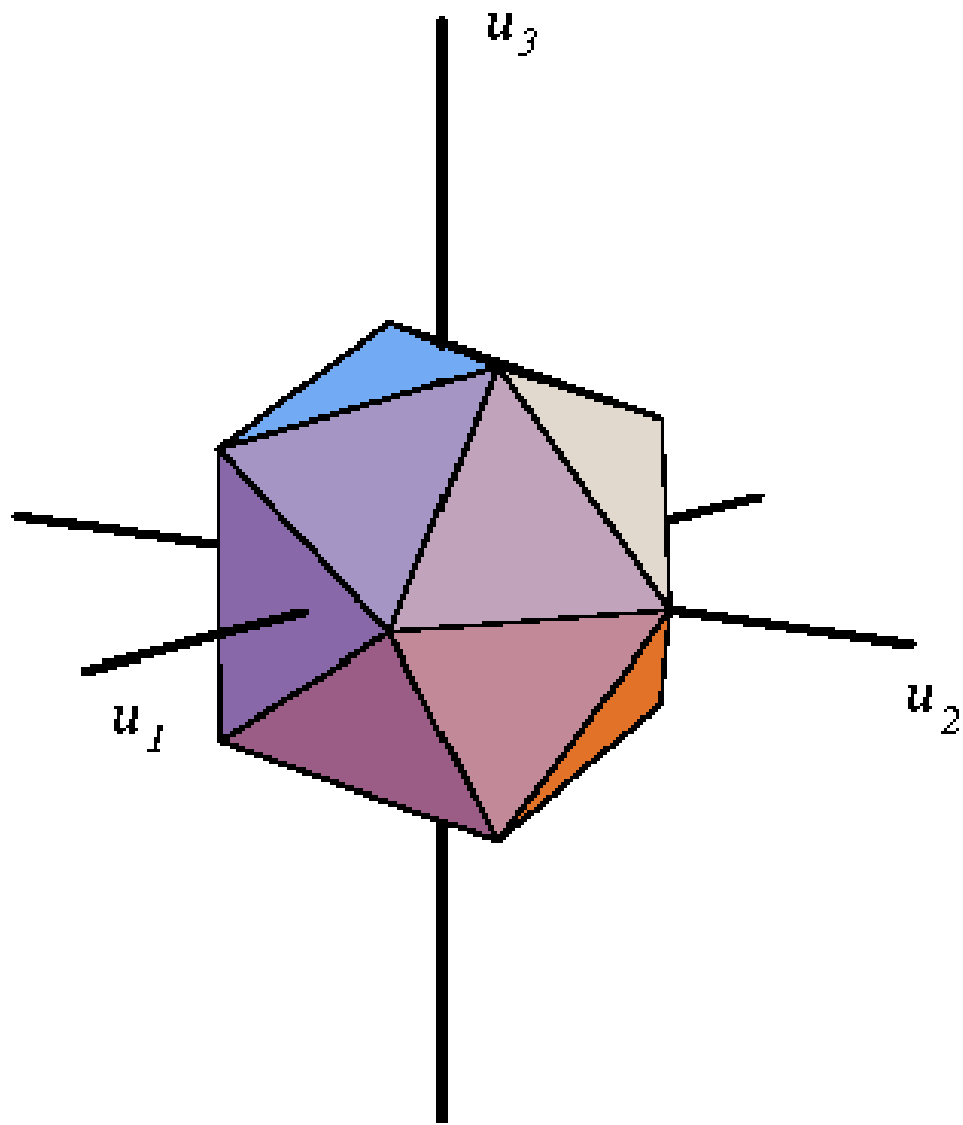}}

\caption{The icosahedron in ``icosahedral"coordinates.}
\label{fig:IcosInIcos}

\end{figure}

Extending \Tet\ to an icosahedral group \I\
requires a projective
transformation $P$ of order five that preserves
\Cu{}\ and point-wise
fixes a pair of antipodal icosahedral vertices. 
One way of producing
such a $P$ is to ``turn" the
icosahedron\footnote{See
Figure~\ref{fig:IcosInIcos}.} of
Figure~\ref{fig:IcosInOct} so that a pair of
antipodal vertices
corresponds to the point $[0,1,0]$, i.e., to the
``affine" points $(0,\pm 1,0)$.
In these ``icosahedral" coordinates
$\{u_{1},u_{2},u_{3}\}$ the desired
transformation of order five is
$$ 
P_{u} = \left( \begin{array}{ccc}
\cos \frac{2\pi}{5}&0&-\sin \frac{2\pi}{5}\\
 0&1&0\\
  \sin \frac{2\pi}{5}&0&\cos \frac{2\pi}{5}
\end{array} \right). 
$$ The change of basis from octahedral to icosahedral
coordinates is
$$
u = Ax = \left( \begin{array}{ccc}
c&-s&0\\s&c&0\\0&0&1
\end{array} \right) 
 \left(\begin{array}{c} x_1\\x_2\\x_3
\end{array}\right) 
$$
where 
$$
c=\sqrt{\frac{5+\sqrt{5}}{10}} \hspace*{20pt}
s=\sqrt{\frac{5-\sqrt{5}}{10}}.
$$
Thus, in icosahedral coordinates, the conic form
preserved by \I\ is
$$ C(u)=u_1^2+u_2^2+u_3^2 $$
while $P$ has the expression
$$ 
P_{x}\ =\ A^{-1} P_{u} A\ =\ \frac{1}{2}
\left( \begin{array}{ccc}
1        &\tau^{-1}&-\tau\\
\tau^{-1}&\tau     &1\\
\tau     &-1       &\tau^{-1}
\end{array} \right),\hspace*{.2in}
\tau=\frac{1+\sqrt{5}}{2}.
$$
Finally, $\I = \left<\Tet,P\right>$.  This produces two
Valentiner groups distinguished by chirality:
$$ 
\V = \left<\I, Q\right> \hspace*{20pt}
\Bar{\V} = \left<\I, \Bar{Q}\right>. 
$$

Use of the terms `octahedral coordinates' and `icosahedral
coordinates' follows that of \cite[pp.
263ff]{Fricke}.  His system of octahedral generators are
nearly those above.  In his icosahedral coordinates, the
conic form is
$$ C_{Fricke}(z) = z_{1} z_{3} + z_{2}^2. $$
The change of coordinates $B$ that yields $C(Bz)
= C_{Fricke}(z)$ is
$$ 
B = \left( \begin{array}{ccc}
 1&0&i\\0&1&0\\1&0&-i
\end{array} \right). 
$$
With $\epsilon = e^{2\pi i/5}$ and $Z=Z_{2}$, the generators
$Z_{u}= A Z_{x} A^{-1}$ and $P_{u}=A P_{x} A^{-1}$ become his
icosahedral generators~\cite[p. 263]{Fricke} $T$ and $S$:
$$ \begin{array}{lclcl}

Z_{z}&=&B Z_{u} B^{-1}&=&\frac{1}{\sqrt{5}} \left(
\begin{array}{ccc}
-\frac{1+\sqrt{5}}{2}&2&\frac{-1+\sqrt{5}}{2}\\[5pt]
1&1&1\\[5pt]
\frac{-1+\sqrt{5}}{2}&2&-\frac{1+\sqrt{5}}{2}
\end{array} \right)\\[.4in]

&&&=&\frac{1}{\sqrt{5}} \left( \begin{array}{ccc}
\epsilon^2 + \epsilon^3&2&\epsilon + \epsilon^4\\
1&1&1\\
\epsilon + \epsilon^4&2&\epsilon^2 + \epsilon^3
\end{array} \right)\\[.4in] 
P_{z}&=&B P_{u} B^{-1}&=&\left(
\begin{array}{ccc}
\epsilon&0&0\\ 0&1&0\\ 0&0&\epsilon^4
\end{array} \right).

\end{array} $$


%% file: sextic4.tex
\subsection{Valentiner geometry}

\subsubsection{Icosahedral conics}

Since a ternary icosahedral group has an orthogonal
representation on
\R{3}, its complexification stabilizes the quadratic form
$$ C(x)=x_1^2 + x_2^2 + x_3^2  $$ 
over \CC{}\ and hence, the conic $\Cu{}=\{C(x)=0\}$ in
\CP{2}.  For reasons that will become clear, the two
systems of six icosahedral groups \Ib{k}\ and
\Iu{k}\ in \V\ receive the designations ``barred" and
``unbarred".  Corresponding to each of the twelve
icosahedral subgroups of \V\ is a quadratic form
\CCb{a} or \CCu{a} and conic \Cb{a} or \Cu{a} that \Ib{a}
or \Iu{a}\ preserves.  Thereby, each such conic possesses
the structure of an icosahedron.

Call the form above \CCb{1}.  The action of
\V\ produces the five
remaining forms.  With $\eta=(3+\sqrt{15}\,i)/4$,
$$ \begin{array}{lll}
\CCb{1}(x)
&=&x_{1}^{2} + x_{2}^{2} + x_{3}^{2}\\[10pt]
\CCb{2}(x)
&=&\CCb{1}(Q^{-1}x)\\[5pt]
&=&x_{1}^{2} + \rho^{2} x_{2}^{2} + \rho\, x_{3}^{2}\\[10pt]
\CCb{3}(x)&=&\CCb{2}(P^{-4}x)\\[5pt]
&=&\eta\left[\frac{\eta}{3}\,
 (x_{1}^{2} + \rho\, x_{2}^{2} + \rho^{2} x_{3}^{2})
  + (\rho^{2}\, x_{1} x_{2} + \rho\, x_{1} x_{3} - x_{2}
x_{3})\right]\\[10pt]
\CCb{4}\,(x)
&=&\CCb{2}(P^{-3}x)\\[5pt]
&=&\eta\left[\frac{\eta}{3}\,
 (\rho\, x_{1}^{2} + \rho^{2} x_{2}^{2} +  x_{3}^{2})
  + (-x_{1} x_{2} + \rho^{2} x_{1} x_{3} + \rho\, x_{2}
x_{3})\right]\\[10pt]
\CCb{5}(x)
&=&\CCb{2}(P^{-2}x)\\[5pt]
&=&\eta\left[\frac{\eta}{3}\,
 (\rho\, x_{1}^{2} + \rho^{2} x_{2}^{2} +  x_{3}^{2})
  - (x_{1} x_{2} + \rho^{2} x_{1} x_{3} + \rho\, x_{2}
x_{3})\right]\\[10pt]
\CCb{6}(x)
&=&\CCb{2}(P^{-1}x)\\[5pt]
&=&\eta\left[\frac{\eta}{3}\,
 (x_{1}^{2} + \rho\, x_{2}^{2} + \rho^{2} x_{3}^{2})
  + (\rho^{2}x_{1} x_{2} - \rho\, x_{1} x_{3} + x_{2}
x_{3})\right].
\end{array} $$
As for the action of \V\ on these:
$$ \begin{array}{rl}
P:&\!
 \CCb{1}\leftrightarrow\CCb{1} \hspace*{20pt}
  \CCb{2}\rightarrow\CCb{6}\rightarrow\CCb{5}\rightarrow
   \CCb{4}\rightarrow\CCb{3}\rightarrow\CCb{2}\\[5pt]
Z=Z_2:&\!
 \CCb{1}\leftrightarrow\CCb{1}\hspace*{20pt}
  \CCb{2}\leftrightarrow\CCb{2}\hspace*{20pt}
   \CCb{3}\leftrightarrow\rho^2\,\CCb{4}\hspace*{20pt}
    \CCb{5}\leftrightarrow\rho\,\CCb{6}\\[5pt]
Q:&\!
 \CCb{1}\leftrightarrow\CCb{2}\hspace*{20pt} 
\CCb{3}\rightarrow\CCb{6}\rightarrow\rho^2\,\CCb{5}\rightarrow
   \rho^2\,\CCb{4}\rightarrow\CCb{3}.
\end{array} $$

Direct calculation yields
\begin{prop}
The quadratic forms \CCb{m} are linearly independent and so
span the six dimensional space of ternary quadratic forms. 
In particular, an unbarred conic form
\CCu{a} is a linear combination of the \CCb{m} that is
invariant under the icosahedral group \Iu{a}.
\end{prop}
To be specific, the ``5-cycle" $P$ belongs to
\Ib{1}\ and to an unbarred icosahedral group---say \Iu{3}\
so that the
indexing agrees with that of Fricke.  (Recall that
the intersection of two non-conjugate \A{5}\
subgroups of \A{6}\ is a \D{5}.)  Relative or
\emph{projective} invariance under 
$P$ requires \CCu{3}\ to take the form
$$ 
\CCu{3}\ =\ 
 \alpha\,\CCb{1}+ \CCb{2}+\CCb{3}+\CCb{4}+\CCb{5}+\CCb{6}. 
$$
To determine the constant, apply to the \CCb{m}\ an element
$T$ of \Iu{3}\ that does not belong to \Ib{1}, e.g., one of
the 20 elements of order three in \Iu{3}.  Recalling the
association between
\Iu{3}\ and the permutations of the six pairs of antipodal
icosahedral vertices labeled according to the action of
$P$, such a transformation corresponds to a double 3-cycle
of the form $(abc)(def)$.  Specifically, the permutation
$(\Bar{164})(\Bar{235})$---to be used below---corresponds
to an element of \Iu{3}.  

The action of the generators on the \emph{conics} \Cb{m} is
given by the permutation of the indices:
$$ \begin{array}{rl}

P:& (\Bar{26543})\\[5pt]
Z:&(\Bar{34})(\overline{56})\\[5pt]
Q:& (\Bar{12})(\overline{3654}).

\end{array} $$
Computation in \A{6}\ yields the correspondence
$$ T=QP^2QPQ^3 : (\overline{164})(\overline{235}). $$
Moreover, the action on the conic forms is
$$ \begin{array}{ll}

T:&
\CCb{1}\rightarrow \CCb{6}\rightarrow
\rho\,\CCb{4}\rightarrow \CCb{1}\\[10pt]

&\CCb{2}\rightarrow \rho\,\CCb{3}\rightarrow \rho^2\,\CCb{5}\rightarrow \CCb{2}

\end{array} $$
so that
$$ \begin{array}{lll}

\CCu{3}(T^{-1}x)&=&\rho^{2}\, \CCb{1}+\rho\
 (\CCb{2}+\CCb{3}+\CCb{4}+\CCb{5})+\alpha\,\CCb{6}\, .

\end{array} $$
Projective invariance under $T$ of the
conic $\Cu{3}=\{\CCu{3}=0\}$ requires $\alpha=\rho$. 
Accordingly,
$$ \begin{array}{lll}
\CCu{3}(T^{-1}x)&=&\rho^{2}\, \CCb{1}+\rho\,
 (\CCb{2}+\CCb{3}+\CCb{4}+\CCb{5}+\CCb{6})\\[10pt]
 &=&\rho\, \CCu{3}(x).
\end{array} $$

Just as the barred forms stem from \CCb{1}, the remaining unbarred conic
forms\footnote{Again, the indices are chosen to agree with Fricke's
labels.} arise from \CCu{3}:
$$ \begin{array}{lll}
\CCu{1}(x) &=& \CCu{2}(P^{-1}x)\\
&=& \CCb{1}\,+\,\rho\,\CCb{2}\,+\, 
\CCb{3}\,+\,\rho^{2}\,\CCb{4}\,+\,\CCb{5}\,+\,\rho\,\CCb{6}\\
&=&-\eta \left(
    \rho^2\,x_1^2\,+\,\frac{4}{3}\,\eta^2\,x_2^2\,+\,\rho\,x_3^2\,
     +\,2\,\rho\,(\rho - 1)\,x_1 x_3 \right)\\[10pt]

\CCu{2}(x) &=& \CCu{3}(Q^{-1}x)\\
&=& \CCb{1}\,+\,\rho\,\CCb{2}\,+\,\rho\,
\CCb{3}\,+\,\CCb{4}\,+\,\rho^{2}\,\CCb{5}\,+\,\CCb{6}\\
&=&-\eta \left(
 \rho^2\,x_1^2\,+\,\frac{4}{3}\,\eta^2\,x_2^2\,+\,\rho\,x_3^2\,
     -\,2\,\rho\,(\rho - 1)\,x_1 x_3 \right)\\[10pt]

\CCu{3}(x) &=&\rho\,\CCb{1}\,
+\,\CCb{2}\,+\,\CCb{3}\,+\,\CCb{4}\,+\,\CCb{5}\,+\,\CCb{6}\\
&=&-\eta\,\rho \left(
 \rho\,x_1^2\,+\,\rho^2\,x_2^2\,+\,\frac{4}{3}\,\eta^2\,x_3^2\,
     +\,2\,\rho\,(\rho - 1)\,x_1 x_2 \right)\\[10pt]

\CCu{4}(x) &=& \CCu{2}(P^{-3}x)\\
&=& \CCb{1}\,+\,\rho^2\,\CCb{2}\,+\, 
\CCb{3}\,+\,\rho\,\CCb{4}\,+\,\rho\,\CCb{5}\,+\,\CCb{6}\\
&=&-\eta \left(
 \rho\,x_1^2\,+\,\rho^2\,x_2^2\,+\,\frac{4}{3}\,\eta^2\,x_3^2\,
     -\,2\,\rho\,(\rho - 1)\,x_1 x_2 \right)\\[10pt]

\CCu{5}(x) &=& \CCu{2}(P^{-2}x)\\
&=& \CCb{1}\,+\,\CCb{2}\,+\,\rho^{2}\,
\CCb{3}\,+\,\CCb{4}\,+\,\rho\,\CCb{5}\,+\,\rho\,\CCb{6}\\
&=&-\eta \left(
 \frac{4}{3}\,\eta^2\,x_1^2\,+\rho\,x_2^2\,+\,\rho^2\,x_3^2\,
     +\,2\,\rho\,(\rho - 1)\,x_2 x_3 \right)\\[10pt]

\CCu{6}(x) &=& \CCu{2}(P^{-4}x)\\
&=& \CCb{1}\,+\,\CCb{2}\,+\,\rho\,\CCb{3}\,+\,\rho\,
\CCb{4}\,+\,\CCb{5}\,+\,\rho^{2}\,\CCb{6}\\
&=&-\eta \left(
 \frac{4}{3}\,\eta^2\,x_1^2\,+\rho\,x_2^2\,+\,\rho^2\,x_3^2\,
     -\,2\,\rho\,(\rho - 1)\,x_2 x_3 \right)\\[10pt].
\end{array} $$
Application of \V\ yields
$$ \begin{array}{ll}
P:&\!
 \CCu{3}\leftrightarrow\CCu{3} \hspace*{20pt}
  \CCu{1}\rightarrow\CCu{5}\rightarrow\CCu{4}\rightarrow
   \CCu{6}\rightarrow\CCu{2}\rightarrow\CCu{1}\\[5pt]
Z:&\!
 \CCu{1} \leftrightarrow \CCu{1} \hspace*{20pt}
  \CCu{2} \leftrightarrow \CCu{2} \hspace*{20pt}
   \CCu{3} \leftrightarrow \rho\,\CCu{4} \hspace*{20pt}
    \CCu{5} \leftrightarrow \CCu{6}\\[5pt]
Q:&\!
 \CCu{5} \leftrightarrow \CCu{6} \hspace*{20pt}
  \CCu{1} \rightarrow \CCu{3} \rightarrow \CCu{2}\rightarrow
   \rho\,\CCu{4} \rightarrow \CCu{1}.
\end{array} $$

\subsubsection{Antiholomorphic symmetry}

The one-dimensional icosahedral group $\G_{60}$ acts on two
sets of five tetrahedra each of which corresponds to a
quadruple of points in \CP{1}. However, no element of the
group sends the tetrahedra of one set to those of the
other.  Such an exchange occurs by means of
anti-holomorphic maps of degree one.  Of these, 15
correspond to reflections through the 15 great circles of
reflective icosahedral symmetry; the remaining 45 are the
various ``odd" compositions of the 15 basic
reflections---e.g., the antipodal map.  Extending the
holomorphic $\G_{60}$ by such an ``anti-involution"
produces the group $\overline{\G}_{120}$ of all 120
symmetries of the icosahedron.  The 15 icosahedral
reflections generate this extended group while their even
products result in $\G_{60}$. In coordinates where one of
the great circles corresponds to the real axis, the
associated anti-involution is complex conjugation---in
homogeneous coordinates:
$$ [x_{1},x_{2}] \rightarrow
[\overline{x_{1}},\overline{x_{2}}]. $$

The Valentiner analogues of the tetrahedra are the two
systems of conics.  Are there ternary
anti-involutions that exchange the barred and unbarred
conics?  If so, can they take the form
$$ 
[x_{1},x_{2},x_{3}]\rightarrow
[\overline{x_{1}},\overline{x_{2}},\overline{x_{3}}]\,? 
$$
To these questions \cite[pp. 270-1,
286-9]{Fricke} provides\footnote{See below for a
combinatorial geometric computation of this additional
symmetry.} affirmative answers.  In the current octahedral coordinates, this
``bar-unbar" map is
$$ 
\Bub[x_{1},x_{2},x_{3}] = [\ \rho^{2}\, \overline{x_{1}} -
 \rho\, \overline{x_{3}}\ ,\  -\rho\, (\rho + \tau)\,
\overline{x_{2}}\ ,\  -\rho\, \overline{x_{1}} -
\overline{x_{3}}\ ]. 
$$
As for the action\footnote{The match between \CCb{a}\ and
\CCu{a} is no accident.  Fricke used this map to dub the
unbarred conics.} on the conic forms:
$$ \begin{array}{lllllll}
\CCu{1}(\Bub(x))&=&\alpha\,\rho^{2}\,\CCb{1}(x)&\hspace*{15pt}&
\CCu{2}(\Bub(x))&=&\alpha\, \rho\,\CCb{2}(x)\\[10pt]
\CCu{3}(\Bub(x))&=&\alpha\,\CCb{3}&\hspace*{15pt}&
\CCu{4}(\Bub(x))&=&\alpha\,\CCb{4}(x)\\[10pt]
\CCu{5}(\Bub(x))&=&\alpha\,\rho^{2}\,\CCb{5}(x)&\hspace*{15pt}&
\CCu{6}(\Bub(x))&=&\alpha\,\rho\,\CCb{6}(x)
\end{array} $$
where $\alpha=(3 + \sqrt{15}\,i)/2$. 
\begin{prop}
The group $\Bar{\V}_{2\cdot 360}=\left<\V,\Bub\right>$ is a
degree two extension of \V.
\end{prop}
\pf  For $T \in \V$, $T^{\prime}=\Bub \circ T \circ \Bub$
is a projective transformation that permutes the conics
within a system.  Thereby,
$T^{\prime}$ belongs to \V.  \fp 
Concerning the form of a bub map, there are coordinate
systems in which its expression \emph{is} conjugation of
each coordinate.  While interesting in their own right,
such coordinates also yield computational benefits.  Some
of the Valentiner structure suggests a means of achieving
this diagonalization.  I will take up the topic once the
relevant framework is in place.

%% file: sextic5.tex
\subsubsection{Special orbits and a Valentiner
nomenclature}

\begin{figure}[hp]

\resizebox{\textwidth}{!}{\includegraphics{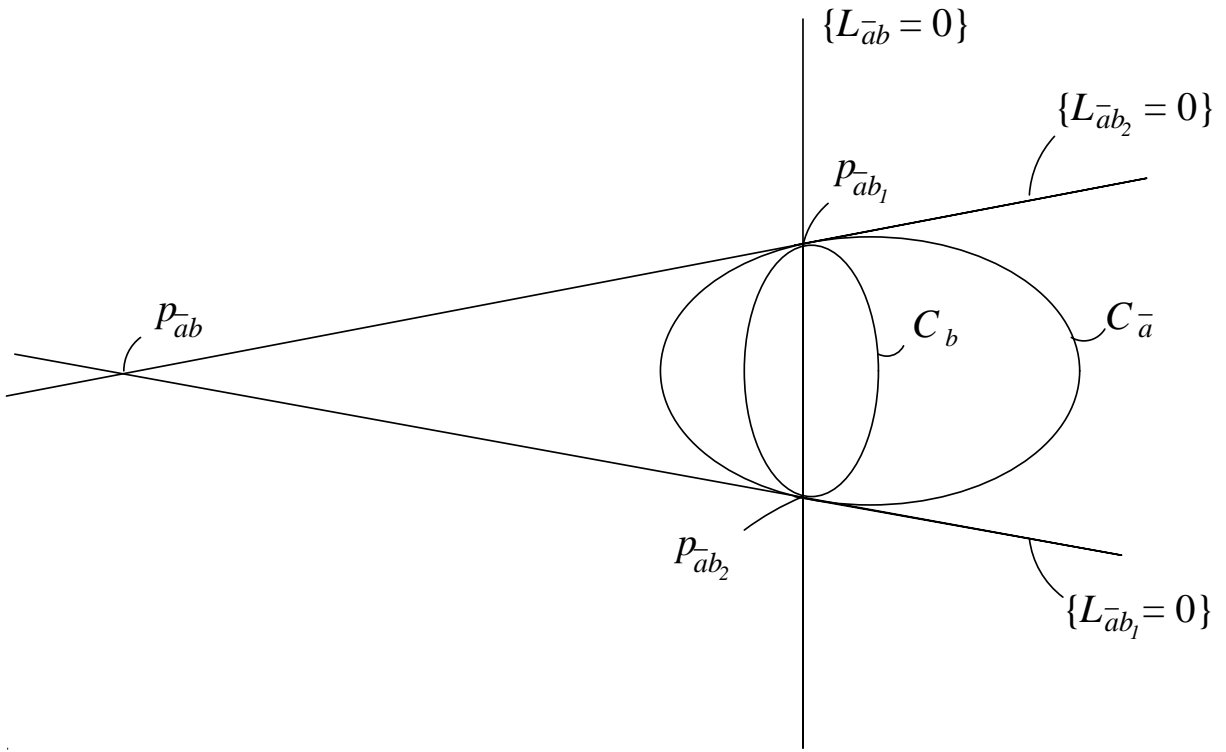}}

\caption{The triangle of one 36-point and two
72-points.}
\label{fig:36a72PointsLines}

\end{figure}

\paragraph{Intersecting conics.}
Some of the special icosahedral points on a conic
\Cb{a} occur at the intersections of \Cb{a}\ and the
other 11 conics.  
\begin{fact} Within a system, \Cb{a}\ meets each
$\Cb{m}\ (\Bar{m}\neq \Bar{a})$ in four tetrahedral
points; this gives the 20 face-centers on \Cb{a}.
\end{fact} 
The overall result is a 60 point \V-orbit
\Orb{\Bar{60}}.  Similarly, the unbarred intersections
yield \Orb{60}.  Alternatively, each member of
\Orb{\Bar{60}}\ (or \Orb{60}) is a point fixed by one
of the 20 barred (or unbarred) cyclic
subgroups\footnote{In
\A{6}, the barred-unbarred splitting manifests itself
in the two structurally distinct sets in order three: 
$(abc)$ and $(abc)(def)$.} of order three in \V.
 \begin{fact} Across systems the intersection of
\Cb{a}\ with the \Cu{b}\ gives six pairs of antipodal
icosahedral vertices
$\{\pP{a}{b}{1}, \pP{a}{b}{2}\}$.
\end{fact} 
These total to $72=6\cdot 12$ points each of
which is fixed under one of 36 order-five cyclic
groups
$\left<\PP{a}{b}\right>$.  Since an icosahedral group
is transitive on its vertices---indeed, some element
\Ib{a}\ of order two exchanges \pP{a}{b}{1}\ and
\pP{a}{b}{2}, the 72 points
$$
\Orb{72}=\{\pP{a}{b}{1}|a,b=1,\ldots,6\} \cup
\{\pP{a}{b}{2}|a,b=1,\ldots,6\} 
$$
form a \V-orbit.  A Valentiner exchange of
\pP{a}{b}{1}\ and
\pP{a}{b}{2}\ also transposes the lines\footnote{The
labeling of these lines in the sub-subscript is
arbitrary and is done so as to agree with the natural
cases in which a point does not reside on its
associated line.} \LP{a}{b}{2}\ and \LP{a}{b}{1}\
tangent to \Cu{a}\ and \Cb{b} at
\pP{a}{b}{1}\ and \pP{a}{b}{2}. Hence, the
intersection of these lines belongs to a 36 point
orbit.\footnote{See Figure~\ref{fig:36a72PointsLines}.}
Each of these\footnote{The practice of referring to
special points and lines in terms of the size of their
orbits will continue throughout.} ``36-points"
\pD{a}{b}\ corresponds to the ``36-line"
$\LD{a}{b}=\{\LLD{a}{b}=0\}$ passing through
\pP{a}{b}{1} and \pP{a}{b}{2}.  Furthermore, a
dihedral $\D{5}$ stabilizes the ``triangle"
$\{\pP{a}{b}{1}, \pP{a}{b}{2}, \pD{a}{b}\}$:
$$
\DD{a}{b}=\St{\{\pD{a}{b}\}}=\St{\{\pP{a}{b}{1},\pP{a}{b}{2}\}}
 =\left<\PP{a}{b},\ZZ{ac}{bd}\right> \simeq \D{5}.
$$
An explanation of the indices attached to the element
\ZZ{ac}{bd}\ of order two occurs below.

As for other special orbits, each of the 45
involutions $Z$ in \V\ is conjugate to
$$ \left( \begin{array}{ccc}
1&0&0\\0&-1&0\\0&0&-1 
\end{array} \right). $$
In \CP{2} each such $Z$ fixes a
point $p_{Z}$ and point-wise fixes a line
$\mathcal{L}_{Z}$.  Furthermore, $Z$ is the
square of an element
$Q$ of order four.  The three fixed points of $Q$
consist of $p_{Z}$ and two
points $(p_Q)_1$ and $(p_Q)_2$ on
$\mathcal{L}_{Z}$.  Each of the
latter points have a \Z{4}\ stabilizer and so
belong to a 90-point orbit
\Orb{90}.  The points $p_{Z}$, having $\D{4}$
stabilizers, give an
orbit \Orb{45}.  Since $Z$ acts trivially on
$\mathcal{L}_Z$, the \D{4}\ action restricted to
$\mathcal{L}_Z$ reduces to that of a Klein-four
group.  Finally, the generic points
\emph{on} a ``45-line" lie in four-point orbits
and, overall, provide \V-orbits
of size 180.

\paragraph{The configuration of 45 lines and
points.}
The intersections of  45-lines yield special orbits of
size less than 180.  Furthermore, the number of lines
meeting at such a site $p$ corresponds to the number
of involutions in the stabilizer of
$p$.  Being \D{5}-stable, a 36-point lies on five of
the 45-lines.  Similarly, four of the 45-lines meet at
a 45-point while three concur at each of the 60 and
$\overline{60}$-points. This accounts for all
intersections of the 45-lines:
$$ 36
\left(\begin{array}{c}5\\2\end{array}\right)\ +\
     45
\left(\begin{array}{c}4\\2\end{array}\right)\ +\
     60
\left(\begin{array}{c}3\\2\end{array}\right)\ +\
     60
\left(\begin{array}{c}3\\2\end{array}\right)\ =\
  \left(\begin{array}{c}45\\2\end{array}\right). 
$$

On a 45-line there are four points of each type
so that the
corresponding clusters of lines complete the
remaining set of 44 lines:
$$ 4\, (\overbrace{4}^{36-\mathrm{points}} +
\underbrace{3}_{45-\mathrm{points}}
 + \overbrace{2}^{60-\mathrm{points}} +
 \underbrace{2}_{\overline{60}-\mathrm{points}})\
=\ 44. $$

Furthermore, only one involution fixes a
90-point.  Hence, such points lie on just one of
the 45-lines.  Since a 72-point has a
\Z{5}\ stabilizer, \Orb{72}\  acquires its
exceptional status as the only special orbit
that is not a subset of the 45-lines.

\paragraph{Valentiner-speak.}
Since the combinatorial relationships among special
objects derive from those at the
level of conics, the special orbits should admit
description in a conic-based terminology.  Already
evident is the natural designation of 36,
72-points\footnote{By the Valentiner ``duality"
between special points and lines, whatever holds for
such points also holds for the associated lines. 
Thus, I will usually supress reference to one or the
other.} in terms of one barred and one unbarred conic:
$$  
\{\pP{a}{b}{1},\pP{a}{b}{2}\}\ =\ \Cb{a}\,\cap\,
\Cu{b} 
$$ 
and \pD{a}{b}\ is the pole of 
$\{\pP{a}{b}{1},\pP{a}{b}{2}\}$ with respect to
\Cb{a}\ and \Cu{b}. Turning to the 45-points, a given
involution $Z$ belongs to two tetrahedral subgroups of
the barred and unbarred icosahedral systems alike.  If
$Z$ preserves conics \Cb{a}, \Cb{b},
\Cu{c}, and \Cu{d}, then it does so
uniquely.\footnote{See below.}  Accordingly,
\ZZ{ab}{cd}, \pZ{ab}{cd}, and \LZ{ab}{cd}\ denote the
transformation, corresponding point and line
respectively.  Alternatively, the 45-line \LZ{ab}{cd}
contains the 36-points \pD{a}{c}, \pD{a}{d}, \pD{b}{c},
\pD{b}{d}\ or by duality, the 45-point \pZ{ab}{cd}
lies at the intersection of the 36-lines \LD{a}{c},
\LD{a}{d}, \LD{b}{c},
\LD{b}{d}.  Furthermore, \LZ{ab}{cd}\ is an
icosahedral axis of order two for \Cb{a}, \Cb{b},
\Cu{c}, and \Cu{d}\ while 
$$ \LZ{ab}{cd} \cap \Cb{a}, \Cb{b}, \Cu{c},\Cu{d} $$
give the antipodal pairs of edge midpoints.  Of
course, the labels for the pair of 90-points---each of
which \ZZ{ab}{cd}\ fixes---should be \pQ{ab}{cd}{1}\
and \pQ{ab}{cd}{2}.

Now, given a 45-line \LZ{ab}{cd}, which four of the
45-points belong to it?  Since there are
$15=(\!\!\begin{array}{l}{\mbox{\scriptsize
6}}\\[-7pt]{\mbox{\scriptsize 2}}\end{array}\!\!)$
choices each for the prefix \Bar{ab} and suffix $cd$,
a $15 \times 15$ array with 45 distinguished entries
depicts the configuration of 45-points and
lines.\footnote{See Figure~\ref{fig:45Array}.} Being
\V-equivalent, the rows and columns each contain three
marked entries.  Associated with each of the 15 barred
rows \Bar{ab}\ and its triple of ``45-things" indexed
by
$\{\Bar{ab}cd,\Bar{ab}ef,\Bar{ab}gh\}$ where
$\{c,d,e,f,g,h\}=\{1,\ldots,6\}$ is a tetrahedral group
$$ \Tet_{\Bar{ab}}\ =\ \Ib{a}\; \cap\; \Ib{b} $$
whose three involutions are \ZZ{ab}{cd},
\ZZ{ab}{ef}, \ZZ{ab}{gh}.  The analogous
state of affairs obtains for the unbarred columns
where the tetrahedral group
$$  \Tet_{cd}\ =\ \Iu{c}\; \cap\; \Iu{d} $$
contains involutions \ZZ{ab}{cd},  \ZZ{\imath
\jmath}{cd}, \ZZ{k \ell}{cd}. Each of these 15
tetrahedral groups extends to an
octahedral subgroup of \V
$$ 
\mathcal{O}_{\Bar{ab}} =
\St{\{\pZ{ab}{cd},\pZ{ab}{ef},\pZ{ab}{gh}\}}
\hspace*{20pt}
\mathcal{O}_{cd} =
\St{\{\pZ{ab}{cd},\pZ{\imath\jmath}{cd},\pZ{kl}{cd}\}}.
$$
Hence, the stabilizer of \pZ{ab}{cd}\ is the
intersection of octahedral groups
$$
\Oct_{\Bar{ab}}\cap
\Oct_{cd}=\St\{\pZ{ab}{cd}\}=\St\{\LZ{ab}{cd}\}\simeq \D{4}. 
$$
Furthermore, the involution \ZZ{ab}{cd}\
associates canonically with the pair of barred and
unbarred tetrahedral groups
$$ 
\Tet_{\Bar{ab}} \cap\Tet_{cd}=\left<\ZZ{ab}{cd}\right> \simeq \Z{2}.  
$$

This ``45-array" encodes a wealth of combinatorial
geometry including an answer to the query of the
preceding paragraph.\footnote{A graphical version of
this array appears in \cite[p. 542]{Wiman}.}  At a
45-point
\pZ{ab}{cd}\ there are four concurrent 45-lines whose
references have the form \LZ{ab}{ef}, \LZ{ab}{gh},
\LZ{mn}{cd}, and
\LZ{rs}{cd}.  To find these lines read along the
\Bar{ab}\ row and the
$cd$ column.  By way of example,
$$ 
\pZ{12}{34}\ \in\ \LZ{12}{12}\; \cap\;\LZ{12}{56}\;
\cap\; \LZ{36}{34}\; \cap\; \LZ{45}{34}.
$$
Duality gives
$$
\{\pZ{12}{12},\pZ{12}{56},\pZ{36}{34},\pZ{45}{34}\}
\subset\LZ{12}{34}.
$$
The \D{5}\ stabilizer of a 36-point \pD{a}{b} contains
five involutions whose indices have a prefix \Bar{a}
and a suffix $b$.  For \pD{3}{5}\ this gives
\Bar{13}35, \Bar{23}15,
\Bar{34}45, \Bar{35}56, \Bar{36}25.  Hence,
$$ \pD{3}{5}\ \in\ \LZ{13}{35}\; \cap\;
\LZ{23}{15}\; \cap\;
\LZ{34}{45}\; \cap\; \LZ{35}{56}\; \cap\;
\LZ{36}{25} $$
and 
$$
\{\pZ{13}{35},\pZ{23}{15},\pZ{34}{45},\pZ{35}{56},\pZ{36}{25}\} 
 \subset \LD{3}{5}. 
$$

The array also supplies a connection to \A{6}.  For
instance, the involution \ZZ{12}{34}\ fixes
\LZ{12}{34}\ pointwise while preserving \LZ{12}{12},
\LZ{12}{56}, \LZ{36}{34}, and
\LZ{45}{34} as sets.  Accordingly, it permutes the
conics by:
$$ \begin{array}{lll}
\Cb{3} \leftrightarrow \Cb{6} &\hspace*{20pt}&
\Cb{4} \leftrightarrow
\Cb{5}\\[15pt]
\Cu{1} \leftrightarrow \Cu{2} && \Cu{5}
\leftrightarrow \Cu{6}.
\end{array} $$

\begin{figure}[hp]

\resizebox{\textwidth}{!}{\includegraphics{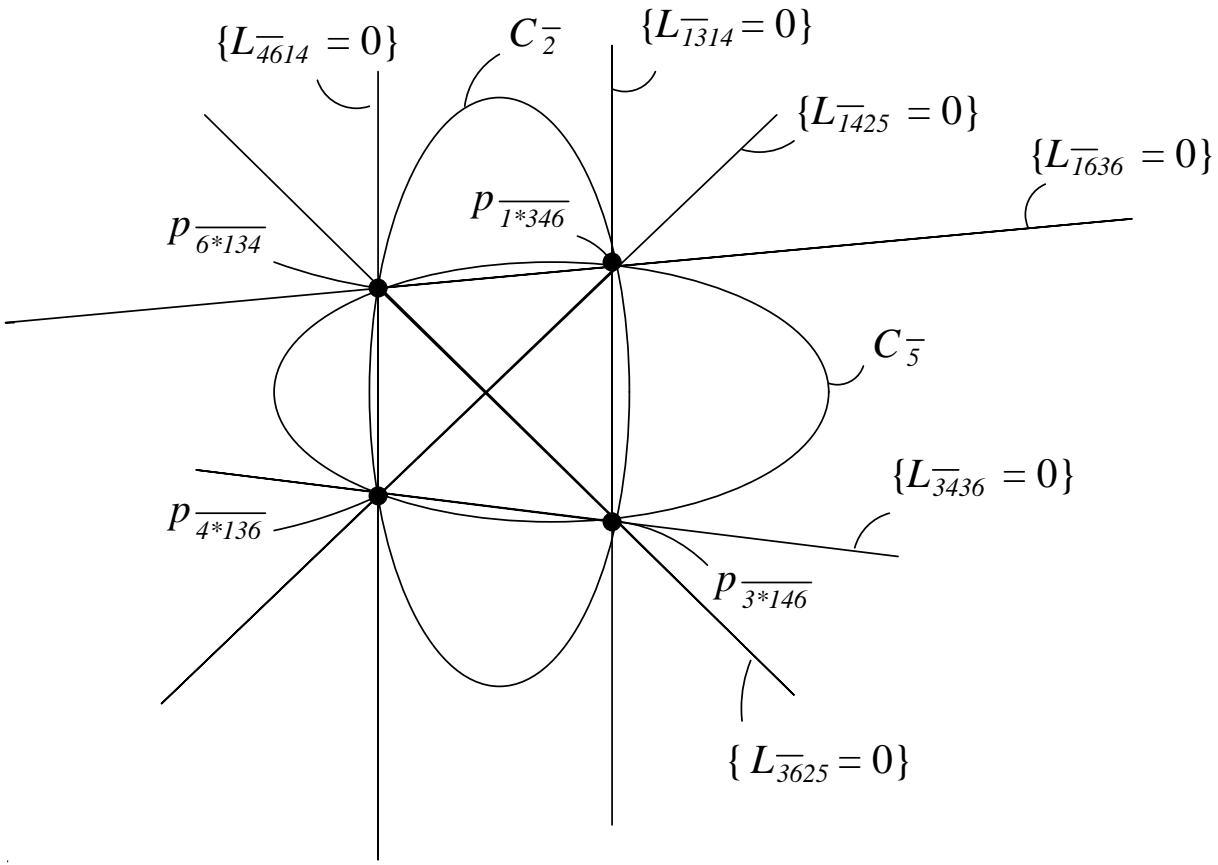}}

\caption{Intersecting conics within a system---a
tetrahedral configuration}
\label{fig:tetPic}

\end{figure}

Finally, which three involutions fix a \Bar{60},
60-point $p$?  (The stabilizer of $p$ is a
\D{3}.) Unlike the other special orbits,
\Orb{\Bar{60}}\ and
\Orb{60}\ have a bias toward one or the other system
of conics.  Recalling that $p$ is an icosahedral
face-center for two conics, say
\Cb{a}\ and \Cb{b}, an involution $Z$ that fixes $p$
cannot preserve the two conics individually; such an
action would have order three.   Hence, $Z$ exchanges
the conics and lacks a prefix \Bar{ab}. For
$\Bar{ab}=\Bar{25}$ the array indicates six such
involutions:
$$ \begin{array}{lll}
\ZZ{13}{14}:\; (\Bar{25})(\Bar{46})
&\hspace*{10pt}& \ZZ{46}{14}:\;
(\Bar{25})(\Bar{13})\\[10pt]
\ZZ{36}{25}:\; (\Bar{25})(\Bar{14})
&\hspace*{10pt}& \ZZ{14}{25}:\;
(\Bar{25})(\Bar{36})\\[10pt]
\ZZ{34}{36}:\; (\Bar{25})(\Bar{16})
&\hspace*{10pt}& \ZZ{16}{36}:\;
(\Bar{25})(\Bar{34}).
\end{array} $$
The six associated lines pass through the four points
in \mbox{$\Cb{2}
\cap \Cb{5}$} as edges of the tetrahedron whose
stabilizer is
$\Tet_{\Bar{25}}$.  They naturally fall into four sets
of three lines; among these triples a given index in
$\{\Bar{1},\Bar{3},\Bar{4},\Bar{6}\}$ appears twice. 
The intersection of the three lines occurs at the
three-fold
\Bar{60}-points\footnote{See Figure~\ref{fig:tetPic}.}
thereby suggesting appropriate names:
$$ \begin{array}{lllllll}
\LZ{13}{14}&\cap&\LZ{14}{25}&\cap&\LZ{16}{36}&=&\{\pTbar{1}{346}\}\\
\LZ{13}{14}&\cap&\LZ{34}{36}&\cap&\LZ{36}{25}&=&\{\pTbar{3}{146}\}\\
\LZ{14}{25}&\cap&\LZ{34}{36}&\cap&\LZ{46}{14}&=&\{\pTbar{4}{136}\}\\
\LZ{16}{36}&\cap&\LZ{36}{25}&\cap&\LZ{46}{15}&=&\{\pTbar{6}{134}\}.
\end{array} $$
Similarly, I will call the 60-points \pT{a}{bcd}.

To finish off the description of the special
Valentiner points on a 45-line:  Which four of the
\Bar{60}, 60-points lie on \LZ{ab}{cd}?  Since
\LZ{ab}{cd} contains 60-points whose indices satisfy
$c\!*\!dst$ and $d\!*\!cxy$, the matter comes down to
finding values of $s,t,x,y$ that are ``Valentiner
consistent."  This means that they fill out the scheme
$$ \begin{array}{lllllll}
\pT{c}{dst}&\in&
\LZ{\cdot\cdot}{cd}&\cap&\LZ{\cdot\cdot}{cs}&\cap&\LZ{\cdot\cdot}{ct}\\
\pT{c}{duv}&\in& 
\LZ{\cdot\cdot}{cd}&\cap&\LZ{\cdot\cdot}{cu}&\cap&\LZ{\cdot\cdot}{cv}\\
\pT{d}{cxy}&\in& 
\LZ{\cdot\cdot}{cd}&\cap&\LZ{\cdot\cdot}{dx}&\cap&\LZ{\cdot\cdot}{dy}\\
\pT{d}{czw}&\in& 
\LZ{\cdot\cdot}{cd}&\cap&\LZ{\cdot\cdot}{dz}&\cap&\LZ{\cdot\cdot}{dw}
\end{array} $$
where each triple of prefixes exhausts
$\{\Bar{1},\ldots ,\Bar{6}\}$ and 
$$ \{s,t,u,v\}\ =\ \{x,y,z,w\}\ =\ \{1, \ldots,
6\} - \{c,d\}. $$
For \LZ{12}{34},
$$ \begin{array}{lllllll}
\pT{3}{456}&\in& 
\LZ{12}{34}&\cap&\LZ{56}{35}&\cap&\LZ{34}{36}\\
\pT{3}{124}&\in& 
\LZ{12}{34}&\cap&\LZ{35}{13}&\cap&\LZ{46}{23}\\
\pT{4}{356}&\in& 
\LZ{12}{34}&\cap&\LZ{34}{45}&\cap&\LZ{56}{46}\\
\pT{4}{123}&\in& 
\LZ{12}{34}&\cap&\LZ{35}{24}&\cap&\LZ{46}{14}
\end{array} $$
and
$$ \begin{array}{lllllll}
\pTbar{1}{236}&\in& 
\LZ{12}{34}&\cap&\LZ{13}{26}&\cap&\LZ{16}{15}\\
\pTbar{1}{245}&\in& 
\LZ{12}{34}&\cap&\LZ{14}{25}&\cap&\LZ{15}{16}\\
\pTbar{2}{136}&\in& 
\LZ{12}{34}&\cap&\LZ{23}{15}&\cap&\LZ{26}{26}\\
\pTbar{2}{145}&\in& 
\LZ{12}{34}&\cap&\LZ{24}{16}&\cap&\LZ{25}{25}\,.
\end{array} $$

Collected below are the various data for \LZ{12}{34}.
$$ \begin{array}{c|c|c}
\mbox{Orbit}&\mbox{Special points on
\LZ{12}{34}}& \mbox{Multiplicity of 45-lines}\\
\hline
&&\\[-5pt]
\Orb{36}&\pD{1}{3}\ \; \pD{1}{4}\ \; \pD{2}{3}\
\; \pD{2}{4}&
 \left(\begin{array}{c}5\\2\end{array}\right) \ =\ 
10\\[15pt]

\Orb{45}&\pZ{12}{12}\ \; \pZ{12}{56}\ \;
\pZ{36}{34}\ \; \pZ{45}{34}&
 \left(\begin{array}{c}4\\2\end{array}\right) \ =\
6\\[15pt]

\Orb{\Bar{60}}&\pTbar{1}{236}\ \; \pTbar{2}{136}\ \;  
 \pTbar{1}{245}\ \; \pTbar{2}{145}&
 \left(\begin{array}{c}3\\2\end{array}\right) \ =\
3\\[15pt]

\Orb{60}&\pT{3}{124}\ \; \pT{4}{123}\ \;
\pT{3}{456}\ \; \pT{4}{356}&
 \left(\begin{array}{c}3\\2\end{array}\right) \ =\
3\\[15pt]

\Orb{90}&\pQ{12}{34}{1}\  \pQ{12}{34}{2}&0
\end{array} $$

For subsequent reference, the following table
summarizes the geometric terminology.
$$ \begin{array}{lcc}
 \mbox{point/line}&\mbox{orbit}&\mbox{stabilizer}\\

\cline{1-1}  \cline{2-2}  \cline{3-3} 

\pD{a}{b}/\LD{a}{b}&\Orb{36}&\D{5}\\

\{\pP{a}{b}{1},\pP{a}{b}{2}\}/
\{\LP{a}{b}{1},\LP{a}{b}{2}\}&\Orb{72}&\Z{5}\\

\pZ{ab}{cd}/\LD{ab}{cd}&\Orb{45}&\D{4}\\

\{\pQ{ab}{cd}{1},\pQ{ab}{cd}{2}\}/
\{\LQ{ab}{cd}{1},\LQ{ab}{cd}{2}\}&\Orb{90}&\Z{4}\\

\pTbar{a}{bcd}/\LTbar{a}{bcd}&\Orb{\Bar{60}}&\D{3}\\

\pT{a}{bcd}/\LT{a}{bcd}&\Orb{60}&\D{3}\\
\end{array} $$


\begin{figure}[hp]

\input{array45}

\caption{A combinatorial scheme for the Valentiner
group}

\label{fig:45Array}

\end{figure}

%% file: array45.tex
\begin{tabular}[ht]{c||c|c|c|c|c|c|c|c|c|c|c|c|c|c|c}

&12&13&14&15&16&23&24&25&26&34&35&36&45&46&56\\[5pt]
\hline
\hline

&&&&&&&&&&&&&&&\\[-12pt]

\Bar{12}&$\bullet$&&&&&&&&&$\bullet$&&&&&$\bullet$\\[5pt]
\hline
&&&&&&&&&&&&&&&\\[-12pt]

\Bar{13}&&&$\bullet$&&&&&&$\bullet$&&$\bullet$&&&&\\[5pt]
\hline
&&&&&&&&&&&&&&&\\[-12pt]

\Bar{14}&&$\bullet$&&&&&&$\bullet$&&&&&&$\bullet$&\\[5pt]
\hline
&&&&&&&&&&&&&&&\\[-12pt]

\Bar{15}&&&&&$\bullet$&$\bullet$&&&&&&&$\bullet$&&\\[5pt]
\hline
&&&&&&&&&&&&&&&\\[-12pt]

\Bar{16}&&&&$\bullet$&&&$\bullet$&&&&&$\bullet$&&&\\[5pt]
\hline
&&&&&&&&&&&&&&&\\[-12pt]

\Bar{23}&&&&$\bullet$&&$\bullet$&&&&&&&&$\bullet$&\\[5pt]
\hline
&&&&&&&&&&&&&&&\\[-12pt]

\Bar{24}&&&&&$\bullet$&&$\bullet$&&&&$\bullet$&&&&\\[5pt]
\hline
&&&&&&&&&&&&&&&\\[-12pt]

\Bar{25}&&&$\bullet$&&&&&$\bullet$&&&&$\bullet$&&&\\[5pt]
\hline
&&&&&&&&&&&&&&&\\[-12pt]

\Bar{26}&&$\bullet$&&&&&&&$\bullet$&&&&$\bullet$&&\\[5pt]
\hline
&&&&&&&&&&&&&&&\\[-12pt]

\Bar{34}&$\bullet$&&&&&&&&&&&$\bullet$&$\bullet$&&\\[5pt]

\hline
&&&&&&&&&&&&&&&\\[-12pt]

\Bar{35}&&$\bullet$&&&&&$\bullet$&&&&&&&&$\bullet$\\[5pt]
\hline
&&&&&&&&&&&&&&&\\[-12pt]

\Bar{36}&&&&&$\bullet$&&&$\bullet$&&$\bullet$&&&&&\\[5pt]
\hline
&&&&&&&&&&&&&&&\\[-12pt]

\Bar{45}&&&&$\bullet$&&&&&$\bullet$&$\bullet$&&&&&\\[5pt]
\hline
&&&&&&&&&&&&&&&\\[-12pt]

\Bar{46}&&&$\bullet$&&&$\bullet$&&&&&&&&&$\bullet$\\[5pt]
\hline
&&&&&&&&&&&&&&&\\[-12pt]

\Bar{56}&$\bullet$&&&&&&&&&&$\bullet$&&&$\bullet$&

\end{tabular}

%% file: sextic6.tex
\subsubsection{Computing a diagonal bub involution}

One way to approach the matter of an antiholomorphic
symmetry that exchanges systems of conics is to look for
points that such a symmetry should fix.  Given three such
points $a,b,c$ in coordinates $y$ where
$$ 
a=[1,0,0] \hspace*{15pt} b=[0,1,0]
\hspace*{15pt} c=[0,0,1], 
$$
the associated bub map would have the diagonal
form
$$  
\Bub[y_{1},\, y_{2},\, y_{3}]\ =\ 
 [\alpha\, \Bar{y_{1}},\, \beta\,
 \Bar{y_{2}},\,\Bar{y_{3}}]. 
$$
Fixing a fourth point determines appropriate values for
the inhomogeneous parameters $\alpha$, $\beta$. 

But, which points should such a map fix?  Moreover, how
many such anti-involutions should there be?

\begin{figure}[hp]

\resizebox{\textwidth}{!}{\includegraphics{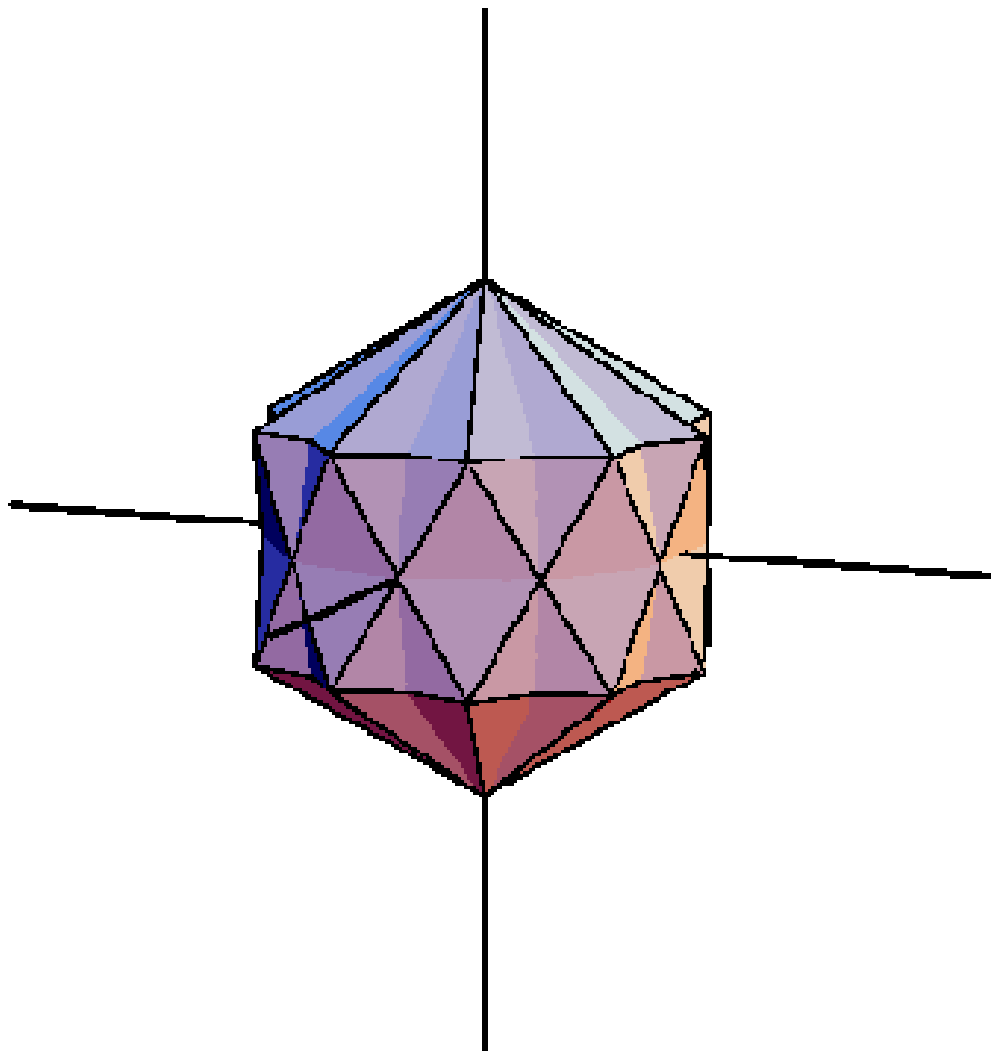}}

\caption{Intersecting icosahedra:  A heuristic for
bub-symmetry}
\label{fig:bubIcos}

The union of a barred and unbarred conic has a
\D{5}\ structure represented by two icosahedra that
``meet" at a pair of antipodal vertices and are turned
away from one another by an angle of $\pi/5$.  The
reflection through the equatorial plane exchanges the
icosahedra and so suggests that for each pair
$\{\Cb{a},\Cu{b}\}$ there is a ``primary" bub involution. 
Also transposing the icosahedra are ``secondary"
reflections through five vertical planes.  These
correspond to primary reflections for five other pairs of
conics.

\end{figure}

\paragraph{A heuristic for bub-symmetry.}
The basic Valentiner item that involves a mixture of
barred and unbarred conics is the \D{5}\ structure
consisting of a pair of conics
$\{\Cb{a},\Cu{b}\}$ that intersect in the pair of
72-points $\{\pP{a}{b}{1},\pP{a}{b}{2}\}$.  A bub map that
exchanges these two conics, must preserve the set
$\{\pP{a}{b}{1},\pP{a}{b}{2}\}$ as well as the associated
36-point \pD{a}{b}.  To put some flesh on the skeletal
configuration of Figure~\ref{fig:36a72PointsLines},
consider two icosahedra that 1) share a five-fold axis and
2) about this axis, are  one-tenth of a revolution away
from each other.\footnote{See Figure~\ref{fig:bubIcos}.}
The poles where the axis passes through the icosahedra
correspond to the pair of five-fold points
\pP{a}{b}{1,2}.  A reflection through the equatorial plane
preserves this arrangement while exchanging the icosahedra
and the poles.  The icosahedra also transpose under
reflection through five planes that include the polar
axis.  In these cases, the two poles are fixed.  This
model hints that for each pair \Cb{a}\ and \Cu{b}, there
is a distinguished bub-involution and five of a secondary
nature.  This makes for a total of 36 maps \bub{a}{b}.

\begin{figure}[hp]

\resizebox{\textwidth}{!}{\includegraphics{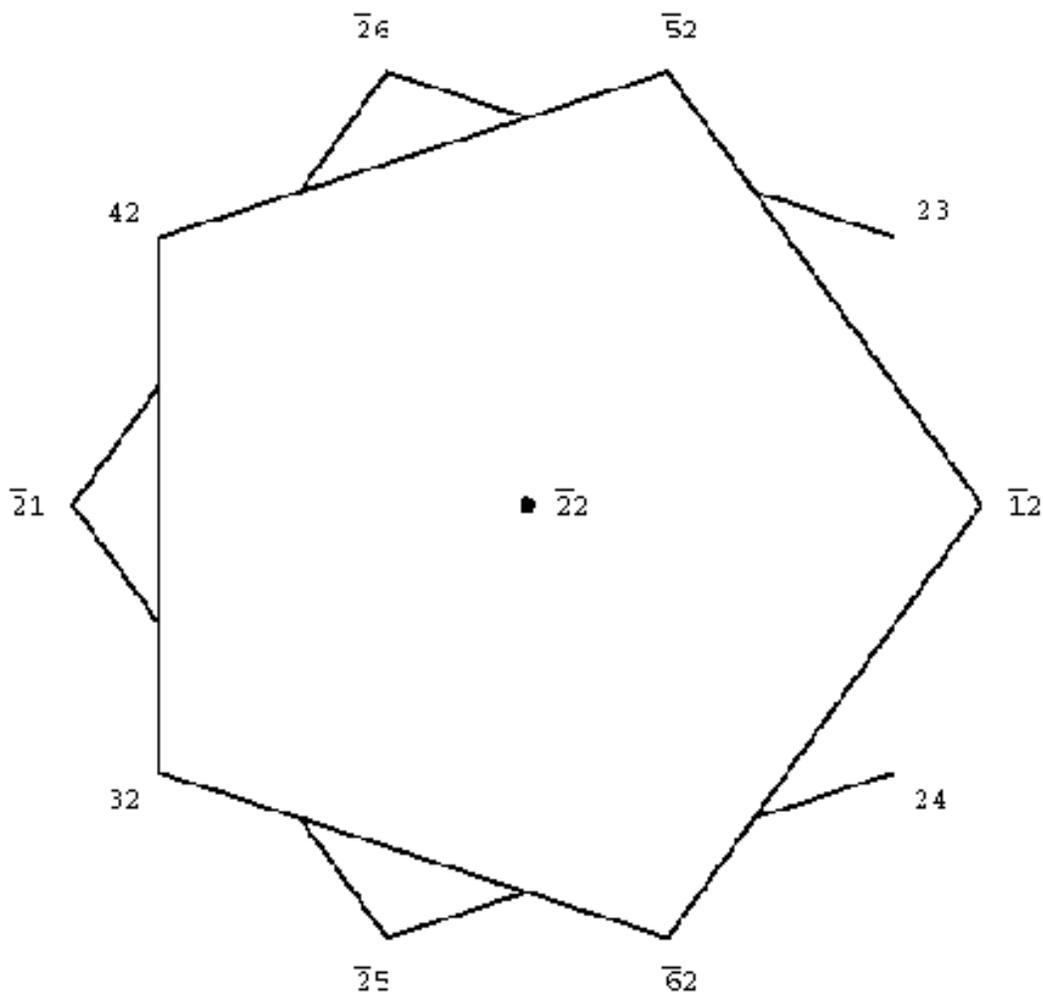}}

\caption{A heuristic for bub-symmetry}
\label{fig:bubPents}

Regard each of the five pairs of antipodal vertices on the
\D{5} union of conics \Cb{2}\ and \Cu{2}\ as a vertex of
one of two pentagons whose arrangement corresponds to that
of the remaining icosahedra. The primary
\bub{2}{2}\ reflection interchanges the pentagons as well
as antipodal vertices.  The secondary reflections are
$\{\bub{a}{a}|a\neq 2\}$ which transpose the vertices
$\Bar{2}a$ and $\Bar{a}{2}$ respectively.

\end{figure}

For the primary reflection relative to the pair
$\{\Cb{2},\Cu{2}\}$, this heuristic demands that
\pP{2}{2}{1}\ and
\pP{2}{2}{2}\ exchange while
\pP{2}{2}\ remains fixed.  What other points should
``\bub{2}{2}" fix?  Since five other bub maps switch
\Cb{2}\ and
\Cu{2}, symmetry requires that \bub{2}{2}\ provide a
secondary reflection for each barred-unbarred pair of
icosahedra associated with the \Bar{2}2 configuration.  As
such, \bub{2}{2}\ fixes the corresponding poles of
72-points.  Accordingly, the correspondence between the
five remaining barred and unbarred conics determines these
five pairs of points.  Now, the five pairs of non-polar
antipodal vertices on the \Cb{2}\ icosahedron correspond
to the points of intersection of \Cb{2}\ with the five
conics \Cu{1},
\Cu{3}, \Cu{4}, \Cu{5}, and \Cu{6} while the non-polar
vertices on the
\Cu{2} icosahedron correspond to the intersections of
\Cu{2}\ with the five conics \Cb{1}, \Cb{3}, \Cb{4},
\Cb{5}, and
\Cb{6}.  Let these sets of five pairs correspond to
vertices of two pentagons that are one-tenth of a
revolution away from each other\footnote{See
Figure~\ref{fig:bubPents} and imagine looking down, from
above a pole, on the intersecting icosahedra of
Figure~\ref{fig:bubIcos}.} with antipodal pairs of points
being
\bub{2}{2}\ symmetric.  The \D{5}\ action
$\DD{2}{2}=\St{\{\Cb{2},\Cu{2}\}}$ determines the specific
arrangement.

One of the elements of order five that belongs to
\DD{2}{2} is
$$ \PP{2}{2}\ =\ Q P Q^{-1}. $$ The associated
five-cyclings of conics are (\Bar{15436}) and (15436). The
matching of the $\Bar{2}k$ and
\Bar{m}2 vertices depends upon the five involutions that
stabilize \Cb{2} and
\Cu{2}.  For example, the generator \ZZ{12}{12}\
associates \Bar{1}\ with 1 while the 5-cycles above
determine the remaining matches of \Bar{a}\ with $a$.  This
information is also readily available in the 45-array. 
The five entries that involve both \Bar{2} and 2 are
$\Bar{2a}2a$.  Indeed, the array's symmetry about the
diagonal $\Bar{ab}ab$ is a combinatorial manifestation of
\bub{2}{2}.


\paragraph{Special coordinates and the bub-\RP{2}.}

Following the clue provided by the above heuristic, make
the parametrized change of icosahedral to octahedral
coordinates
$$ \begin{array}{lll}
A&=&
\left(\begin{array}{l|l|l}
a\,\pP{1}{1}{1}^t\;&\;b\,\pP{1}{1}{2}^t\;&\;\pD{1}{1}^t
\end{array} \right)\\[15pt] 

&=&
\left(\begin{array}{l|l|l}
a\left(\begin{array}{c}
 \frac{1 - \sqrt{5}}{2}\\[5pt] 
 -\sqrt{\frac{5 - \sqrt{5}}{2}}\, i\\[5pt] 1
\end{array} \right)
\:&\:
b\left(\begin{array}{c}
 \frac{1 - \sqrt{5}}{2}\\[5pt] 
 \sqrt{\frac{5 - \sqrt{5}}{2}}\, i\\[5pt] 1
\end{array} \right)
\:&\:
\left(\begin{array}{c}
 \frac{1 + \sqrt{5}}{2}\\[5pt] 0\\[5pt] 1
\end{array} \right)
\end{array} \right). 
\end{array} $$
In these icosahedral coordinates $y$, the \Bar{1}1
triangle is
$$ 
\pP{1}{1}{1} = [1,0,0] \hspace*{15pt}
\pP{1}{1}{2} = [0,1,0] \hspace*{15pt}
\pD{1}{1}    = [0,0,1] .
$$
The candidate \bub{2}{2} map $K(y)=\Bar{y}$ fixes each of
these. In octahedral coordinates the \Bar{2}2 triangle
consists of
$$ \begin{array}{lll}
\pP{2}{2}{1}&=& 
 \left[\frac{1 + \sqrt{5}}{2}\, \rho^{2},\
   \sqrt{\frac{5 + \sqrt{5}}{2}}\,i \rho,\
1\right]\\[15pt]
\pP{2}{2}{2}&=& 
 \left[\frac{1 + \sqrt{5}}{2}\, \rho^{2},\
   -\sqrt{\frac{5 + \sqrt{5}}{2}}\,i \rho,\
1\right]\\[15pt]
\pD{2}{2}&=&
 \left[\frac{1 - \sqrt{5}}{2}\, \rho^{2},\  0,\
1\right].
\end{array} $$
The hope here is that, when transformed to $y$
coordinates, some choice of $a,b$ results in
$$ 
K(A^{-1}\, \pP{2}{2}{1})\ =\ \pP{2}{2}{2}
\hspace*{20pt}
   K(A^{-1}\, \pP{2}{2}{2})\ =\ \pP{2}{2}{1}
\hspace*{20pt}
   K(A^{-1}\, \pD{2}{2}) = \pD{2}{2}. 
$$
Satisfying these conditions are the values 
$$ 
a=b=\frac{\sqrt{3}-\sqrt{5}\,i}{8}.
$$
The change of coordinates becomes
$$ \begin{array}{lll}
A &=&
\left( \begin{array}{ccc}
\frac{\left(1-\sqrt{5}\right)\,\left(3-\sqrt{5}\,i\right)}{4\sqrt{2}}&
\frac{\left(1-\sqrt{5}\right)\,\left(3-\sqrt{5}\,i\right)}{4\sqrt{2}}&
\frac{1+\sqrt{5}}{2}\\[5pt]
\frac{\sqrt{5-\sqrt{5}}\,\left(3-\sqrt{5}\,i\right)}{4}&
\frac{\sqrt{5-\sqrt{5}}\,\left(3-\sqrt{5}\,i\right)}{4}&0\\[5pt]
\frac{3-\sqrt{5}\,i}{2\sqrt{2}}&\frac{3-\sqrt{5}\,i}{2\sqrt{2}}&1
\end{array} \right).
\end{array} $$
As for the conic forms, they satisfy the desired
condition:
$$ \CCb{k}(y) = \Bar{\CCu{k}(\bar{y}}). $$

A further change of coordinates given by a real diagonal
matrix leaves
$\bub{2}{2}(y)=\Bar{y}$ undisturbed.  In $y$ coordinates
the one-point orbit
\pD{2}{2} under the \D{5} for \Bar{2}2 is\footnote{Recall
that $\tau=(1+\sqrt{5})/2$.}
$$ 
\left[1,1,\frac{\sqrt{6}}{\tau} \right].
$$
For a final simplification, the additional coordinate change
$$ 
B\ =\ \left( \begin{array}{ccc}
1&0&0\\ 0&1&0\\ 0&0&\frac{\sqrt{6}}{\tau}
\end{array} \right) 
$$
arranges for this point to be $[1,1,1]$. In these adjusted
$y$ coordinates, the \Bar{1}1 and \Bar{2}2 triangles
are\footnote{Recall that $\eta=(3+\sqrt{15}\,i)/4$.}
$$ \begin{array}{llllll}
\pP{1}{1}{1}&=&[1,0,0] \hspace*{20pt}& 
 \pP{2}{2}{1}&=&[3,2\,\eta^2,-\eta]\\[10pt]
\pP{1}{1}{2}&=&[0,1,0] \hspace*{20pt}& 
 \pP{2}{2}{2}&=&[3,2\,\Bar{\eta}^2,-\Bar{\eta}]\\[10pt]
\pD{1}{1}&=&[0,0,1] \hspace*{20pt}& 
 \pD{2}{2}&=&[1,1,1]
\end{array} $$
while, in these ``\bub{2}{2}\ coordinates" the normalized
conic forms\footnote{The unwieldy expressions for the
remaining forms are not recorded here.} for \Bar{1} and 1
are
$$ \begin{array}{lll}
\CCb{1}(y) &=&
\left(\frac{2}{3}\,\Bar{\eta}\right)^2\,y_1 y_2 +
y_3^2\\[10pt]
\CCu{1}(y) &=&
\left(\frac{2}{3}\,\eta\right)^2\,y_1 y_2 + y_3^2.
\end{array} $$

\begin{figure}[hp]

\resizebox{\textwidth}{!}{\includegraphics{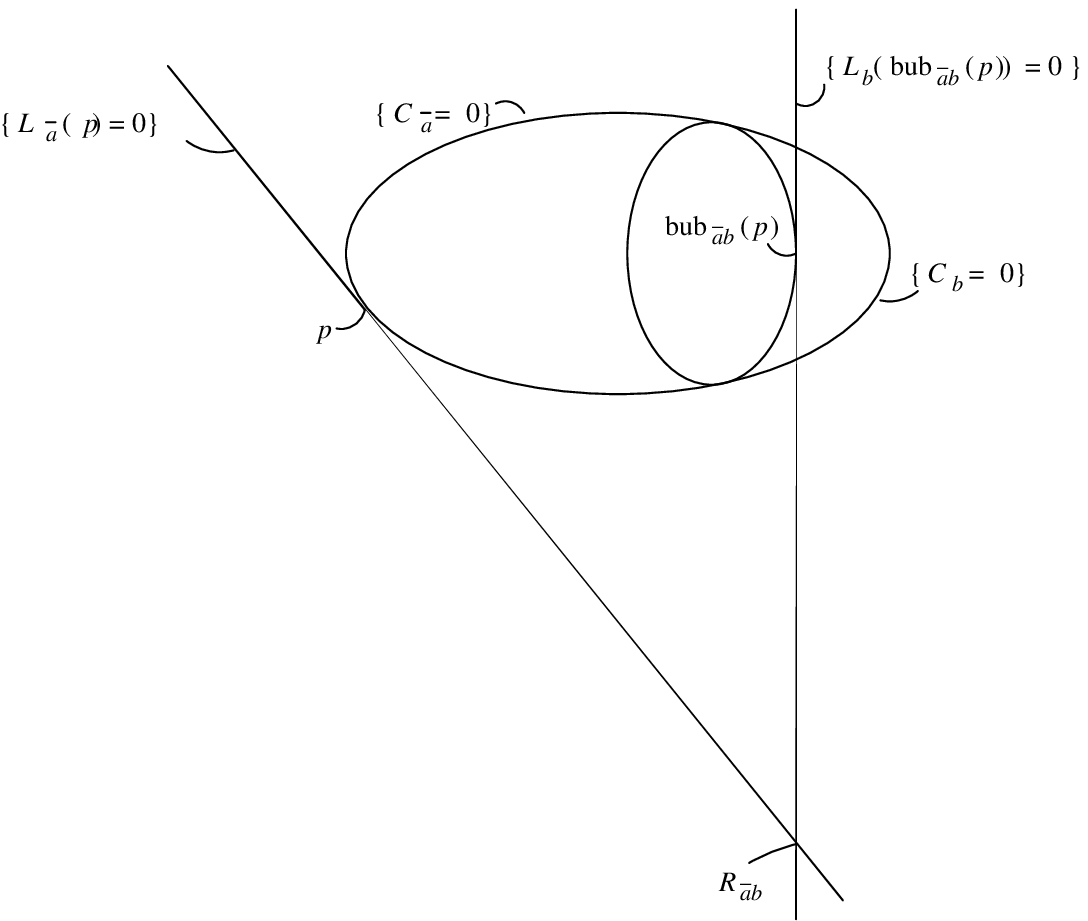}}

\caption{A geometric interpretation of bub-maps}
\label{fig:bubConstruct} Since \bub{a}{b} interchanges the
pair of conics \Cb{a}\ and \Cu{b}\ as well as the lines 
$$
\LL_{\Bar{a}}(p)=\{L_{\Bar{a}}(p)=0\} \hspace*{.3in} 
\LL_{b}(\bub{a}{b}(p))=\{L_{b}(p)=0\}
$$ 
that are tangent to \Cb{a}\ at $p$ and to \Cu{b} at
\bub{a}{b}(p), it also fixes the \RP{2}'s worth of points
$$ 
\RR{a}{b}=  \{\LL_{\Bar{a}}(p) \cap
\LL_{b}(\bub{a}{b}(p))|{p \in \Cb{a}}\}. 
$$

\end{figure}

Since \bub{2}{2}\ restricts to the identity on
the \RP{2}\ given by
$$ 
\RR{2}{2}\, =\, 
 \{[t_{1},t_{2},t_{3}]\,|\,t_{k}\, \in \R{} \}, 
$$
symmetry provides for such a fixed set \RR{a}{b}\ for each
of the 36 maps \bub{a}{b}.  Figure~\ref{fig:bubConstruct}
provides a geometric interpretation of these \RP{2}s.  One
consequence of the extra symmetry is that
$\Bar{\V}_{2\cdot 360}$-invariant forms and
$\Bar{\V}_{2\cdot 360}$-equivariant maps are, when
expressed in special bub-coordinates, given by polynomials
with real and even, in special cases, rational
coefficients.\footnote{See Section~\ref{sec:invariants}. 
Although \cite[pp. 548-50]{Wiman} and \cite[pp.
286-9]{Fricke} mention these coordinates, they seem not to
have made much use of them.}


%% file: sextic7.tex
\subsection{Invariant structure}

For a group action on a vector space the Molien series
provides one of the basic tools of classical invariant
theory.  Given a finite group \G\ acting faithfully on
\CC{n}, the dimension of the space $\CC{}[x]_m^{\G}$ of 
invariant homogeneous polynomials of degree
$m$ appears as the coefficient of the $m$th degree term in
the
\emph{Molien series} for \G:
$$
M(\CC{}[x]^{\G})\ =\ \sum_{m=0}^{\infty}\left(
\dim\CC{}[x]_m^{\G}\right) t^m. 
$$
In the Valentiner case the space is \CC{3} while the group
is a 1-to-3 lift of \V\ to a subgroup \Va\ of
$\mbox{SU}_{3}$.  As a result of the character
$\left<\rho\right>$ that appears under \V's action on the
icosahedral conic forms, this lift of \V\ to a
linear group has minimal order.\footnote{A lift of a
projective group \G\ to a linear group $\G^\prime$ has
\emph{minimal} order if for every lift
$\mathcal{H}$ of \G\, $|\mathcal{H}|\geq |\G^\prime|$.  See
\cite[pp. 267-8]{Fricke} for details.}  A further
consequence of minimality is that the Molien series for the
\Va\ gives complete information concerning the invariants
of the projective group \V.
\begin{prop} The invariants of \V\ and \Va\ are in
one-to-one correspondence. 
\end{prop}
\pf\ Trivially, a \Va-invariant gives a
\V-invariant.  Conversely, let $F(x)$ be a
\V-invariant with
$$F(T^{-1}x)=\alpha(T)\,F(x) $$
for $T\in \Va$.  The kernel of the multiplicative
character
$$\alpha:\Va \longrightarrow \CC{}-\{0\} $$
is the normal subgroup $\St(F) \subset \Va$ that
stabilizes $F$.  Since the projective image
$[\St(F)]\simeq \St(F)/\left<\rho\right>$ is
normal in the simple group
$\V\simeq \Va/\left<\rho\right>$,
$[\St(F)]$ is either trivial or \V.  In the
former case, $\St(F)$ would be either trivial or
$\left<\rho\right>$ so that $\Va/\St(F)$ would be
non-abelian.  Thus, $[\St(F)]=\V$. Since \Va\
includes no subgroup of order 360,
$\St(F)=\Va$ and $F$ is \Va-invariant. 
\fp
Also, \V\ lifts 1-to-6 to a so-called unitary reflection
group \Vb\ generated by 45 involutions on
\CC{3}.\cite[pp. 278, 287]{ST}  The elements of \Vb\
satisfy $\det\,T = \pm 1$ while those of \Va\ satisfy
$\det\,T=1$.

\subsubsection{Molien's theorem and its application to \V}
\label{sec:MolThm}

By projecting $\CC{}[x]_m$ onto $\CC{}[x]_m^{\G}$, one
arrives at a generating function for the Molien series.
\begin{thm} 
For a finite group action \G\ on \CC{n},
\begin{eqnarray*} M(\G) = \frac{1}{|\G|} \sum_{C_{T}
\subset \G}
 \frac{|C_{T}|}{\det\, (I\ -\ t\,T^{-1})}
 \end{eqnarray*} where $C_T$ are conjucay classes.
\end{thm}
\pf  See \cite[pp. 21-22]{Benson}. \fp
\begin{prop}
For the Valentiner groups \Va\ and \Vb, the
Molien series are given by
\begin{eqnarray*}
M(\Va) & = & \frac{1\, +\, t^{45}}{(1\, -\,
t^{6})\, (1\, -\, t^{12})\, (1\, -\,t^{30})}\\
&=& 1\, + t^{6} + 2\, t^{12} + 2\, t^{18} +
\ldots + t^{45} + \ldots\\[10pt]
M(\Vb) & = & \frac{1}{(1\, -\, t^{6})\, (1\, -\, t^{12})\, (1\, -\,t^{30})}\\
&=& 1 + t^{6} + 2\, t^{12} + 2\, t^{18} + \ldots.
\end{eqnarray*}
\end{prop}
\pf  With $k=0,1,2$, the matrices
$$
\pm\rho^k\,I,\ \pm\rho^k\,P,\ \pm\rho^k\,Z,\ \pm\rho^k\,Q,
\ \pm\rho^k\,P Z,\ \pm\rho^k\,T
$$
represent distinct conjugacy classes in \Vb.  For \Va\ the
three matrices of each type corresponding to the $+\rho^k$
do the job.  Substitution into the formula of Molien's
theorem produces the indicated generating functions.  
\fp

\subsubsection{The basic invariants themselves} 
\label{sec:invariants}

From the theory of complex reflection groups there are
three algebraically independent ``basic" forms that
generate the ring  of \Vb-invariants. The generating
function for the Molien series indicates that these occur
in degrees 6, 12, and 30.  Techniques of classical
invariant theory provide for the computation of the forms
in degrees 12 and 30 from that of degree 6.  But, how does
the latter arise? Although \Vb\ permutes the
\emph{conics}, its action on the conic forms is not
``simple"---a non-trivial character appears.  However, the
cubes of the forms do receive simple treatment by \Vb. 
Hence, summing the cubes of either system of conic forms and
normalizing the coefficients yields a
\Vb-invariant:
\begin{eqnarray*}
F(x)&=& \alpha\,\sum_{m=1}^{6}
\CCb{m}(x)^{3}\\[10pt]
&=& \alpha\,\sum_{m=1}^{6}
\CCu{m}(x)^{3}\\[10pt] 
&=& x_{1}^{6} +  x_{2}^{6} + x_{3}^{6}+ 
        3\,\left( 5 - \sqrt{15}\,i
\right)x_{1}^{2} x_{2}^{2} x_{3}^{2}\\[5pt] 
&&+\frac{3}{4}\,\left( 2\,\sqrt{5} -
     \left(5 - \sqrt{5} \right)\,\rho \right)\,
       \left( x_{1}^{4} x_{2}^{2} +
x_{2}^{4} x_{3}^{2} + x_{1}^{2} x_{3}^{4} 
       \right)\\[5pt] 
    &&-\frac{3}{4}\,\left( 2\,
       \sqrt{5} + \left(5 + \sqrt{5}
\right)\,\rho^{2} \right)\,
       \left( x_{1}^{4} x_{3}^{2} + x_{1}^{2}
x_{2}^{4} +x_{2}^{2} x_{3}^{4} 
       \right).  
\end{eqnarray*}
By uniqueness, $F$ is also $\Bar{\V}_{2\cdot
360}$-invariant.  Expressed and normalized in
\bub{2}{2} coordinates,
$$
F(y)=10\,y_{1}^{3} y_{2}^{3} + 9\,y_{1}^{5} y_{3}
 + 9\,y_{2}^{3} y_{3} - 45\,y_{1}^{2} y_{2}^{2}
y_{3}^{2}
 - 135\,y_{1} y_{2} y_{3}^{4} + 27\,y_{3}^{6}. 
$$

The forms $\Phi$ and $\Psi$ of degrees 12 and 30 arise
respectively from the determinants of the Hessian $H_{F}$ of
$F$ and the ``bordered Hessian" $BH(F,\Phi)$ of $F$ and
$\Phi$:
$$ \begin{array}{lll}

\Phi(y) 
&=& \alpha_{\Phi}\,|H(F(y))|\\[10pt]

&=& 6\,y_{1}^{11} y_{2} - 38\, y_{1}^{6}
y_{2}^{6} +
    6\,y_{1} y_{2}^{11} + 90\, y_{1}^{8}
y_{2}^{3} y_{3} +
    90\,y_{1}^{3} y_{2}^{8} y_{3} - 9\,y_{1}^{10}
y_{3}^{2} -\\[5pt]
&& 468\,y_{1}^{5} y_{2}^{5} y_{3}^{2} -  
   9\,y_{2}^{10} y_{3}^{2} + 
   1080\,y_{1}^{7} y_{2}^{2} y_{3}^{3} +
   1080\,y_{1}^{2} y_{2}^{7} y_{3}^{3} +\\[5pt]
&& 3375\,y_{1}^{4} y_{2}^{4} y_{3}^{4} - 324\,
y_{1}^{6} y_{2}
   y_{3}^{5} - 324\,y_{1} y_{2}^{6} y_{3}^{5} -
   1080\,y_{1}^{3} y_{2}^{3} y_{3}^{6} +
   2916\,y_{1}^{5} y_{3}^{7} +\\[5pt]
&& 2916\,y_{2}^{5} y_{3}^{7} + 1215\,y_{1}^{2}
y_{2}^{2} y_{3}^{8} +
    4374\,y_{1} y_{2} y_{3}^{10} +
729\,y_{3}^{12}\\[20pt]
\Psi(y) &=&\alpha_{\Phi}\,|BH(F(y),\Phi(y))|\\[10pt]

&=& \alpha_{\Phi}\,\left| \begin{array}{c|c}
H(\Phi(y))& \begin{array}{c}
F_{y_{1}}\\F_{y_{2}}\\F_{y_{3}} \end{array}\\
\hline
F_{y_{1}}\ F_{y_{2}}\ F_{y_{3}}&0
\end{array} \right|\\[30pt]

&=& 3\, y_1^{30} + \ldots\ + 3\, y_2^{30} +
\ldots + 57395628\,y_3^{30}.

\end{array} $$
The constants $\alpha_{\Phi}=-1/20250$ and
$\alpha_{\Psi}=1/24300$
remove the highest common factor among the
coefficients.  

Finally, the product of the 45 linear forms that correspond
to the generating involutions is a relative \Vb-invariant
but an absolute \Va-invariant and hence, a projective
\V-invariant.  As a specific instance of a general
result~\cite[p. 283]{ST}, this degree 45 form is given by
the Jacobian determinant
$$ \begin{array}{lll}

X(y)&=&
\alpha_X\,|J(F(y),\Phi(y),\Psi(y))|\\[10pt]
&=&\alpha_X\,\left| \begin{array}{ccc}
F_{y_{1}}&F_{y_{2}}&F_{y_{3}}\\[5pt]
\Phi_{y_{1}}&\Phi_{y_{2}}&\Phi_{y_{3}}\\[5pt]
\Psi_{y_{1}}&\Psi_{y_{2}}&\Psi_{y_{3}}
\end{array} \right| \\[30pt]
&=& \beta\, \prod \LLZ{ab}{cd}(y)\\[5pt]
&=& y_1^{45} + \ldots - y_2^{45} + \ldots +
3570467226624\,y_2^5 y_3^{40}

\end{array} $$
where $\alpha_X=-1/4860$ and $\beta$ is a
constant.  Being \Vb-invariant, $X^{2}$ is a
polynomial in $F$, $\Phi$, $\Psi$:
\begin{eqnarray} \label{eq:X^2}
 3^9\cdot X^{2}\ 
&=&4\,F^{13} \Phi + 80\,F^{11} \Phi^2 +
   816\,F^9 \Phi^3 + 4376\,F^7 \Phi^4 +
13084\,F^5 \Phi^5 +\\ \nonumber
&&12312\,F^3 \Phi^6 + 5616\,F \Phi^7 + 18\,F^{10}
\Psi + 198\,F^8
  \Phi \Psi + 954\,F^6 \Phi^2 \Psi -\\ \nonumber
&&198\,F^4 \Phi^3 \Psi - 5508\,F^2 \Phi^4 \Psi -
1944\,\Phi^5 \Psi
  - 162\,F^5 \Psi^2 - 1944\,F^3 \Phi \Psi^2 -\\  \nonumber
&& 1458\,F \Phi^2 \Psi^2 + 729 \Psi^3.  
\end{eqnarray}

\subsubsection{\V-symmetric maps and the sextic}

The system of invariants provides a foundation on which to
construct mappings of \CC{3}\ or \CP{2}\ that are symmetric
or \emph{equivariant} under the action of \Va\ or
\V.  Algebraically, this means that the map commutes with
the action.  Given such a map that also possesses nice
dynamical properties, the sixth-degree equation has an
iterative solution.

%% file: sextic8.tex
\section{Rational Maps with Valentiner Symmetry}

An iterative solution to the sextic utilizes a
parametrized family of dynamical systems having
\A{6} symmetry.  In practice, a given sixth-degree
polynomial with galois group \A{6} specifies a projective
transformation
$$S_{p}:\CP{2} \rightarrow \CP{2}$$ 
and thereby ``hooks-up" to a rational map 
$$S_p^{-1}\circ f\circ S_p$$
that has \A{6} symmetry.  Accordingly, the fixed map $f$ is
the centerpiece of a sextic-solving algorithm.

\subsection{Finding equivariant maps}

A linear group \G\ acts on the exterior algebra
$\Lambda(\CC{n})$ by
$$ (T(\alpha))(x) = \alpha(T^{-1}x) $$
where $T\in\,\G$, $\alpha$ is a $p$-form, and 
$x\in\CC{n}$.  As in the case of 
\V-invariant polynomials, \V-invariant $p$-forms associate
one-to-one with \Va-invariant $p$-forms.  Hence, the
search for symmetric maps can take place within the regime
of the linear action. 

When looking for equivariants, the \Va\ action on
$\Lambda(\CC{3})$ provides guidance.  Such utility is due
to a correspondence between \Va-equivariants and
\Va-invariant 2-forms.  Let $dX_2=(dx_2 \wedge dx_3,dx_3
\wedge dx_1,dx_1 \wedge dx_2 )$ and `$\cdot$' signify a
formal dot product.
\begin{prop}
For a given finite action $\G \subset
\mathrm{U}_3$ and a
\G-invariant 2-form 
\begin{eqnarray*}
 \phi(x)
&=&f_{1}(x)\, dx_2 \wedge dx_3 + f_{2}(x)\, dx_3 \wedge dx_1 + f_{3}(x)\, dx_1
\wedge dx_2 \\
&=&f(x)\cdot dX_2,
\end{eqnarray*}
the map \mbox{$f=(f_1,f_2,f_3)$} is relatively
$\G$-equivariant. If $\G \subset
\mathrm{SU}_3$, then the equivariance of $f$ is
absolute.
\end{prop}
\pf  Given $T \in \G$, let
$$ T = \left(
\begin{array}{ccc} 
t_{11}&t_{12}&t_{13}\\
t_{21}&t_{22}&t_{23}\\
t_{31}&t_{32}&t_{33}
\end{array}
\right) \hspace*{20pt}  
T^{-1} = |T|^{-1} \left(
\begin{array}{ccc} 
C_{11}&C_{12}&C_{13}\\
C_{21}&C_{22}&C_{23}\\
C_{31}&C_{32}&C_{33}
\end{array}
\right). $$
For an invariant 2-form $\phi$, 
\begin{eqnarray*} 
f(x)\cdot dX_2 
&=&\phi(x)\\
&=&T(\phi(x))\\
&=&f(T^{-1}x) \cdot |T|^{-2}(\\
&&\hspace*{-0pt}(C_{22}C_{33}-C_{23}C_{32})\,dx_2 \wedge dx_3 +
   (C_{23}C_{31}-C_{21}C_{33})\,dx_3 \wedge dx_1 +\\
&&\hspace*{10pt} (C_{21}C_{32}-C_{22}C_{31})\,dx_1 \wedge dx_2,\\
&&\hspace*{-0pt}(C_{32}C_{13}-C_{33}C_{12})\,dx_2 \wedge dx_3 +
   (C_{33}C_{11}-C_{31}C_{13})\,dx_3 \wedge dx_1 +\\
&&\hspace*{10pt}   (C_{31}C_{12}-C_{32}C_{11})\,dx_1 \wedge dx_2,\\
&&\hspace*{-0pt}(C_{12}C_{23}-C_{13}C_{22})\,dx_2 \wedge dx_3 +
   (C_{13}C_{21}-C_{11}C_{23})\,dx_3 \wedge dx_1 +\\
&&\hspace*{10pt}   (C_{11}C_{22}-C_{12}C_{21})\,dx_1 \wedge dx_2 )\\
&=&
f(T^{-1}x) \cdot \\
&&\hspace*{-50pt} |T|^{-2} \left( \begin{array}{ccc}
  C_{22}C_{33}-C_{23}C_{32}&
  C_{23}C_{31}-C_{21}C_{33}&
  C_{21}C_{32}-C_{22}C_{31}\\
  C_{32}C_{13}-C_{33}C_{12}&
  C_{33}C_{11}-C_{31}C_{13}&
  C_{31}C_{12}-C_{32}C_{11}\\
  C_{12}C_{23}-C_{13}C_{22}&
  C_{13}C_{21}-C_{11}C_{23}&
  C_{11}C_{22}-C_{12}C_{21}
\end{array} \right)
\left(\begin{array}{c}
dx_2 \wedge dx_3\\
dx_3 \wedge dx_1\\
dx_1 \wedge dx_2
\end{array} \right)\\ 
&=& f(T^{-1}x)^t\,|T|^{-2}\,|T|\,T^t\,dX_2^t\\
&=& |T|^{-1} (T f(T^{-1}x))^t\,dX_2^t\\
&=& |T|^{-1}\,T\,f(T^{-1}x)) \cdot dX_2.
\end{eqnarray*}           
Hence, 
$$f(x) = |T|^{-1}\,T\,f(T^{-1}x)$$ 
and $f$ is relatively equivariant. In case $T \in
\mathrm{SU}_3$, $|T|=1$ and absolute equivariance
occurs.  \fp 
Conversely, an absolute equivariant corresponds to a
relatively invariant 2-form, with absolute invariance
holding in case $\G \subset \mathrm{SU}_3$.

For invariant exterior forms, there is a 2-variable
``exterior Molien series" $M(\Lambda^{\G})$ in which the
variables $s$ and $t$ index respectively the rank of the
form and the polynomial degree:
$$ 
M(\Lambda^{\G}) = \sum_{p=0}^{n} \left(\,
  \sum_{m=0}^{\infty} 
\left(\dim\,\Lambda _p\left(\CC{n} \right)_m^\G \right)t^m 
 \right) s^p 
$$
where $\Lambda_p(\CC{n})_m^\G$ are the $\G$-invariant
homogeneous $p$-forms of degree $m$. (\cite[p. 62]{Benson}
or \cite[pp. 265ff]{Smith})  Projection of
$\Lambda_p(\CC{n})_m$ onto
$\Lambda_p(\CC{n})^\G_m$ yields the analogue to Molien's
theorem.
\begin{thm}
The exterior Molien series for a finite group action \G\
is given by the generating function
\begin{eqnarray*}
M(\Lambda^{\G})&=&
\frac{1}{\left|\G\right|}\  
 \sum_{\mathcal{C}_T \subset \G}
  \left|\mathcal{C}_T \right|\, 
   \frac{\det \left(I + s\,T^{-1}\right)}
   {\det \left(I - t\,T^{-1}\right)}
\end{eqnarray*}
where $\mathcal{C}_T$ are conjugacy classes. 
\end{thm}
As in the case of invariants, this result produces
\begin{prop} 
For \Va, the exterior Molien series is given by
$$\begin{array}{lll}
M(\Lambda^{\Va})&=&
\frac{1 + t^{45} + 
(t^5 + t^{11} + t^{20} +
t^{26} + t^{29} + t^{44})\,s +
(t + + t^{16} + t^{19} + t^{25} + t^{34} +
t^{40})\,s^2 + 
(1 + t^{45})\,s^3}
{(1 - t^6)(1 - t^{12})(1 - t^{30})}\\[10pt] 

&=&1 + {t^6} + 2\,{t^{12}} + 2\,{t^{18}} +
3\,{t^{24}} + 4\,{t^{30}} + \dots\\

&&+(t^5 + 2\,t^{11} + 3\,t^{17} + t^{20} + 4\,t^{23} + 
2\,t^{26} + 6\,t^{29} + \ldots)\,s\\

&&+(t + t^7 + 2\,t^{13} + t^{16} + 3\,t^{19} + t^{22} + 
5\,t^{25} + 2\,t^{28} + \ldots)\,s^2\\

&&+(1 + {t^6} + 2\,{t^{12}} + 2\,{t^{18}} +
3\,{t^{24}} + 4\,{t^{30}} + \dots)\,s^3.
\end{array}$$
\end{prop}

\subsection{A query on finite reflection groups}
 \label{sec:query}

For a reflection group \G\ that acts on \CC{n}\ the $n$
basic invariant 0-forms are algebraically independent.
\cite[pp. 282ff]{ST}  Multiplication of an invariant
$p$-form $\alpha$ of degree $\ell$ by an invariant 0-form
$F$ of degree $m$ promotes $\alpha$ to an invariant
$p$-form \mbox{$F\cdot \alpha$} of degree $\ell+m$.  In
the series $M(\Lambda^{\mathcal{H}})$ for a subgroup
$\mathcal{H}
\subset\:\G$ the contribution of the free algebra
generated by the basic 0-forms disappears upon division of
$M(\Lambda^{\mathcal{H}})$ by
$M(\Lambda_{0}^{\G})=M(\CC{}[x]^{\G})$.  The resulting
\emph{polynomial} in two variables displays the degrees of
the generating $\mathcal{H}$-invariant forms.  In the
cases of the 0, $n$-forms, which have identical series,
what remains are the terms corresponding to non-constant
polynomials that are
$\mathcal{H}$-invariant but not \G-invariant.
\begin{prop} For the Valentiner group, the ``exterior
Molien quotient" is
$$\begin{array}{ccl}

M(\Lambda^{\Va})/M(\Lambda_{0}^{\Vb})
&=&(1+t^{45})+(t^{5}+t^{11}+t^{20}+t^{26}+t^{29}+t^{44})s\,+\\

&&(t+t^{16}+t^{19}+t^{25}+t^{34}+t^{40})s^{2}+(1+t^{45})s^{3}.

\end{array}$$
\end{prop}
Notice the duality in degree 45 between 0 and 3-forms:
$$\begin{array}{ccc}

s^{0}: & 1=t^{0} & t^{45}\\[5pt] s^{3}: & t^{45}  & 1=t^{0}
   
\end{array}$$
and between 1 and 2-forms:

$$\begin{array}{ccccccc} 

s^{1}:&t^{5}&t^{11}&t^{20}&t^{26}&t^{29}&t^{44}\\[5pt]

s^{2}:&t^{40}&t^{34}&t^{25}&t^{19}&t^{16}&t.

\end{array}$$
By uniqueness, up to scalar multiplication, of the 3-form
$X\cdot\mathit{vol}$ associated with the 45 complex
planes of reflection, the exterior product of
``dual" forms must yield a multiple of this form.

This duality between invariant $p$ and $(n-p)$-forms also
appears in the \emph{ternary} icosahedral group
$\I_{60}$.  The generating 0-forms for the full reflection
group $\I_{2 \cdot 60}$ have degrees{\footnote{See line 23
of the table in \cite[p. 301]{ST}.} 2, 6, and 10. From the
discussion of the Valentiner conics the invariant of
degree 2 is familiar, while those of degrees 6 and 10 are
products of linear forms that $\I_{60}$ preserves.  Here,
the duality occurs in degree 15 which is the number of
reflection-planes for the icosahedron:
$$\begin{array}{ccl}

M(\Lambda^{\I_{60}})/ M(\Lambda_{0}^{\I_{2 \cdot 60}})
&=&(1+t^{15})+(t+t^{5}+t^{6}+t^{9}+t^{10}+t^{14})s\,+\\

&&(t+t^{5}+t^{6}+t^{9}+t^{10}+t^{14})s^{2}+(1+t^{15})s^{3}. 
 
\end{array}$$
Which finite reflection groups have series with this
property?  Is this duality connected to that described by
\cite[p. 286]{Orlik}?

\subsubsection{On \V-equivariance and special orbits}

Suppose $a$ is fixed by an element $T \in \G$.  Since a
\G-equivariant map
$f:\CP{n}\!\rightarrow\!\CP{n}$ satisfies 
$$f(a)=f(Ta)=Tf(a),$$ 
$T$ also fixes $f(a)$.  Hence, special orbits map to
special orbits. For the Valentiner action, the points
fixed by an involution are a 45-point and its line while
the 3, 4, 5-fold fixed points come in triples. 
Figure~\ref{fig:fix} summarizes the matter. Thus, under a
\V-equivariant $f$, a 45-line \LZ{ab}{cd}\ maps either to
itself or to its point \pZ{ab}{cd}.  In the former case,
$f$ preserves the pair of 90-points
$\{\pQ{ab}{cd}{1},\pQ{ab}{cd}{2}\}$.  Since \Orb{45}
cannot map to \Orb{90}, $f$ must fix the
45-points.\footnote{The possibility of $f$'s being
projectively undefined at \pZ{ab}{cd} is real. See below
for a case-study.}  Concerning the 36-72 triples the
matter stands just as in the case of the 45-90 points so
that $f$ either fixes or exchanges the 72-points
\pP{a}{b}{1}, \pP{a}{b}{2}\ and fixes the 36-points
\pD{a}{b}.  What about a triple of 60 points?  Symmetry
forces $f$ to permute the three points.  Since the
Valentiner group does not distinguish between 60-points, an
equivariant action must fix the orbit pointwise.

\subsection{The degree 16 map}

Returning to the exterior Molien series for \Va, the
coefficient of
$t^m s^2$ gives the dimension of the space of degree $m$
equivariants.  The series in $t$ begins
$$t+t^{7}+2\,t^{13}+t^{16}+\ldots.$$
The first term $t$ is due to the
identity\footnote{Although the identity \emph{map} is
always absolutely equivariant, its 2-form counterpart 
$$ x_1\,dx_2\wedge dx_3 + x_2\,dx_3\wedge dx_1 + x_3\,dx_1\wedge dx_2, $$
will not be absolutely invariant if the group is not in
$\mathrm{SU}_3$.} map while
$t^{7}$ occurs through promotion of the identity to the
degree 7 map
\mbox{$F\cdot
\id$}.  Two dimensions worth of invariants in degree 12
account for the $2\,t^{13}$ term.  The occurrence in
degree 16 of the first non-trivial equivariant finds
explanation in exterior algebra.  Since exterior
differentiation and multiplication preserve invariance,
the 2-form \mbox{$dF \wedge d\Phi$} is invariant and
hence, corresponds to an equivariant map whose coordinate
functions are given by the coefficients of the 2-form
basis 
$$ 
\{\mbox{$dx_2 \wedge dx_3$}, \mbox{$dx_3 \wedge dx_1$},
\mbox{$dx_1\!\wedge\!dx_2$}\}.
$$
\begin{prop}
Up to scalar multiplication, the unique \Va-invariant
2-form of degree 16 is 
$$ dF(x) \wedge d\Phi(x)=(\nabla F(x) \times
\nabla\Phi(x) ) \cdot dX_{2}. $$ Consequently, the unique
degree 16 \V-equivariant is 
$$ \psi_{16}(x)=\nabla F(x) \times \nabla\Phi(x) . $$
\end{prop}
Here, $\nabla$ is a ``formal" gradient $\nabla F =
(\frac{\partial F}{\partial x_1},\frac{\partial
F}{\partial x_2},\frac{\partial F}{\partial x_3})$ and
`$\times$' is the cross-product.  Geometrically,  $\psi$
associates with a point $y\!\in\!\CP{2}$ the intersection
of the pair of lines
$$ \psi(y)=\{\nabla F(y) \cdot x = 0\} \cap \{\nabla \Phi(y) \cdot x =0\}. $$


\begin{figure}[ht]    

$$ \begin{array}{l|l|l}
\mbox{Order}&\mbox{Number and Type}&\mbox{Notation}\\

\cline{1-1} \cline{2-2} \cline{3-3} 

\mbox{2-fold:}&\mbox{one 45-point, one 45-line} &
 \{\pZ{ab}{cd},\LZ{ab}{cd}\}\\[5pt]

\mbox{3-fold:}&\mbox{three 60-points}
 &\{\pT{a}{def},\pT{b}{def},\pT{c}{def}\}\\[5pt]

\mbox{3-fold:}&\mbox{three $\overline{60}$-points}
&\{\pTbar{ade}{f},\pTbar{b}{def},\pTbar{c}{def}\}\\[5pt]

\mbox{4-fold:}&\mbox{one 45-point, two 90-points} &
  \{\pZ{ab}{cd},\pQ{ab}{cd}{1},\pQ{ab}{cd}{2}\}\\[5pt]

\mbox{5-fold:}&\mbox{one 36-point, two 72-points} &
  \{\pP{a}{b},\pP{a}{b}{1},\pP{a}{b}{2}\}. 

\end{array}$$

\caption{Fixed points of \V}

\label{fig:fix}

\end{figure}


\subsubsection{Dynamics of the 16-map}

Given a generic point $a$ that lies on \emph{just one}
45-line \LZ{ab}{cd}\, the lines 
$$ \{\nabla F(a)\cdot x=0\}\ \mbox{and}\ \{\nabla \Phi(a)\cdot x=0\} $$
each pass through the point \pZ{ab}{cd}.  Thus, $\psi$
collapses \LZ{ab}{cd}\ to its companion point.  Taking
$\pZ{ab}{cd}=[1,0,0]$ and $\LZ{ab}{cd}=\{x_1=0\}$,
$$\psi(a)=[\psi_1(a),0,0]$$
where $\psi_1(a)$ is degree 16 in the homogeneous
coordinates $a=[0,a_{2},a_{3}]$.  The 16 roots of
$\psi_1(a)$ correspond to the 16 points where \LZ{ab}{cd}\
intersects the 44 remaining 45-lines.  These occur at the
36, 45, 60, \Bar{60}-points of which there are four each on
\LZ{ab}{cd}.  The ``blowing-down" of the 45-lines forces
the ``blowing-up" of their intersections:  To which
45-point does the intersection go?

Their collapsing behavior makes $\psi$ critical on the
45-lines.  Since the Jacobian determinant $|J_\psi|$ has
degree $3\cdot(16-1)=45$, $\{X=0\}$ is exactly the
critical set.  Thus, the 45-lines are superattracting. 
But, in approaching \LZ{ab}{cd}\ a trajectory inevitably
gets carried near the point \pZ{ab}{cd}\ which is blowing
up onto
\LZ{ab}{cd}; this means that the \CP{1}\ of directions
through
\pZ{ab}{cd}\ maps, by the agency of the Jacobian
transformation
$J_{\psi}$, to points on \LZ{ab}{cd}.\footnote{See
below.}  (Conversely, $J_{\psi}$ associates a point on
\LZ{ab}{cd} with a direction through
\pZ{ab}{cd}.) The next two iterations of $\psi$ first
return the trajectory to the vicinity of \pZ{ab}{cd} and
then send it back to
\LZ{ab}{cd} \emph{somewhere}.

Observation of trajectories that start near \pZ{ab}{cd}\
or \LZ{ab}{cd}\ reveals a rapid attraction on every other
iteration. What about the 3 and 5-fold intersections of
45-lines?  Are they attracting?  While these points do
indeed blow-up, their ``image" under $\psi$ is not a curve
that blows-down to the point;  such is the sole
propriety of the 45-points.  Hence, $\psi$ draws a generic
point near a 60, \Bar{60}, 36-point into one of three or
five point-line cycles.

This ``every-other" dynamics poses a problem:  which
every-other iterate do we watch?  For some initial
postions, the even iterates converge to a 45-point while
the trajectory spends the odd times around the
corresponding 45-line.  For others, the process is
reversed.  Moreover, experiment indicates that the
dynamics at the 45-line eventually settles down;
trajectories end up at one of the 45-points on the line. 
Hence, the iteration outputs a pair of 45-points each of
which lies on the line of the other.  But, there are, for
each 45-point, four possible pairs of this sort so that
taking every other iterate amounts to a neglect of
information.

The dynamics of $\psi$ appears to come down to what takes
place on the critical 45-lines. Given a point $x$ on a
45-line
\LZ{ab}{cd}, the derivative $J_{\psi}$ associates with $x$
a
\CP{1}\ through \pZ{ab}{cd}, namely, the image
$\mathcal{L}_{\psi}(x)$ of
$J_{\psi}(x)$.  In turn, $J_{\psi}(\pZ{ab}{cd})$ sends the
line $\mathcal{L}_{\psi}(x)$ to the point
$[J_{\psi}(\pZ{ab}{cd})]\mathcal{L}_{\psi}(x)$ on
$\LZ{ab}{cd}$.  The degree 15 map
$$ x \rightarrow [J_{\psi}(\pZ{ab}{cd})]\mathcal{L}_{\psi}(x) $$
gives $\psi$ on \LZ{ab}{cd}.  In
Appendix~\ref{sec:Seeing}, the dynamics on a 45-line
appears in a basins-of-attraction plot.\footnote{See
figure~\ref{fig:psi16Basin}.}  Here, the union of the
patches in a given color indicates a basin of attraction
for one of the four attracting 45-points.  The dark regions
correspond to points that the 45-points \emph{might} not
attract.  Does this set have interior or positive measure?

On a 45-line the map is not critically finite;  in
particular, its critical points are not
periodic.\footnote{On critical finiteness, see \cite[pp.
223ff.]{FS}.}  Indeed, there might be a wandering critical
point there.\footnote{Correspondence with Curt McMullen.}
Hence, establishing convergence almost everywhere would be
difficult should the map even possess this property. 
Moreover,  this behavior hardly reveals the geometric
elegance whose prospective discovery motivates the present
enterprise.

Unlike the 1-dimensional case of the icosahedral group for
which the non-trivial equivariant of lowest degree
provides an elegant dynamical system\footnote{See
\cite[pp. 152-3]{DM}.} for the purposes of solving the
quintic, the higher dimensional Valentiner action fails to
bear similar fruit.  The failure occurs in spite
of the 16-map's being obtained by a procedure analogous to
that employed by \cite{DM} in producing the degree-11
icosahedral map:
$$\begin{array}{ccc}

\underline{\mathrm{1\ dimension}}&\hspace*{20pt}&
\underline{\mathrm{2\ dimensions}}\\[5pt]

\mbox{$F$:  \I-invariant of degree 12}&&
\mbox{$F, \Phi$:  $\mathcal{V}$-invariants of degrees 6,
12}\\[5pt]

f_{11}=\left| \begin{array}{ll}
 \hat{x}&\hat{y}\\
  F_{x}&F_{y}
\end{array} \right|
&&
\psi_{16}=\left| \begin{array}{lll} 
\hat{x}&\hat{y}&\hat{z}\\
 F_{x}&F_{y}&F_{z}\\
\Phi_{x}&\Phi_{y}&\Phi_{z}
\end{array} \right|

\end{array}$$
where $\hat{x}, \hat{y} ,\hat{z}$ represent unit
coordinate vectors. However, in one dimension \emph{all}
\I-symmetric maps arise as combinations of three others
each of which are constructed in the manner of
$f_{11}$ but with different basic invariants standing in
for $F$.  A richer stock of equivariants inhabits the
Valentiner waters.  In contrast, the next higher degree
offers promise as well as a bit of mystery.

\subsection{A family of 19-maps}

The Molien series for Valentiner equivariants
$$ t+t^{7}+2\,t^{13}+t^{16}+3\,t^{19}\ldots $$
specifies three dimensions worth of maps in degree 19 of
which two are due to promotion of the identity by degree 18
invariants.  Hence, there are, as the exterior Molien
quotient\footnote{See Section~\ref{sec:query}.} indicates,
non-trivial \V-symmetric maps in degree 19.  How do these
arise?  Since there is no apparent exterior algebraic
means of producing such a map, the more practical matter
of computing them takes priority.\footnote{Peter Doyle and
Anne Shepler appear to have found a ``differential" method
of generating all invariant 2-forms.  The geometric and
dynamic consequences remain to be explored.}

\subsubsection{19 = 64 - 45}

Multiplication of a degree 19 equivariant $f$ by $X_{45}$
elevates $f$ to the 14 dimensional space of 64-maps. 
There are 14 ways of promoting the maps $\psi_{16}$,
$\phi_{34}$, and $f_{40}$ to degree 64:

\begin{description}

\item{1)} 7 dimensions of degree 48 invariants to promote $\psi_{16}=\nabla F \times \nabla
\Phi$

\item{2)} 4 dimensions of degree 30 invariants to promote $\phi_{34}=\nabla F \times \nabla
\Psi$
  
\item{3)} 3 dimensions of degree 24 invariants to promote $f_{40}=\nabla \Phi \times \nabla
\Psi$.

\end{description}
\begin{prop}
These 14 maps span the space of degree 64 equivariants.
\end{prop}
\pf If not, then for some $\alpha_k$ not all 0,
\begin{eqnarray*}
(\alpha_1\,F^8 + \ldots + \alpha_7\,F\,\Phi\,\Psi)\,dF \wedge d\Phi + 
(\alpha_8\,F^5 + \ldots + \alpha_{11}\,\Psi)\,dF \wedge d\Psi +&&\\
(\alpha_{12}\,F^4 + \alpha_{13} F^2\,\Phi + \alpha_{14}\,\Phi^2)\,d\Phi \wedge d\Psi&=&0.
\end{eqnarray*}
Since $dF \wedge d\Phi \wedge d\Psi =\beta\,X\cdot
\mathit{vol}$ for some constant $\beta\neq 0$, either
$$(\alpha_1\,F^8 + \ldots + \alpha_7\,F\,\Phi\,\Psi)\beta\,X = 0, $$
$$(\alpha_8\,F^5 + \ldots + \alpha_{11}\,\Psi)\beta\,X = 0, $$
or
$$(\alpha_{12}\,F^4 + \alpha_{13} F^2\,\Phi + \alpha_{14}\,\Phi^2)\beta\,X=0.$$
In at least one case, there is a non-zero $\alpha_k$. 
But, then the invariants $F$, $\Phi$, $\Psi$ are not
algebraically independent, contrary to the theory of
complex reflection groups. \cite[p. 282]{ST}  \fp In this
event, \mbox{$X\! \cdot \! f$} is a combination of maps
whose computation is straightforward.

Reasoning in the other direction, a 64-map
$$ f_{64} = F_{48}\cdot\psi_{16} + F_{30}\cdot\phi_{34} + F_{24}\cdot f_{40} $$
that ``vanishes" on the 45-lines---i.e., each coordinate
function of
$f_{64}$ vanishes---must have a factor of $X$.  The
quotient is a degree 19 equivariant
$$f_{19}=\frac{f_{64}}{X_{45}}.$$  Arranging for the vanishing of
$f_{64}$ on the 45-lines requires consideration of only
one line;  symmetry tends to the remaining 44.  Forcing
$f_{64}$ to vanish at 12 ``independent"
points\footnote{The points are \emph{independent} in the
sense that the 12 resulting linear conditions in the 14
undetermined coefficients of $f_{64}$ are independent.} on
a 45-line yields a 2-parameter family of 64-maps each
member of which vanishes on $\{X=0\}.$  The two
\emph{inhomogeneous} parameters reflect the three
dimensions (i.e., homogeneous parameters) of degree 19
\V-equivariants. In \bub{2}{2} coordinates, setting these
two parameters equal to 0 and normalizing the coefficients
yields
\begin{eqnarray}  \label{eq:f64} 
 f_{64}(y) 
&=&\left[10\,F(y)^6\,\Phi(y) + 100\,F(y)^4\,\Phi(y)^2 + 45\,F(y)^2\,\Phi(y)^3 +
         \right. \nonumber\\
&&\left. 156\,\Phi(y)^4 + 39\,F(y)^3\,\Psi(y) + 51\,F(y) \Phi(y)\,\Psi(y)\right]\cdot \psi(y) -        
\nonumber\\
&& 27\,\Psi(y) \cdot \phi(y) + 54\,\Phi(y)^2 \cdot f(y). 
\end{eqnarray}
The 2-parameter family of non-trivial 19-maps is then
\begin{eqnarray}  \label{eq:19-maps}
 g_{19}(y;a,b)\ =\ f_{19}(y)\:+\:  \left(a\,F(y)^{3} +
  b\,F(y)\,\Phi(y)\right) \cdot y.
\end{eqnarray}
Are any of these maps dynamically ``special"?  Indeed,
what might it mean to be special in this sense?

\subsubsection{Extended symmetry in degree 19}

Since $f_{19}$ is a non-trivial $\Bar{\V}_{2 \cdot
360}$-equivariant---note the integer coefficients
in (\ref{eq:f64}), each member of the two \emph{real}
parameter family
$$ f_{19}\ +\ F\left(a\,B_{12}+\Bar{a}\,U_{12}\right)\cdot \id  $$
is impartial towards the two systems of conics and so,
enjoys the additional symmetry. Here $B_{12} = \prod
C_{\overline{k}}$ and $U_{12} = \prod C_{k}$ are the
degree 12 invariants given by the product of the
respective six conic forms.  To honor the doubled symmetry
a member of this family must preserve each \RR{a}{b}\ that
\bub{a}{b}\ fixes point-wise.

\subsection{\emph{The} 19-map}

\subsubsection{The icosahedron again}

An intriguing aspect of the 19-maps is the degree itself.
Since 19 is one of the special equivariant
numbers\footnote{See \cite[p. 166]{DM}.} for the  binary
icosahedral group, there arises the prospect of finding a
\V-equivariant that restricts to self-mappings of the
conics.  By symmetry, an equivariant that
fixes\footnote{Here, `fix' means ``set-wise".} a conic of
a given system also fixes the other five.  Might there be a
map that preserves each of the 12 conics?

From \cite[p. 163]{DM} comes a geometric description of
the canonical degree 19 icosahedral mapping of the round
Riemann sphere:  stretch each face \textsf{F} over the 19
faces in the complement of the face antipodal to
\textsf{F} while making a half-turn in order to place the
three vertices and edges of \textsf{F} on the three
antipodal vertices and edges.   By symmetry, the 20
face-centers are fixed and repelling.  Since the resulting
map is critical only at the 12 period-2 vertices, the map
has nice dynamics. \cite[p. 156]{DM} A conic-fixing
\V-equivariant of degree 19 would, when restricted to a
conic \Cb{a} or \Cu{b}, give the unique map in degree 19
with binary icosahedral symmetry. This would determine its
effect on the special icosahedral orbits:

\begin{itemize}

\item[1)]  fix the face-centers:  \Bar{60}, 60-points of the
 appropriate system

\item[2)]  exchange antipodal vertices:  pairs of 72-points

\item[3)]  exchange antipodal edge-midpoints:  intersections of conics \Cb{a}, \Cu{b}\ with a
45-line indexed by \Bar{a} and $b$.

\end{itemize}

General \V-equivariants satisfy\footnote{While some maps
can be set to blow-up at the 60-points, such a
circumstance is rare.} condition 1).  Since
$$\Orb{72}=\{F=0\} \cap \{\Phi=0\},$$
the image of a 72-point under each $g_{19}$ in
(\ref{eq:19-maps}) is the same as that under $f_{19}$. 
Propitiously,  $f_{19}$ exchanges pairs of 72-points. 
Finally, a 19-map cannot blow-down a 45-line; this would
make it critical there---a condition precluded by the
critical set's being an invariant of degree
\mbox{$54\,=\,3\,(19-1)$}.  Consequently, each 45-line
must map to itself so that arranging for 3) costs one
parameter for each system of conics.  Fortunately, there
are two parameters to spend and their expenditure purchases
a canonical \V-equivariant $h_{19}$ that maps each of the
12 conics onto itself.\footnote{See below for an explicit
expression.}
\begin{fact}
There is a unique degree-19 \V-equivariant $h_{19}$ that
preserves each of the 12  icosahedral conics.
\end{fact}  
By favoring neither system of conics, $h_{19}$ possesses
\Bub-symmetry and so, self-maps each of the \RR{a}{b}. 
Expressing the family of 19-maps by
$$
g_{19}\ =\ h_{19}\:+\:F\left(a\,B_{12} +
b\,U_{12}\right) \cdot \id
$$  
makes evident the 1-parameter collections that fix the
barred $(b=0)$ and the unbarred $(a=0)$ conics.  

Unlike the 16-map $\psi_{16}$, $h_{19}$ does not blow up
somewhere.
\begin{prop}
The conic-preserving map $h_{19}$ is holomorphic
on \CP{2}.
\end{prop}
\pf\  By equivariance, the set of points on which
$h_{19}$ blows up is empty or a union of
\V-orbits. Direct calculation yields that 
$$h_{19}(p)\neq 0$$  
for 
$$
p\in \Orb{36}\cup \Orb{45}\cup \Orb{60}\cup
\Orb{\Bar{60}}\cup \Orb{90}.
$$
The remaining possibilities are that $h_{19}=0$ on a 180 or
360 point orbit.  

First, take the case of a 180 point orbit and
recall that each such point belongs to one
45-line.  Also, let
$$h_{19}=[h_1,h_2,h_3].$$   Since $h_{19}$
preserves each 45-line \LL, the only way that
$h_1=h_2=h_3=0$ is for the coordinates of the
restriction $h_{19}|_\LL$ to have a common
factor.  In \bub{2}{2} coordinates
$\LZ{12}{12}=\{y_1-y_2=0\}$ so that
$$h_{19}|_{\LZ{12}{12}}=[f,f,g].$$   But, the
resultant of $f$ and $g$ does not vanish.  Hence,
$f$ and $g$ do not have a common factor.

Finally, suppose that $h_{19}=0$ at a 360 point
orbit and that [0,0,1] is a 36-point $p_{36}$. 
Since 
$$|\{h_1=0\} \cap \{h_2=0\} \cap \{h_3=0\}| \leq
19\cdot 19=361,$$ there is only one member of
$h_{19}^{-1}(p_{36})$ in \CP{2}.  Of course, this
holds for every 36-point of which there are four
on a 45-line \LL.  Moreover,
$h_{19}(p_{36})=p_{36}$.  Thus, the
one-dimensional rational map $h_{19}|_\LL$ has
four exceptional points---a state of affairs that
requires the map to be degree one. \cite[p.
52]{Beardon}  Since the restricted map is not
degree one, $h_{19}$ does not vanish at a 360
point orbit.  (Indeed, no degree 19 equivariant
can blow up a 360 point orbit.)
\fp

\begin{fact}
The conic-fixing equivariant has the \bub{2}{2}\
expression:

\vspace{10pt}
\small
\input{h19}
\end{fact}
\normalsize

\subsubsection{Dynamical behavior}

The discovery of $h_{19}$ supplies the unique degree-19
\V-equivariant that self-maps, in addition to the 45-lines
and the 36 \Bub-\RP{2}s, the 12 conics.  The dynamics
\emph{on} each conic is well-understood.\footnote{A
basins-of-attraction plot appears in
Figure~\ref{fig:h19Conic}.}  
\begin{prop}
Under $h_{19}$, the trajectory of almost any point on an
icosahedral conic tends to an antipodal pair of the
superattracting vertices.
\end{prop}
Moreover, the conics themselves are attracting.
\begin{prop}
The Jacobian $J_{h_{19}}$ has rank one at the
superattracting 72-points. Thus, $h_{19}$ attracts on a
full \CP{2}\ neighborhood of such a point. Furthermore,
the Fatou components of the  restricted map
$h_{19}|_{\Cb{a}}$ are the intersections with \Cb{a} of
Fatou components of the map on \CP{2}.  
\end{prop}
Is this attracting behavior of the conics pervasive in the
measure-theoretic sense?  What about the ``restricted"
dynamics on the 45-lines and $\RP{2}$s?

Perhaps the place to begin is at a 72-point, say
\pP{1}{1}{1},  which lies at the hub of much Valentiner
activity.  Passing through
\pP{1}{1}{1} are many special objects. To enumerate:
 
\begin{itemize}
  
\item[1)]  the pair $\{\Cu{1}, \Cb{1}\}$ of conics, which meet tangentially
 
\item[2)]  the 36-line \LD{1}{1}, which gives a 10-fold
$\D{5}$ axis about which \mbox{ $\Cu{1} \cup \Cb{1}$}\
``turns"
  
\item[3)]  the 72-line \LP{1}{1}{2}, which is stable under
the cyclic half of the $\D{5}$ stabilizer of \LD{1}{1}\
and thereby tangent to
\Cb{1}\ and \Cu{1}
  
\item[4)]  the sixth-degree curve $\{F=0\}$, which is
tangent to \LD{1}{1}
 
\item[5)]  the twelfth-degree curve $\{\Phi=0\}$, which is
tangent to \LP{1}{1}{2}
 
\item[6)]  the five $\RP{2}$s $\{$\RR{2}{2}, \RR{3}{4},
\RR{4}{3}, \RR{5}{6}, \RR{6}{5}$\}$, each of which
intersect \Cu{1}\ and
 \Cb{1}\ only at \pP{1}{1}{1} and \pP{1}{1}{2}.  

\end{itemize}

In addition, a 72-point situates itself at the
intersection of two components $\{F_{6}=0\}$ and
$\{G_{48}=0\}$ of $h_{19}$'s critical set.  
\begin{fact}
The Jacobian determinant $|J_{h_{19}}|=F_{6}\,G_{48}$
where the invariant $G_{48}$ distinguishes itself by
lacking a $\Psi$ term when decomposed into an expression
in the basic invariants $F,\ \Phi,\ \Psi$:
\begin{eqnarray*}
 G_{48}(y)
&=&-13718\,\left[
14\,F(y)^8 + 180\,F(y)^6\,\Phi(y) + 
1701\,F(y)^4\,\Phi(y)^2 + \right.\\  
&&\left. 3402\,F(y)^2\,\Phi(y)^3 +
5103\,\Phi(y)^4 \right].
\end{eqnarray*}
\end{fact}
Since $G_{48}$ is a polynomial in $F$ and $\Phi$ alone,
its curve satisfies 
$$ \{F=0\} \cap \{G_{48}=0\}\ =\ \{F=0\} \cap \{\Phi=0\}\
=\ \Orb{72}. $$ Furthermore, the special invariant
structure of $G_{48}$ has an alternative expression in
terms of $B_{12}$ and $U_{12}$ alone.  From the identities
\begin{eqnarray*}
 B_{12}(y)&=&2\,(5-\sqrt{15}\,i)\,\rho\,F(y)^2 - 
  2\,\sqrt{15}(1+\sqrt{15}\,i)\,\rho^2\,\Phi(y)\\
 U_{12}(y)&=&2\,(5+\sqrt{15}\,i)\,\rho^2\,F(y)^2 - 
  2\,\sqrt{15}(1-\sqrt{15}\,i)\,\rho\,\Phi(y)
\end{eqnarray*}
it follows that
\begin{eqnarray*}
G_{48}(y)
&=&2^3 3^{24} 5^8 19^3 \left(-6\,(3 - \sqrt{15}\,i)\,\rho\,B_{12}(y)^4 \right. +\\
&&4\,(32 + 3\,\sqrt{15}\,i)\,\rho^2\,B_{12}(y)^3\,U_{12}(y)  - 333\,B_{12}(y)^2\,U_{12}(y)^2
+\\
&&\left.4\,(32 - 3\,\sqrt{15}\,i)\,\rho\,B_{12}(y)\,U_{12}(y)^3 -6\,(3 +
\sqrt{15}\,i)\,\rho^2\,U_{12}(y)^4 \right).
\end{eqnarray*}
Consequently, the degree 48 component of the critical set
meets the conics \emph{only} at the 72-points:
\begin{eqnarray*}
  \{G_{48}=0\} \cap \{B_{12}=0\} &=& \{G_{48}=0\} \cap \{U_{12}=0\}\\
                     &=& \{B_{12}=0\} \cap \{U_{12}=0\}\\
                     &=& \Orb{72}
\end{eqnarray*}
Accounting for the multiplicity at these eighth-order
intersections is the singularity of $\{G_{48}=0\}$ at
\Orb{72}, a result that follows directly from the
invariant decomposition.

\subsubsection{\RP{2}\ dynamics}

On each of the five \Bub-\RP{2}s that are mutually tangent
at
\pP{1}{1}{1} and \pP{1}{1}{2}, these 72-points are
superattracting for the restricted maps
$$ 
h_{19}|_{\RR{a}{b}},\ \ \ \Bar{a}b = \Bar{2}2,\,
\Bar{3}4,\,
 \Bar{4}3,\, \Bar{5}6,\, \Bar{6}5. 
$$
Are there attracting sites on \RR{a}{b} other than the five
pairs of 72-points?  Are there sets of positive measure or
open sets on which $h_{19}|_{\RR{a}{b}}$ fails to converge
to a pair of 72-points? The experimental evidence strongly
suggests that 
\begin{itemize}

\item[1)] the 72-points are the only attractors

\item[2)] there is no region with thickness or positive measure that remains outside of their influence

\item[3)] the set of 45-lines $\{X=0\}$ is repelling.

\end{itemize}
Appendix \ref{sec:Seeing} exhibits basin plots of
$h_{19}|_{\RR{2}{2}}$ which, of course, is dynamically
equivalent to each $h_{19}|_{\RR{a}{b}}$. What
significance does the
\RP{2}-dynamics hold for that on \CP{2}?  Extensive
trials\footnote{At \emph{http://math.ucsd.edu/\~{}scrass}
there are \emph{Mathematica} notebooks and supporting
files with which to iterate $h_{19}$ and, from the output,
to approximate a solution to a sixth-degree equation.} in
\CP{2}\ have not revealed behavior contrary to that
observed on the \RP{2}s.

\subsubsection{Conjectures}

\begin{conj}  \label{conj:h19} The only attracting
periodic points for $h_{19}$ are the elements of
\Orb{72}.  Moreover, the union of the basins of attraction
for \Orb{72}\ has full measure in  \CP{2}.
\end{conj}

The 45-lines present a problem in that they map to
themselves but do not contain the
72-points.\footnote{Again, this is a feature peculiar to
the 72-points.  They form the only special
\V-orbit that does not lie on the 45-lines.} The basin
plot in Figure~\ref{fig:h19RP2} reveals repelling behavior
along the \RP{1}\ where the \RP{2}\ meets one of the
45-lines that \bub{2}{2}\ fixes set-wise.  
\begin{conj} On the $45$-lines $h_{19}$ is repelling and,
hence, $\{X=0\}$ resides in the Julia set
$J_{h_{19}}$.
\end{conj} Is $J_{h_{19}}$ the closure of the backward
orbit of the 45-lines?

%% file: h19.tex
\noindent $h_{19}(y) = 1620\,F(y)^3 \cdot [y_1,y_2,y_3] + f_{19}(y) = $
\vspace{10pt}
\begin{sloppypar} \raggedright \noindent
$[-3591\,y_1^{15} y_2^4 - 
   5263\,y_1^{10} y_2^9 + 
   9747\,y_1^5 y_2^{14} - 81\,y_2^{19} + 
   17955\,y_1^{12} y_2^6 y_3 +  
   10260\,y_1^7 y_2^{11} y_3 - 
   7695\,y_1^2 y_2^{16} y_3 - 
   107730\,y_1^{14} y_2^3 y_3^2 - 
   74385\,y_1^9 y_2^8 y_3^2 +  
   161595\,y_1^4 y_2^{13} y_3^2 - 
   969570\,y_1^{11} y_2^5 y_3^3 + 
   1292760\,y_1^6 y_2^{10} y_3^3 - 
   46170\,y_1 y_2^{15} y_3^3 -   
   2346975\,y_1^8 y_2^7 y_3^4 - 
   807975\,y_1^3 y_2^{12} y_3^4 - 
   3587409\,y_1^{10} y_2^4 y_3^5 + 
   10277442\,y_1^5 y_2^9 y_3^5 +   
   13851\,y_2^{14} y_3^5 - 
   969570\,y_1^{12} y_2 y_3^6 - 
   3986010\,y_1^7 y_2^6 y_3^6 - 
   1939140\,y_1^2 y_2^{11} y_3^6 -  
   5263380\,y_1^9 y_2^3 y_3^7 - 
   28117530\,y_1^4 y_2^8 y_3^7 + 
   831060\,y_1^{11} y_3^8 + 
   2423925\,y_1^6 y_2^5 y_3^8 +   
   4363065\,y_1 y_2^{10} y_3^8 - 
   24931800\,y_1^8 y_2^2 y_3^9 + 
   43630650\,y_1^3 y_2^7 y_3^9 -
   31123197\,y_1^5 y_2^4 y_3^{10}+   
   9598743\,y_2^9 y_3^{10} + 
   14959080\,y_1^7 y_2 y_3^{11} +
   23269680\,y_1^4 y_2^3 y_3^{12} - 
   26178390\,y_1^6 y_3^{13}+   
   52356780\,y_1 y_2^5 y_3^{13} + 
   18698850\,y_1^3 y_2^2 y_3^{14} + 
   20194758\,y_2^4 y_3^{15} + 
   22438620\,y_1^2 y_2 y_3^{16}+   
   7479540\,y_1 y_3^{18},$\\[5pt] 
$-81\,y_1^{19} + 9747\,y_1^{14} y_2^5 - 
   5263\,y_1^9 y_2^{10} - 
   3591\,y_1^4 y_2^{15} -   
   7695\,y_1^{16} y_2^2 y_3 + 
   10260\,y_1^{11} y_2^7 y_3 + 
   17955\,y_1^6 y_2^{12} y_3 + 
   161595\,y_1^{13} y_2^4 y_3^2 -   
   74385\,y_1^8 y_2^9 y_3^2 - 
   107730\,y_1^3 y_2^{14} y_3^2 - 
   46170\,y_1^{15} y_2 y_3^3 +
   1292760\,y_1^{10} y_2^6 y_3^3-   
   969570\,y_1^5 y_2^{11} y_3^3 - 
   807975\,y_1^{12} y_2^3 y_3^4 -
   2346975\,y_1^7 y_2^8 y_3^4 + 
   13851\,y_1^{14} y_3^5+   
   10277442\,y_1^9 y_2^5 y_3^5 - 
   3587409\,y_1^4 y_2^{10} y_3^5 - 
   1939140\,y_1^{11} y_2^2 y_3^6 - 
   3986010\,y_1^6 y_2^7 y_3^6 -   
   969570\,y_1 y_2^{12} y_3^6 - 
   28117530\,y_1^8 y_2^4 y_3^7 - 
   5263380\,y_1^3 y_2^9 y_3^7 + 
   4363065\,y_1^{10} y_2 y_3^8+   
   2423925\,y_1^5 y_2^6 y_3^8 + 
   831060\,y_2^{11} y_3^8 + 
   43630650\,y_1^7 y_2^3 y_3^9 - 
   24931800\,y_1^2 y_2^8 y_3^9 +
   9598743\,y_1^9 y_3^{10} -
   31123197\,y_1^4 y_2^5 y_3^{10} + 
   14959080\,y_1 y_2^7 y_3^{11} + 
   23269680\,y_1^3 y_2^4 y_3^{12} +   
   52356780\,y_1^5 y_2 y_3^{13} - 
   26178390\,y_2^6 y_3^{13} + 
   18698850\,y_1^2 y_2^3 y_3^{14} + 
   20194758\,y_1^4 y_3^{15}+   
   22438620\,y_1 y_2^2 y_3^{16} + 
   7479540\,y_2 y_3^{18},$\\[5pt]
$-1026\,y_1^{17} y_2^2 - 
   3078\,y_1^{12} y_2^7 - 
   3078\,y_1^7 y_2^{12} - 
   1026\,y_1^2 y_2^{17}-   
   5130\,y_1^{14} y_2^4 y_3 + 
   113240\,y_1^9 y_2^9 y_3 - 
   5130\,y_1^4 y_2^{14} y_3 + 
   3078\,y_1^{16} y_2 y_3^2 -  
   272916\,y_1^{11} y_2^6 y_3^2 - 
   272916\,y_1^6 y_2^{11} y_3^2 + 
   3078\,y_1 y_2^{16} y_3^2 + 
   215460\,y_1^{13} y_2^3 y_3^3 +   
   687420\,y_1^8 y_2^8 y_3^3 + 
   215460\,y_1^3 y_2^{13} y_3^3 + 
   4617\,y_1^{15} y_3^4 + 
   937251\,y_1^{10} y_2^5 y_3^4 +   
   937251\,y_1^5 y_2^{10} y_3^4 + 
   4617\,y_2^{15} y_3^4 + 
   290871\,y_1^{12} y_2^2 y_3^5 + 
   4813992\,y_1^7 y_2^7 y_3^5 +   
   290871\,y_1^2 y_2^{12} y_3^5 - 
   1454355\,y_1^9 y_2^4 y_3^6 - 
   1454355\,y_1^4 y_2^9 y_3^6 + 
   2520882\,y_1^{11} y_2 y_3^7 +   
   8812314\,y_1^6 y_2^6 y_3^7 + 
   2520882\,y_1 y_2^{11} y_3^7 + 
   19876185\,y_1^8 y_2^3 y_3^8 + 
   19876185\,y_1^3 y_2^8 y_3^8 -  
   2036097\,y_1^{10} y_3^9 + 
   5623506\,y_1^5 y_2^5 y_3^9 - 
   2036097\,y_2^{10} y_3^9 + 
   5235678\,y_1^7 y_2^2 y_3^{10} +   
   5235678\,y_1^2 y_2^7 y_3^{10} + 
   37813230\,y_1^4 y_2^4 y_3^{11} - 
   2617839\,y_1^6 y_2 y_3^{12} - 
   2617839\,y_1 y_2^6 y_3^{12} -   
   2908710\,y_1^3 y_2^3 y_3^{13} + 
   6357609\,y_1^5 y_3^{14} + 
   6357609\,y_2^5 y_3^{14} - 
   5983632\,y_1^2 y_2^2 y_3^{15} - 
   4487724\,y_1 y_2 y_3^{17} - 
   1023516\,y_3^{19}]$.
\end{sloppypar}

%% file: sextic9.tex
\section{Solving the Sextic}  \label{sec:solve}

By means of various algebraic manipulations, a general
sixth-degree polynomial reduces to a member of a
2-parameter family of ``Valentiner resolvents".  Such a
reduction requires the extraction of square and cube
roots.   Furthermore, a certain set of sextics transforms
into a special 1-parameter collection of resolvents.   
These resolvents are especially suited for solution by an
iterative algorithm that exploits Valentiner symmetry and
symmetry-breaking.

\subsection{General sixth-degree Valentiner resolvents}

At the core of Klein's program for equation-solving is
the ``form problem" relative to a particular action of a
given equation's symmetry group:  for prescribed values
$a_1, \ldots, a_n$ of the generating invariants $F_1,
\ldots, F_n$ find a point $p$ common to the inverse
images $F_1^{-1}(a_1), \ldots, F_n^{-1}(a_n)$.  As with
the quintic, solving the general sextic is
tantamount to solving the corresponding form-problem.
(\cite[pp. 308-10]{Fricke} and \cite{Coble})
This circumstance has a ``projectively-equivalent"
formulation in terms of rational functions in the basic
invariants.  In the Valentiner setting this concerns
$$ 
Y_1 = \alpha\, \frac{\Phi}{F^2}\hspace*{25pt} Y_2 =
\beta\,\frac{\Psi}{F^5} 
$$  
where $\alpha$ and $\beta$
are chosen so that $Y_1=1$ and $Y_2=1$ at a
36-point.\footnote{In \bub{2}{2} coordinates
$\alpha = 1$ and $\beta = 1/4$.}  Given values $a_1$ and
$a_2$ of $Y_1$ and $Y_2$, the task is to find a point $z$
in \CP{2}\ that belongs to the \V-orbit
$$ Y_1^{-1}(a_1) \cap Y_2^{-1}(a_2). $$
Accordingly, the general 6-parameter sextic $p(x)$
reduces to a resolvent that depends on the two parameters
$Y_1$ and $Y_2$.  Such a reduction requires the
extraction of a cube root \cite[p. 285]{Fricke} in
addition to the square root of $p$'s discriminant.  This
cube root is a so-called ``accessory irrationality"---its
adjunction to the coefficient field does not reduce the
galois group.  The 1-to-3 correspondence between the
projective and linear Valentiner groups \V\ and \Va\
accounts for its appearance.  In the 1-dimensional
icosahedral case, the projective group lifts 1-to-2 to a
linear group thereby producing the need for an accessory
square root. \cite[pp. 172-3]{Klein}

As for the \emph{derivation} of a 2-parameter resolvent,
the map
$$ Y: \CP{2}-\{F(z)=0\} \rightarrow (\CP{2}-\{F(z)=0\})/\V $$
given by
$$ Y(z)=[F(z)^3\,\Phi(z),\Psi(z),F(z)^5]=[Y_1(z),Y_2(z),1] $$
provides the \V-quotient of $\CP{2}-\{F(z)=0\}$ in that
the fibers are \V-orbits.  The exceptional status of the
sixth-degree curve is due to its being the fiber above
the single point
$[0,1,0]$. Furthermore, under the icosahedral function
$$ \Ub{1}(z) = \frac{\CCb{1}(z)^3}{F(z)} $$
a fiber $Y^{-1}[a_1,a_2,1]$ maps to six points
$$ 
\left\{\left. \Ub{n}(z) =
\frac{\CCb{n}(z)^3}{F(z)}\,\right|\,
 n = 1, \ldots, 6 \right\}  
$$
where $z \in Y^{-1}[a_1,a_2,1]$.  The $\Ub{n}(z)$ are
the roots of the sixth-degree polynomial
$$R_z(u) = \prod_{n=1}^{6} (u\, -\, \Ub{n}(z)). $$
As $z$ varies in $\CP{2}-\{F(z)=0\}$, $R_z(u)$ yields a
family of sextic resolvents.  Since \Va\ permutes the
$\CCb{n}(z)^3$ simply---no multiplicative character
appears, $R_z$ is
\V-invariant in $z$ and hence, so is each
$u$-coefficient.  Expressing the coefficients in terms of
the basic invariants $F(z)$, $\Phi(z)$, $\Psi(z)$ and
then converting to $Y_1$ and $Y_2$ yields, in \bub{2}{2}
coordinates, the resolvents
$$ \begin{array}{lll}
R_Y(u) &=& R_{(Y_1,Y_2)}(u)\\[5pt]

&=&u^{6}\:+\:\frac{-5\,+\,\sqrt{15}\,i}{90}\, u^5\,+\\[5pt]

&&\frac{11\,(1\,-\,\sqrt{15}\,i)\,-\,
 3\,(3\,+\,\sqrt{15}\,i)\,Y_1}{2^2 3^5 5^2}\, u^4
  \:+\:\frac{(100\,+\,57\,\sqrt{15}\,i)\,+\,
   9\,(30\,+\,\sqrt{15}\,i)\,Y_1}{3^9 5^4}\, u^3\,+\\[5pt]

&&\frac{-(152\,+\,17\,\sqrt{15}\,i)\,+\,
 18\,(-21\,+\,4\,\sqrt{15}\,i)\,Y_1\,+\,
  27\,(-4\,+\,\sqrt{15}\,i)\,Y_1^2}{2^2 3^{11} 5^5}\,u^2\,+\\[5pt]

&&\frac{(425\,+\,103\,\sqrt{15}\,i)\,+\,
 6\,(75\,+\,193\,\sqrt{15}\,i)\,Y_1\,+\,
  27\,(-25\,+\,33\sqrt{15}\,i)\,Y_1^2\,-\,7776\,\sqrt{15}\,i\,Y_2}
   {2^3 3^{14} 5^8}\,u\,+\\[5pt]

&&\frac{-(5\,+\,3\,\sqrt{15}\,i)\,+\,
 9\,(15\,-\,7\,\sqrt{15}\,i)\,Y_1\,+\,
  81\,(25\,-\,\sqrt{15}\,i)\,Y_1^2\,+\,
   81\,(45\,+\,11\sqrt{15}\,i)\,Y_1^3}{2^4 3^{18} 5^8}\,.
\end{array} $$
This makes explicit the fact that the solution of $R_Y(u)$ follows from inversion of $Y$.  

For the unbarred functions
$$U_n(z) = \frac{\CCu{n}(z)^3}{F(z)}, $$
one obtains the associated resolvents $S_Y$ from $R_Y$ by
complex conjugation of the
$u$-coefficients:
$$ S_Y(u) = \Bar{R_Y(\Bar{u})}. $$

\subsection{Special sixth-degree resolvents}

For the resolvents $R_Y$ the parameter space is an affine
plane
$[Y_1,Y_2,1]$ that lifts to $\CP{2}-\{F(z)=0\}$.  There
is a complementary set of resolvents parametrized by a
\CP{1}\ that lifts to
$\{F(z)=0\}$.

The sixth-degree \V-invariant curve $\{F(z)=0\}$ is a
genus 10 surface that contains three special Valentiner
orbits:
$$ \begin{array}{lll}

\Orb{72}&=&\{F(z)=0\} \cap \{\Phi(z)=0\}\\[10pt]

\Orb{90}&=&\{F(z)=0\} \cap \{\Psi(z)=0\}\\[10pt]

\Orb{180}&=&\{F(z)=0\} \cap \{X(z)=0\} - \Orb{90}.

\end{array} $$
On $\{F(z)=0\}$ the rational map
$$ V(z) = \alpha\, \frac{\Phi(z)^5}{\Psi(z)^2} $$
gives\footnote{The constant $\alpha$ is chosen so that
$V=1$ at a 180-point. In \bub{2}{2} coordinates,
$\alpha = 3/8$.} the 2-4-5 quotient of $\{F(z)=0\}$ under
\V:
$$ V:\{F(z)=0\} \rightarrow \{F(z)=0\}/\V. $$
Furthermore, the icosahedral function
$$ S_{\Bar{1}}(z) = \frac{\Phi(z)^2}{\Psi(z)}\,\CCb{1}(z)^3 $$
divides $\{F(z)=0\}$ by the icosahedral subgroup \Ib{1}:
$$ S_{\Bar{1}}:\{F(z)=0\} \rightarrow \{F(z)=0\}/\Ib{1}. $$
A value $V_0 \neq 0, 1, \infty$ of $V$ has an inverse
image on $\{F(z)=0\}$ that consists of a \V-orbit of size
360 while the image of $V^{-1}(V_0) \cap \{F(z)=0\}$ under
$S_{\Bar{1}}$ is the set of six points
$$ 
\left\{\left. S_{\Bar{n}}(z) =
\frac{\Phi(z)^2\,\CCb{n}(z)^3}{\Psi(z)}\,
  \right|\, n = 1, \ldots, 6 \right\} 
$$
where $z \in V^{-1}(V_0) \cap \{F(z)=0\}$.  The
$S_{\Bar{n}}(z)$ supply roots of a \V-parametrized family
of sextic resolvents 
$$ T_z(s) = \prod_{n=1}^{6} (s\, -\, S_{\Bar{n}}(z)). $$
As above, $T_z$ is \V-invariant in $z$ and hence, so is
each $s$-coefficient. Expressing the coefficients in
terms of the basic invariants $F(z)$,
$\Phi(z)$, and $\Psi(z)$, restricting to $\{F(z)=0\}$,
and converting to $V$ gives the one-parameter resolvents
$$ \begin{array}{lcl}
T_V(s) &=& s^{6}\: -\:

 \frac{-3\,+\,\sqrt{15}\,i}{2^5 3^3 5^2}\,V s^{4}\: -\: 

 \frac{4\,+\,\sqrt{15}\,i}{2^8 3^6 5^5}\,V^{2} s^{2}\: +\:

 \frac{\sqrt{15}\,i}{2^6 3^7 5^8}\, V^{2} s\: +\:

 \frac{45\,-\,11\,\sqrt{15}\,i}{2^{13} 3^{11} 5^8}\,V^{3}.
\end{array} $$

Again, with the unbarred functions 
$$ S_n(z) = \frac{\Phi(z)^2}{\Psi(z)}\, \CCu{n}(z)^3, $$
conjugation of the coefficients of $T_V$ yields the
unbarred resolvents.

\subsection{Parametrized families of Valentiner groups}  
 \label{sec:ParamVals}

The algorithms that solve given resolvents $R_Y$ or $T_V$
employ an iteration of a dynamical system $\h{Y}(w)$ or
$\h{V}(w)$ that belongs to a family of maps
\begin{itemize}

\item[1)] parametrized by  $Y={(Y_1,Y_2)}$ or $V$

\item[2)] each member of which is conjugate to $h_{19}(y)$.

\end{itemize} 
The first task is to parametrize by $Y$ and $V$ families
of Valentiner groups.  Each such group supports a
conic-fixing 19-map the computation of which follows that
of its conjugate,
$h_{19}(y)$.

\subsubsection{Invariant building-blocks}

Success in finding a 19-map for almost every
value\footnote{The singular values $Y(\{X(z)=0\})$ of Y
and 0, 1, $\infty$ of $V$ provide exceptions.} of $Y$ or
$V$ requires provision only of basic invariant forms
$$
F_Y, \Phi_Y, \Psi_Y, X_Y\ \mbox{or}\ 
F_V, \Phi_V, \Psi_V, X_V
$$
that are parametrized by $Y$ and $V$.  In
turn, the latter three of each type depend on the single
forms $F_Y$ and $F_V$.  

Much of this development amounts to keeping track of
coordinates.  The Valentiner actions $\V_z$ and $\V_y$ on
the respective planes $\CP{2}_z$ and $\CP{2}_y$ are the
same---the parameter $z$ merely replaces $y$.  Think of
these as a parameter and reference space respectively.  

To obtain a parametrized sixth-degree form in the general
case:
\begin{itemize}
\item[1)] Compose $F(z)$ with a certain family of maps
$\tau_z(w)$ each of which is $\V_z$-equivariant and
linear in $w$.  (The $w$-space $\CP{2}_w$ is the
iteration space.)
  
\item[2)] Express the coefficients of the $w$ monomials
in terms of $F(z)$, $\Phi(z)$, and
$\Psi(z)$.

\item[3)] Convert these coefficients to expressions in
$Y_1$ and $Y_2$---as in the derivation of
$R_Y$.
 \end{itemize}
The special case requires more care.
\begin{itemize}
\item[1)] Compose $F(z)$ with a select family of maps
$\sigma_z(w)$ each of which is $\V_z$-equivariant and
linear in $w$.

\item[2)] Restrict to $\{F(z)=0\}$ and express the
coefficients of the
$w$ monomials in terms\footnote{The choice of
$\sigma_z(w)$ becomes significant at this stage;  see
below.} of $\Phi(z)$, $\Psi(z)$.

\item[3)] Divide through by any overall factors in $\Phi$
and $\Psi$ to obtain a polynomial whose degree in $z$ is
a multiple of 60---the degree of V---and then express the
result in terms of $V$.
\end{itemize}

\subsubsection{Sixth-degree forms in the 2-parameter
case.}  

Consider the family of maps
$$ y = \tau_z(w) = [F(z)^4\cdot z]\,w_1\,+\,[F(z)\cdot h_{19}(z)]\,w_2\, +\, k_{25}(z)\, w_3 $$
that are degree 25 in $z$ and projective transformations
in $w$.  Here, $k_{25}(z)$ is the equivariant whose
expression in ``Hermitian\footnote{The elements $T
\in \Va$ satisfy $T \Bar{T^{t}} = I$.} coordinates" is 
$$ k_{25}(z) = \Bar{\nabla F(\Bar{\nabla F(z)})}. $$

With $\tau_z$, one constructs a $z$-parametrized family
of Valentiner groups
$\V_w=\tau_z^{-1}V_y\tau_z$ each member of which acts on
$\CP{2}_w$.  By construction, this family possesses an
equivariance property: for $T \in \V_z, \V_y$,
$$ \tau_{Tz}(w) = T \tau_z(w). $$
Hence, 
$$ F(\tau_{Tz}(w)) = F(\tau_z(w)) $$
so that the $w$-coefficients of $F(\tau_z(w))$ are
$\V_z$-invariant and thereby expressible in terms of the
basic forms $F(z)$, $\Phi(z)$, $\Psi(z)$, and $X(z)$. 
However, since the degree \emph{in z} of each
$w$-coefficent is
$6\cdot25=150$, an odd power of $X(z)$ cannot appear in
the decomposition of these coefficients into polynomials
in the basic invariants.  Being \Va-invariant, $X^2$
decomposes\footnote{See (\ref{eq:X^2}).} into a
polynomial in $F$, $\Phi$, and $\Psi$.  Thus, each
coefficient is a combination of the forms of degrees 6,
12, and 30. After division by an appropriate power of
$F(z)$ as well as a simplifying numerical factor
$\alpha$, the result is expressible in terms of $Y_1$ and
$Y_2$:
\begin{equation} \label{eq:FInY} 
 F_Y(w) = \frac{F(\tau_z(w))}{\alpha\, F(z)^{25}}. 
\end{equation} The coefficients of the $Y$ monomials in a
$w$-coefficient are solutions to a system of linear
equations.  An expression for this fundamental form
appears in Appendix~\ref{sec:paramExp}.

An important matter concerns the degeneration of
$\tau_z(w)$ where the determinant $|\tau_z|$ vanishes. 
Taking $z^t$, $h_{19}(z)^t$, and
$k_{25}(z)^t$ to be column vectors
$$ \begin{array}{lll}
|\tau_z| & = &  

\left| \begin{array}{ccccc}
&|&&|&\\

F(z)^4\cdot z^t &|& F(z)\cdot h_{19}(z)^t &|&
k_{25}(z)^t\\

&|&&|&

\end{array} \right| \\[30pt]

&=& F(z)^5\:

\left| \begin{array}{ccccc}

&|&&|&\\

z^t &|&  h_{19}(z)^t &|& k_{25}(z)^t\\

&|&&|&
\end{array} \right| \\[30pt]

&=& -1458\,F(z)^5\,X(z).
\end{array} $$
The final equality follows by uniqueness of $X$ as a
degree 45 invariant and evaluation of $|\tau_z|$, $F(z)$,
and $X(z)$ at a single point.  Thus, the square of
$|\tau_z|$ is expressible in $F(z)$,
$\Phi(z)$, and $\Psi(z)$ alone.  In terms of $Y$,
\begin{eqnarray}
|\tau_z|^2 
&=& 1458^2\, F^{10} X^2 \nonumber\\
&=& 432\,F^{25}\,(Y_1 + 20\,Y_1^2 + 204\,Y_1^3
 + 1094\,Y_1^4 + 3271\,Y_1^5 + 3078\,Y_1^6 +\nonumber\\
&& 1404\,Y_1^7 + 18\,Y_2 + 198\,Y_1  Y_2 + 954\,Y_1^2  Y_2 -
 198\,Y_1^3  Y_2 - 5508\,Y_1^4  Y_2 +\nonumber\\
&& 1944\,Y_1^5  Y_2 - 648\,Y_2^2 - 7776\,Y_1  Y_2^2 - 
 5832\,Y_1^2  Y_2^2 + 11664\,Y_2^3).  \label{eq:tauDet}
\end{eqnarray}

\subsubsection{The special case.}

Take the family of maps
$$
\sigma_z(w) = 
[72\,\Phi(z)^4\cdot z]\,w_1 + 
[\Psi(z)\cdot h_{19}(z)]\,w_2 + 
[24\,\Phi(z)^2\cdot k_{25}(z)]\,w_3. 
$$
with $z$-degree 49 and $w$-degree one.  The integer
coefficients have been chosen so that, in
\bub{2}{2}-coordinates, the point
$$[w_1,w_2,w_3]=[1,1,1]$$
corresponds to the map 
$$\nabla F(z) \times \nabla X(z)$$
associated with the 2-form $dF\wedge dX$.  As in the
general situation,
$$ \sigma_{Tz}(w) = T \sigma_z(w) $$
so that
$$ F(\sigma_{Tz}(w)) = F(\sigma_z(w)). $$
Thus, the $w$-coefficients are $\V_z$-invariant and
thereby expressible in terms of the basic forms $F(z)$,
$\Phi(z)$, $\Psi(z)$, and $X(z)$.  Since the degree in
$z$ of $F(\sigma_z(w))$ is
$6\cdot49=294$, odd powers of $X(z)$ cannot take part in
the basic invariant decomposition of these coefficients. 
Furthermore, restriction of the parameter space
$\CP{2}_z$ to $\{F(z)=0\}$ yields coefficients in
$\Phi(z)$ and $\Psi(z)$ alone.  Finally, since
$294=22\cdot12+30=2\cdot12+9\cdot30$,
\mbox{$F(\sigma_z(w))|_{\{F(z)=0\}}$} is divisible by
\mbox{$\Phi(z)^{2}\,\Psi(z)$}.  Hence,
$$
F(\sigma_z(w))|_{\{F(z)=0\}} =
\eta\,\Phi(z)^{2} \Psi(z)^9 \F{V}(w) 
$$
where $\eta$ is a simplifying numerical factor and
$\F{V}(w)$ is a polynomial that is degree four in $V$ and
degree six in $w$.  The expression for $\F{V}(w)$ also
appears in Appendix~\ref{sec:paramExp}.

Since the parameter space gets restricted to
$\{F(z)=0\}$, $|\sigma_z|$ should not vanish there.  In
fact, 49 is the lowest degree in which this fails to
occur for three (projectively) distinct maps.  Explicitly,
$$ \begin{array}{lll}

|\sigma_z| &=&  

\left| \begin{array}{ccccc}

&|&&|&\\

72\,\Phi^{4}(z)\cdot z &|& \Psi(z)\cdot h_{19}(z) &|& 
 24\,\Phi^{2}(z)\cdot k_{25}(z)\\

&|&&|&

\end{array} \right| \\[30pt]

&=& 24\cdot 72\,\Phi^{4}(z)\,\Psi(z)\,\Phi^{2}(z)\:

\left| \begin{array}{ccccc}

&|&&|&\\

z &|&  h_{19}(z) &|& k_{25}(z)\\

&|&&|&

\end{array} \right| \\[30pt]

&=& 2^7 3^9\,\Phi(z)^6\,\Psi(z)\,X(z).

\end{array} $$
Furthermore, on {\{F(z)=0\}, the expression
(\ref{eq:X^2}) for $X^2$ reduces to
\begin{eqnarray*} 
X^2 
 &=& -\frac{1}{81}\, (8\, \Phi^5 \Psi - 3\,\Psi^3)\\[10pt] 
 &=& -\frac{\Psi^3}{81}\,(8\,\frac{\Phi^5}{\Psi^2} - 3)\\[10pt] 
 &=& -\frac{\Psi^3}{81}\,(8\,\frac{3\,V}{8} - 3)\\[10pt]
 &=& -\frac{\Psi^3}{27}\, (V-1).
 \end{eqnarray*} 
Consequently,
\begin{eqnarray*} 
 |\sigma_z|^2 &=&  -2^{14} 3^{15} \Phi^{12} \Psi^5 (V-1)\\
 &=& -2^{14} 3^{15} \Phi^2 \Psi^9\, \frac{\Phi^{10}}{\Psi^4}
     (V-1)\\
 &=& -2^{14} 3^{15} \Phi^2 \Psi^9\, \left(\frac{3\,V}{8}\right)^2 (V-1)\\
 &=& -2^8 3^{17} \Phi^2 \Psi^9 V^2 (V-1). 
\end{eqnarray*}

\subsubsection{The remaining basic invariants.} 

The forms of degrees 12, 30, and 45 arise from the
sixth-degree invariant as before.  However, a
parametrized change of coordinates requires special
handling.  Under $y=Ax$ the Hessian $H_x(F(x))$, Bordered
Hessian $BH_x(F(x),G(x))$, and Jacobian
$J_x(F(x),G(x),K(x))$ transform as
$$ \begin{array}{rll}

H_x(F(y)) &=& A^t H_y(F(y)) A\\[10pt]

BH_x(F(y),G(y)) &=& 

\left( \begin{array}{l|l} A^t&0\\ \hline 0&1 \end{array}
\right)

BH_y(F(y),G(y))

\left( \begin{array}{l|l} A&0\\ \hline 0&1 \end{array}
\right)\\[15pt]

J_x(F(y),G(y),K(y)) &=& J_y(F(y),G(y),K(y)) A

\end{array} $$
where the subscript indicates the differentiation
variable.  As for transformation of the respective
determinants:
$$ \begin{array}{rll}

|H_x(F(y))| &=& |A|^2\,|H_y(F(y))|\\[15pt]

|BH_x(F(y),G(y))| &=&  |A|^2\,|BH_y(F(y),G(y))|\\[15pt]

|J_x(F(y),G(y),K(y))| &=& |A|\,|J_y(F(y),G(y),K(y))|.

\end{array} $$
For the parametrized change of coordinates $y=\tau_z(w)$,
let\footnote{Recall the constants $\alpha_{\Phi},
\alpha_{\Psi}, \alpha_{X}$ from
Section~\ref{sec:invariants}.}
\begin{eqnarray*}
\Phi_Y(w)&=&\alpha_{\Phi}\,|H_w(F_Y(w))|\\[5pt]
\Psi_Y(w)&=&\alpha_{\Psi}\,|BH_w(F_Y(w),\Phi_Y(w))|\\[5pt]
X_Y(w)&=&\alpha_{X}\,|J_w(F_Y(w),\Phi_Y(w),\Psi_Y(w))|.
\end{eqnarray*}
Then
\begin{eqnarray} 
\Phi(y)
&=&\alpha_{\Phi}\,|H_y(F(y))| \nonumber \\
&=&\alpha_{\Phi}\,|H_y(F(\tau_z(w)))| \nonumber \\
&=&\alpha_{\Phi}\,|\tau_z|^{-2}\,|H_w(\alpha\,F(z)^{25}
F_Y(w))|
     \nonumber \\ 
&=&\frac{\alpha_{\Phi}}{|\tau_z|^2}\,(\alpha\,
    F(z)^{25})^3\,|H_w(F_Y(w))| \nonumber \\ 
&=&\frac{(\alpha\,F(z)^{25})^3}{|\tau_z|^2}\,\Phi_Y(w)
     \label{eq:PhiInY} \\[20pt]
\Psi(y)
&=&\alpha_{\Psi}\,|BH_y(F(y),\Phi(y))| \nonumber \\
&=&\alpha_{\Psi}\,|BH_y(F(\tau_z(w)),\Phi(\tau_z(w)))|
     \nonumber \\
&=&\alpha_{\Psi}\,|\tau_z|^{-2}\,|BH_w(\alpha\, F(z)^{25}
    F_Y(w),
\frac{(\alpha\,F(z)^{25})^3}{|\tau_z|^2}\,\Phi_Y(w))| 
      \nonumber\\ 
&=&\frac{(\alpha\,F(z)^{25})^8}{|\tau_z|^6}\,\Psi_Y(w)
     \label{eq:PsiInY} \\[20pt]
X(y)
&=&\alpha_{X}\,|J_y(F(y),\Phi(y),\Psi(y))| \nonumber \\
&=&\frac{(\alpha\,F(z)^{25})^{12}}{|\tau_z|^9}\,X_Y(w)
     \label{eq:ChiInY}.
\end{eqnarray}
With $y=\sigma_z(w)$, similar calculations lead to the one-parameter forms:
\begin{eqnarray*}
\Phi(y) &=&\frac{(\eta\,\Phi(z)^2\,\Psi(z)^9)^3}
                {|\sigma_z|^2}\,\Phi_V(w)\\[20pt]
\Psi(y) &=&\frac{(\eta\,\Phi(z)^2\,\Psi(z)^9)^8}
                {|\sigma_z|^6}\,\Psi_V(w)\\[20pt]
   X(y) &=& \frac{(\eta\,\Phi(z)^2\,\Psi(z)^9)^{12}}{|\sigma_z|^9}\,X_V(w).
\end{eqnarray*}

\subsubsection{The 19-maps}

With an invariant system in place for each (non-singular)
value of $Y$ (or $V$), production of the degree 19 map
that preserves all 12 of the conics proceeds as before:  

\begin{itemize}
\item[1)] determine a 64-map $(f_{64})_Y(w)$ that
vanishes at $\{X_Y(w)=0\}$

\item[2)] compute $(f_{19})_Y(w) = (f_{64})_Y(w)/X_Y(w)$

\item[3)] compute the conic-fixing map $h_Y(w)$.
\end{itemize}
In fact, the previous calculations in the coordinates
$\{y_1,y_2,y_3\}$ provide a framework for those at hand. 
Into the
$y$-expressions involving $F(y)$, $\Phi(y)$, $\Psi(y)$,
and $X(y)$ as well as the maps $\psi_{16}(y)$,
$\phi_{34}(y)$, and $f_{40}(y)$ substitute the
appropriate ones in $Y$ and $w$, namely, $F_Y(w)$,
$\Phi_Y(w)$, $\Psi_Y(w)$, $X_Y(w)$, $\psi_Y(w)$,
$\phi_Y(w)$, and
$f_Y(w)$.  The substitutions for the invariants appear
above.\footnote{See (\ref{eq:FInY}), (\ref{eq:PhiInY}),
(\ref{eq:PsiInY}), and (\ref{eq:ChiInY}).}  Concerning
maps, they transform as the coefficients of 2-forms.  In
terms of the cross-product, for $u=Ax$,
\begin{eqnarray*} 
\nabla_x(F(u)) \times \nabla_x(G(u)) 
&=&A^t \nabla_u F(u) \times A^t \nabla_u G(u)\\
&=&|A^t|\,((A^t)^{-1})^t\, 
[\nabla_u F(u) \times \nabla_u G(u)]\\  
&=&|A|\,A^{-1} 
\left(\nabla_u F(u) \times \nabla_u G(u)\right). 
\end{eqnarray*}
Accordingly,
\begin{eqnarray*}
\psi(y) 
&=&  \nabla_y F(y) \times \nabla_y \Phi(y)\\
&=&\frac{1}{|\tau_z|}\,\tau_z
    \left[\nabla_w\left(\alpha\,F(z)^{25} F_Y(w)\right)
     \times
      \nabla_w\left(
       \frac{(\alpha\,F(z)^{25})^3}{|\tau_z|^2}\,
        \Phi_Y(w)\right)\right]\\
&=&\frac{(\alpha\,F(z)^{25})^4}{|\tau_z|^3}\,
    \tau_z \nabla_w F_Y(w) \times \nabla_w\ \Phi_Y(w)\\
&=&\frac{(\alpha\,F(z)^{25})^4}{|\tau_z|^3}\, 
    \tau_z(\psi_Y(w)) \\[20pt]
\phi(y) &=&  \nabla_y F(y) \times \nabla_y \Psi(y)\\
&=&\frac{1}{|\tau_z|}\,\tau_z
  \left[\nabla_w\left(\alpha\,F(z)^{25} F_Y(w)\right) 
   \times 
    \nabla_w\left(\frac{(\alpha\,F(z)^{25})^8}{|\tau_z|^6}\,
     \Psi_Y(w)\right)\right]\\
&=&\frac{(\alpha\,F(z)^{25})^9}{|\tau_z|^7}\, 
  \tau_z(\phi_Y(w)) \\[20pt]
f(y) 
&=&  \nabla_y \Phi(y) \times \nabla_y \Psi(y)\\
&=&\frac{1}{|\tau_z|}\,\tau_z
  \left[\nabla_w\left(\frac{(\alpha\,F(z)^{25})^3}{|\tau_z|^2}\,
   \Phi_Y(w)\right) 
    \times 
     \nabla_w\left(\frac{(\alpha\,F(z)^{25})^8}{|\tau_z|^6}\,
      \Psi_Y(w)\right)\right]\\
&=&\frac{(\alpha\,F(z)^{25})^{11}}{|\tau_z|^9}\,\tau_z(f_Y(w)).
\end{eqnarray*}
The one-parameter maps transform as follows:
\begin{eqnarray*}
\psi(y)
&=&\frac{(\eta\,\Phi(z)^2\,\Psi(z)^9)^4}{|\sigma_z|^3}\,\sigma_z(\psi_Y(w))\\[20pt]
\phi(y)
&=&\frac{(\eta\,\Phi(z)^2\,\Psi(z)^9)^9}{|\sigma_z|^7}\,\sigma_z(\phi_Y(w))\\[20pt]
f(y)
&=&\frac{(\eta\,\Phi(z)^2\,\Psi(z)^9)^{11}}{|\sigma_z|^9}\,\sigma_z(f_Y(w)).
\end{eqnarray*}

Making the substitutions into the expression for
$f_{64}(y)$ yields a collection of maps
$(f_{64})_Y(w)$ each of which is divisible by $X_Y(w)$.
Finally, substitution into the formula for the canonical
map $h_{19}(y)$ supplies the conic-fixing family
$h_Y(w)$.  The explicit calculations appear in
Appendix~\ref{sec:paramExp}.  

\subsubsection{General expressions for the canonical
19-map.} 

Each $w$-coefficient in
$$ h_z(w)=\tau_z^{-1}(h_{19}(\tau_z(w))) $$
is a polynomial in $Y_1$ and
$Y_2$ with $z$-degree  
$$
\deg_z h_z(w) \leq 
21\cdot 25=19\cdot 25 + 2\cdot 25 = 12\cdot 40 +
45=30\cdot 16+45.
$$   
The factor of $|\tau_z|^{-1}$ due to $\tau_z^{-1}$ does
not affect the map on \CP{2}\ and so, is neglected. 
After dividing away a factor of
$X(z)$ the polynomial map $h_Y(w)$ satisfies 
$$ 
\deg_{Y_1}h_Y(w) \leq 40  \hspace*{20pt}
\deg_{Y_2}h_Y(w) \leq 16. $$
Hence, one finds the coefficients of the $w$ monomials by
solving, for each term, 353 linear equations.  

In the 1-parameter case, 
$$ h_z(w)=\sigma_z^{-1}(h_{19}(\sigma_z(w))) $$
and  
$$\deg_z h_z(w) \leq 21\cdot49 = 12\cdot 82 + 45 = 
12\cdot 2 + 30\cdot 32 + 45.
$$  
Thus, on
$\{F(z)=0\}$, $h_z(w)$ is divisible by
\mbox{$\Phi(z)^2\,X(z)$} so that
$$
\deg_z \left(\frac{h_z(w)}{\Phi(z)^2\,X(z)}\right) =
60\cdot 16 \hspace*{20pt}
\deg_V h_V(w) \leq 16. 
$$

\subsection{Symmetry lost, a root found}

Under Conjecture~\ref{conj:h19}, the trajectory
$\{h_{19}^k(y_0)\}$ converges to a pair of 72-points for
almost any $y_0 \in \CP{2}_y$. Being conjugate to
$h_{19}(y)$, the maps $h_Y(w)$ share this property for
points in $\CP{2}_w$. Breaking the \A{6}\ symmetry of
$R_Y(u)=0$ qualifies $h_Y(w)$ for a role in
root-finding.\footnote{An analogous treatment applies in
the one-parameter case.}

\subsubsection{Root-selection}

Consider the rational function
$$ 
\Bar{J}_z(w) =
\frac{\Bar{\Gamma}_z(w)^3}{F(z)\,\Psi(\tau_z(w))} 
$$ 
where
$$ 
\Bar{\Gamma}_z(w) = \sum_{m=1}^6 \prod_{n \neq m}
\CCb{n}(\tau_z(w))\cdot\CCb{m}(z). 
$$ 
At a 72-point pair
$\{q_1,q_2\}=\tau_z^{-1}(\{\pP{a}{b}{1},\pP{a}{b}{2}\})$
in $w$ space, five of the six terms in
$\Bar{\Gamma}_z(w)$ vanish.  The result is a ``selection"
of one of the six roots $\Ub{a}(z)$ of $R_Y(u)$:
\begin{eqnarray*} 
 \Bar{J}_z(\{q_1,q_2\})&=& 
  \frac{[\prod_{n \neq a} \CCb{n}(\{\pP{a}{b}{1},\pP{a}{b}{2}\})]^3}
       {\Psi(\{\pP{a}{b}{1},\pP{a}{b}{2}\})}\,\frac{\CCb{a}(z)^3}{F(z)}\\
&=&\frac{\Ub{a}(z)}{\mu}\,.
\end{eqnarray*}
Here, $\mu$ is the value of  
$$ 
\frac{\Psi(\{\pP{a}{b}{1},\pP{a}{b}{2}\})}
     {[\prod_{n \neq a}
\CCb{n}(\{\pP{a}{b}{1},\pP{a}{b}{2}\})]^3}\ 
$$ 
which is constant on \Orb{72}:
$$ \mu = \frac{6561\,(279 + 145 \sqrt{15}\,i)}{2}. $$

In light of the $\V_z$-invariance of $\Bar{\Gamma}_z(w)$,
$\Bar{J}_z$ also enjoys this property and so, presents a
second face which expresses itself in $Y$ and $w$.  Since
$$
\deg_z \Bar{\Gamma}_z(w) = 252 = 10\cdot 25 + 2 
= 6\cdot 42,
$$
it transforms by
$$ 
\Bar{\Gamma}_Y(w) =
\frac{\Bar{\Gamma}_z(w)}{\beta\,F(z)^{42}} 
$$ 
where $\beta$ is a simplifying factor for the coefficients
over $w$.  Each coefficient is a polynomial in $Y_1$ and
$Y_2$ whose coefficients satisfy a system of linear
equations.\footnote{The polynomial
$\Bar{\Gamma}_Y(w)$ is rather large and  takes an
explicit albeit partial form in
Appendix~\ref{sec:paramExp}.}   Now, let $T_Y$ be the
polynomial\footnote{See (\ref{eq:tauDet}).} in
$Y_1$ and $Y_2$ that satisfies
$$ |\tau_z|^2 = F(z)^{25}\,T_Y. $$
Then, from (\ref{eq:PsiInY}),
\begin{eqnarray*}
\Psi(\tau_z(w))&=&\frac{(\alpha\,F(z)^{25})^8}
                      {(F(z)^{25} T_Y)^3}\,\Psi_Y(w)\\[5pt]
&=&\frac{\alpha^8\, F(z)^{125}}{T_Y^3}\,\Psi_Y(w).
\end{eqnarray*}
Finally,
$$ \Bar{J}_Y(w) = \frac{\beta^3}{\alpha^8}\,
            \frac{(T_Y\,\Bar{\Gamma}_Y(w))^3}{\Psi_Y(w)}\,.$$

In the one-parameter case, the above development goes through with
$\Bar{J}_z(w)$ and $\Bar{\Gamma}_z(w)$ replaced by
$$
\Bar{K}_z(w) =
\frac{\Phi(z)^2\,\Theta_z(w)^3}{\Psi(z)\,\Psi(\tau_z(w))}
$$
and
$$ 
\Bar{\Theta}_z(w) = \sum_{m=1}^6 \prod_{n \neq m}
     \CCb{n}(\sigma_z(w))\cdot\CCb{m}(z). 
$$
First of all, there is the transformation\footnote{In
Appendix~\ref{sec:paramExp} there appears a partial
expression for $ \Bar{\Theta}_V(w)$.} of
$\Bar{\Theta}_z(w)$ under $\sigma_z$ on $\{F(z)=0\}$:
$$ 
\Bar{\Theta}_V(w) =  
  \frac{\Bar{\Theta}_z(w)|_{\{F(z)=0\}}}
       {\gamma\,\Phi(z) \Psi(z)^{16}}\,.
$$ 
Furthermore, restricting to $\{F(z)=0\}$ yields, on the
one hand, a root of $T_V(s)$
$$  \Bar{K}_z(\{q_1,q_2\}) = \frac{S_{\Bar{a}}(z)}{\mu} $$
and, on the other,
\begin{eqnarray*}
 \Bar{K}_V(w)
&=&\frac{\Phi(z)^2\,[\gamma\,\Phi(z)\,\Psi(z)^{16}\,\Bar{\Theta}_V(w)]^3\,
         |\sigma_z|^6}
        {(\eta\,\Phi(z)^2\Psi(z)^9)^8\,\Psi_V(w)}\\[5pt]
&=&-\frac{2^9 3^{50}\,\gamma^3\,V^5(V-1)^3\,\Bar{\Theta}_V(w)}
         {\eta^8\,\Psi_V(w)}\,.
\end{eqnarray*}

\subsubsection{The algorithm}

Now within reach are the ingredients required for
preparation of a root-finding
algorithm.\footnote{Available at
\emph{http://math.ucsd.edu/\~{}scrass} are
\emph{Mathematica} notebooks and supporting files that
implement this method.}  To summarize the procedure:

\begin{itemize}

\item[1)]  Select a value $A=(A_1,A_2)$ of $Y=(Y_1,Y_2)$
and, thereby, a sixth-degree resolvent $R_A(u)$. (For
sake of description, let $z \in Y_1^{-1}(A_1)\cap
Y_2^{-1}(A_2)$.  The algorithm actually finds a root
without explicitly inverting $Y_1$ or $Y_2$.)

\item[2)]  From an initial point $w_0 \in \CP{2}_w$,
iterate the map $h_A(w)$ to convergence:
\[ h^n_A(w_0) \longrightarrow \{q_1,q_2\} \in
(\Orb{72})_w \subset \CP{2}_w. \] As output take the pair
of approximate 72-points in $\CP{2}_w$
\[\{p_1,p_2\}\approx \{q_1,q_2\} =
\{\tau^{-1}_z(\pQ{a}{b}{1}),\tau^{-1}_z(\pQ{a}{b}{2})\}.
\]

\item[3)]  Using either $p_1$ or $p_2$ approximate a root
of $R_A(u)$:
\begin{eqnarray*}
\Ub{a}(z)&\approx&\mu\,\Bar{J}_A(p_1)\\
&=&\frac{\mu\,\beta^3\,(T_A\,\Bar{\Gamma}_A(p_1))^3}{\alpha^8\,\Psi_A(p_1)}\\
&=&2^{79}3^{94}(11+3\,\sqrt{15}\,i)\,
   \frac{T_A\,\Bar{\Gamma}_A(p_1))^3}{\Psi_A(p_1)}\,.
\end{eqnarray*}

\end{itemize} 
Performing the corresponding steps in the special case
produces an approximate solution to the
$V$-resolvent $T_{V_0}(s)$ given an initial choice $V_0$
of $V$:
\begin{eqnarray*}
S_{\Bar{a}}(z)&\approx&\mu\,\Bar{K}_{V_0}(p_1)\\
&=&-\frac{2^93^{50}\,\mu\,\gamma^3\,{V_0}^5({V_0}-1)^3\,\Bar{\Theta}_{V_0}(p_1)}
         {\eta^8\,\Psi_{V_0}(p_1)}\\
&=&-2^55^{-96}(11+3\,\sqrt{15}\,i)\,
\frac{{V_0}^5({V_0}-1)^3\,\Bar{\Theta}_{V_0}(p_1)}{\Psi_{V_0}(p_1)}\,.
\end{eqnarray*}

\subsection{Getting \emph{all} the roots}

Because a pair of 72-points lies on \emph{one} barred
conic, the algorithm above determines just one of the six
roots of the selected resolvent.  This manifests the
iteration's failure to break the galois symmetry
completely; the stabilizer of a 72-point is a \D{5}. 
From the coefficient field of the resolvent $R_A(u)$, the
algorithm leads to  an extended field $K$ whose galois
group $G(\Sigma/K)$ is \D{5}.  Of course, $\Sigma$ is the
splitting field for $R_A$.

  Finding all six roots calls for a dynamical system that
converges to 6-cycles in a 360 point orbit and a
root-selector function that gives equations in all the
roots.  Needless to say, this would likely complicate the
associated formulas.

%% file: appA.tex
\section{Seeing is Believing}   \label{sec:Seeing}

The following gallery of pictures provides
empirical dynamical information for some of the
special maps discussed in the text.   The program
\emph{Dynamics}, running on a Silicon Graphics
Indigo-2, created the basin plots.\footnote{The
software is the work of H. Nusse and J. Yorke
while E. Kostelich is responsible for the Unix
implementation.  See their manual \cite{NY}.}  

Complete basin images are the product of the
``BAS" routine.  In the plots for $h_{19}$, this
procedure colors a grid-cell in the event that
the trajectory of the cell's center gets close
enough to a 72-point to guarantee ultimate
attraction to the associated period-2 cycle.  The
color depends upon the destination of its
center.  If, in a specified number of iterations,
the center's trajectory fails to converge to a
pair of 72-points, the cell's color is black. 
The presence of such points would compromise
Conjecture~\ref{conj:h19}.  In the case of
$\psi_{16}$ ``restricted" to a 45-line, the
``test" attractor consists of the 45-points. 
Partial plots come by way of the ``BA"
algorithm's coloring of whole ``trajectories" of
cells, thereby manifesting some aspects of the
dynamics.\footnote{See
\cite[pp. 269ff]{NY} for a more thorough
description.}  All plots have the maximum
resolution available: a \mbox{$720 \times 720$}
grid of cells.  \emph{Mathematica} produced the
sketches of the various curves that appear.

In the basin plots for $h_{19}$ restricted to an
affine plane in one of the \Bub-\RP{2}s, the
chosen coordinates make evident the map's $D_{5}$
symmetry.  Specifically selected here is
\RR{2}{2}\ with the 1-point orbit \pD{2}{2}\ at
the origin and the 1-line orbit \LD{2}{2}\ at
infinity.  Distributed along the unit circle is
the 10-point orbit of 72-points \pP{a}{a}{1,2}
with $a=1,3,4,5,6$.  The five lines of reflective
symmetry passing through (0,0) are affine lines
in the five \RP{1} intersections with \RR{2}{2} of
both the 45-lines
\LZ{2a}{2a}\ and the
\RR{a}{a}, $a=1,3,4,5,6$.

\newpage


\begin{figure}[hp]

\resizebox{\textwidth}{!}{\includegraphics{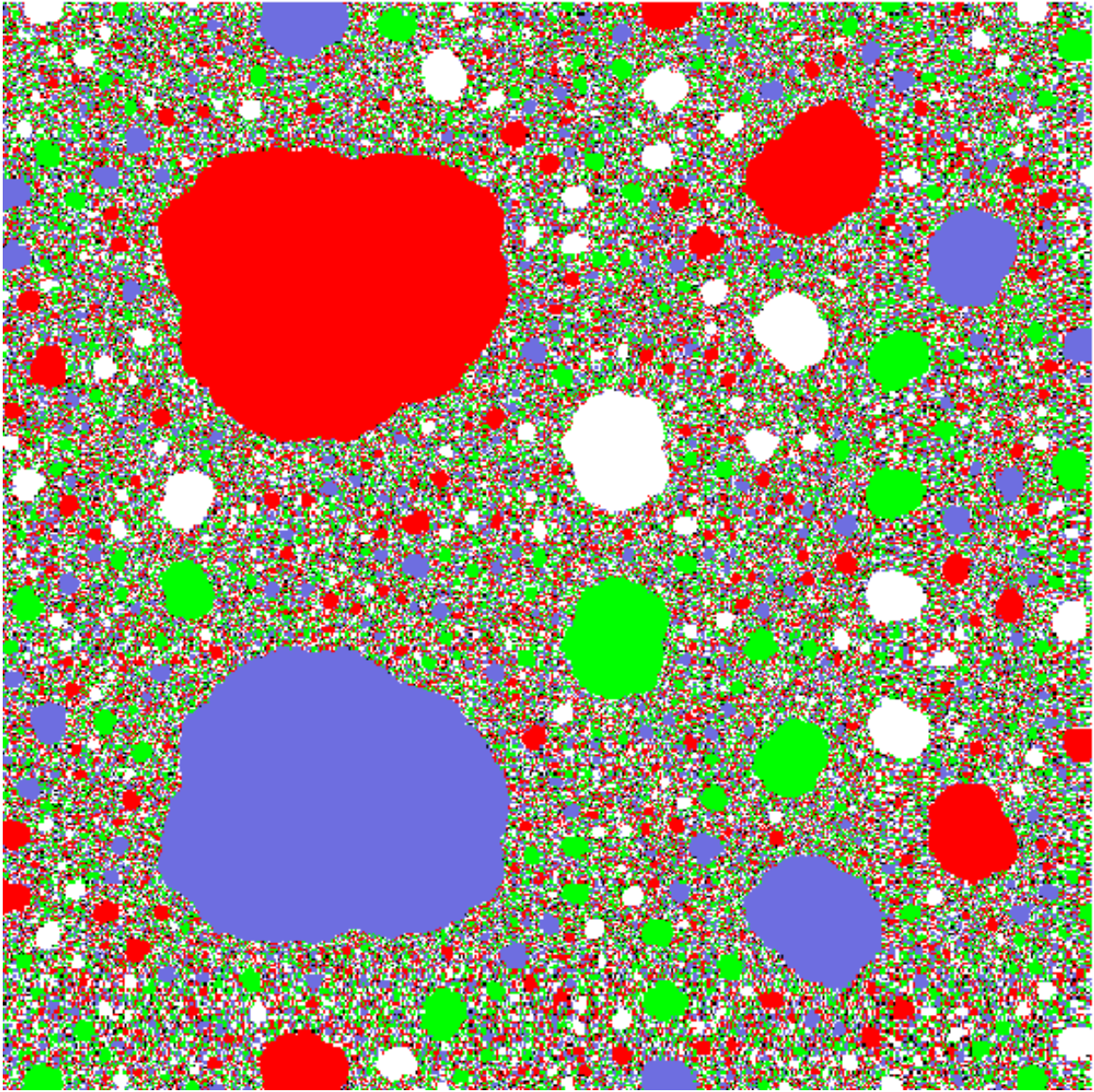}}

\caption{Dynamics of the 16-map}

\label{fig:psi16Basin}

When ``restricted" to a 45-line, the degree 16
map $\psi_{16}$ ``mostly" converges to one of the
45-points on the line.  Does this occur for
almost every point on the line?  Do the black
specks contain sets of positive measure whose
forward orbits fail to converge to one of the four
attractive 45-points that lie in the large
``central" basins? The BAS algorithm checked 60
iterates before concluding that a trajectory did
not converge.

\end{figure}

\newpage

\begin{figure}[hp]

\resizebox{\textwidth}{!}{\includegraphics{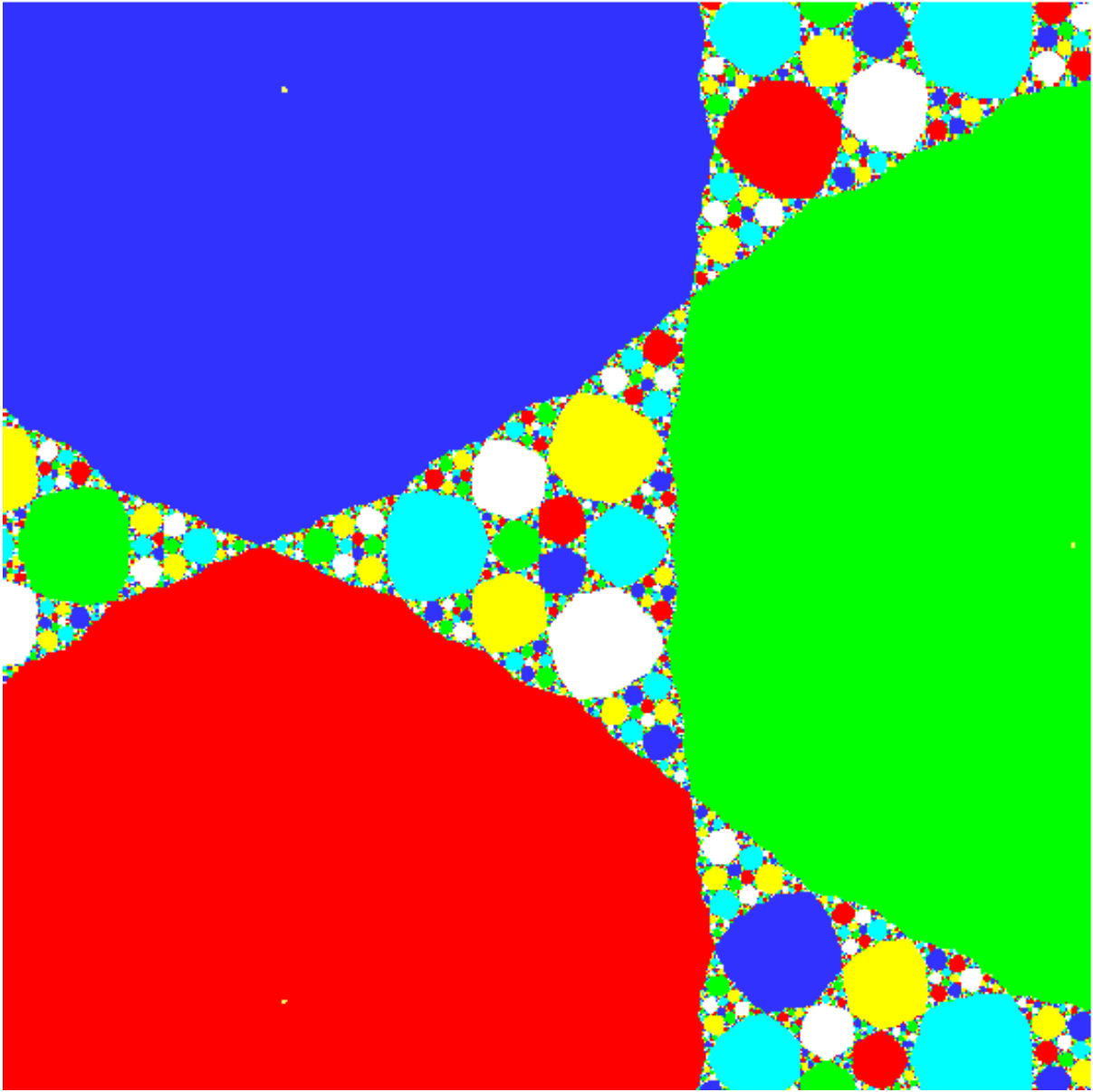}}

\caption{Icosahedral dynamics}

\label{fig:h19Conic}

The degree 19 map with icosahedral symmetry
attracts almost all points in $\CC{} \cup
\{\infty\}$ to an antipodal pair of vertices. 
Each of the six colors corresponds to such a pair
and the three large basins each contain a vertex. 
For the conic-fixing
$h_{19}$, the basin plot on each conic is
conjugate to this one.  Moreover, each basin is
the 1-dimensional intersection of a 2-dimensional
basin in \CP{2}.  Does the backward orbit of these
basins fill out \CP{2} in measure?

\end{figure}

\newpage

\begin{figure}[hp]

\resizebox{\textwidth}{!}{\includegraphics{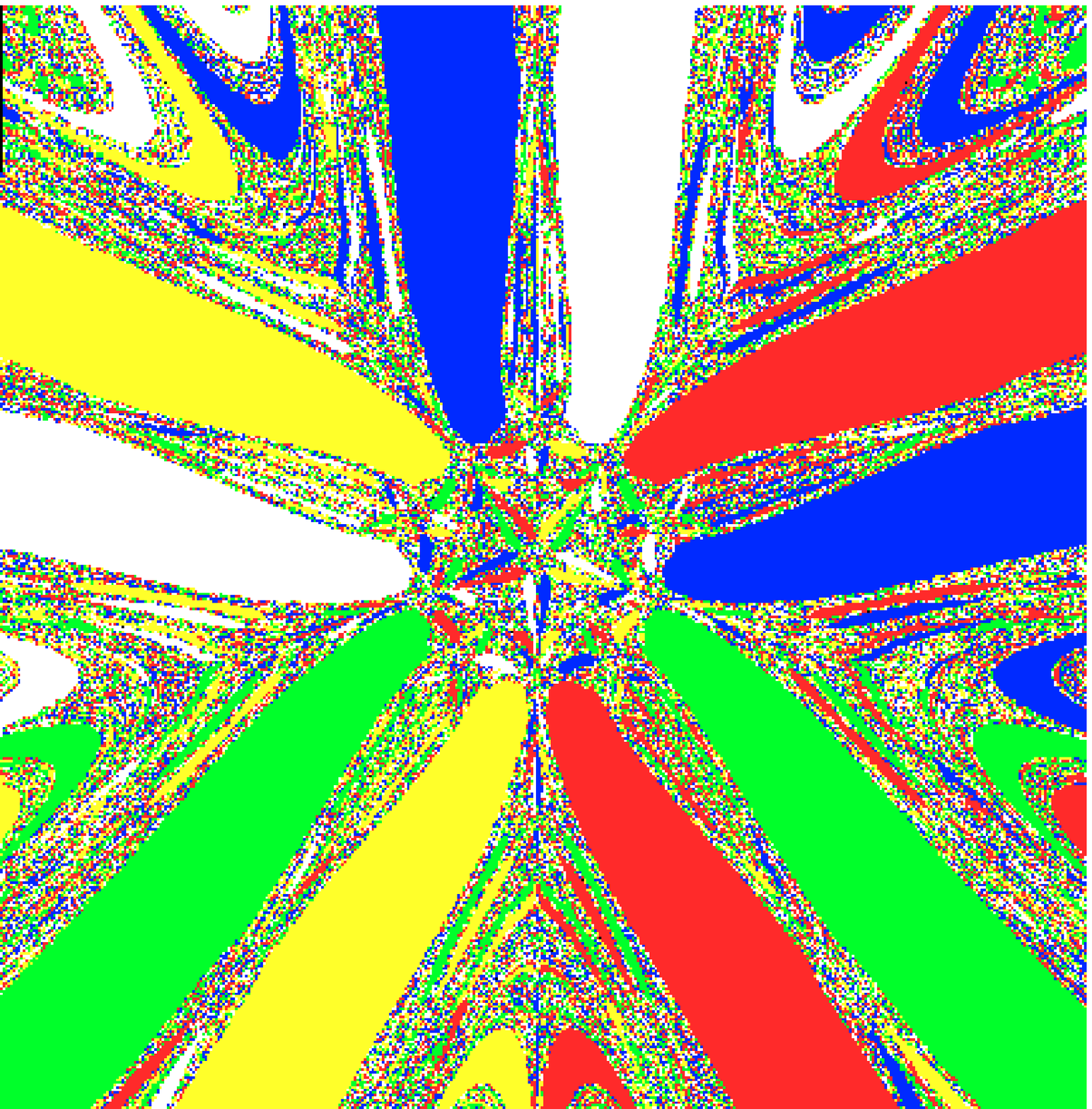}}

\caption{\RP{2} dynamics of $h_{19}$}

\label{fig:h19RP2}

For this plot the vertical and horizontal scales
are roughly from -2 to 2.  The large ``radial"
basins are immediate, i.e., each contains one of
the 72-points and come in pairs as do the period-2
attractors.  Notice the repelling behavior along
the 45-lines and particularly at their
intersection in the 36-point \pD{2}{2}. 

\end{figure}

\newpage

\begin{figure}[hp]

\resizebox{5.2in}{!}{\includegraphics{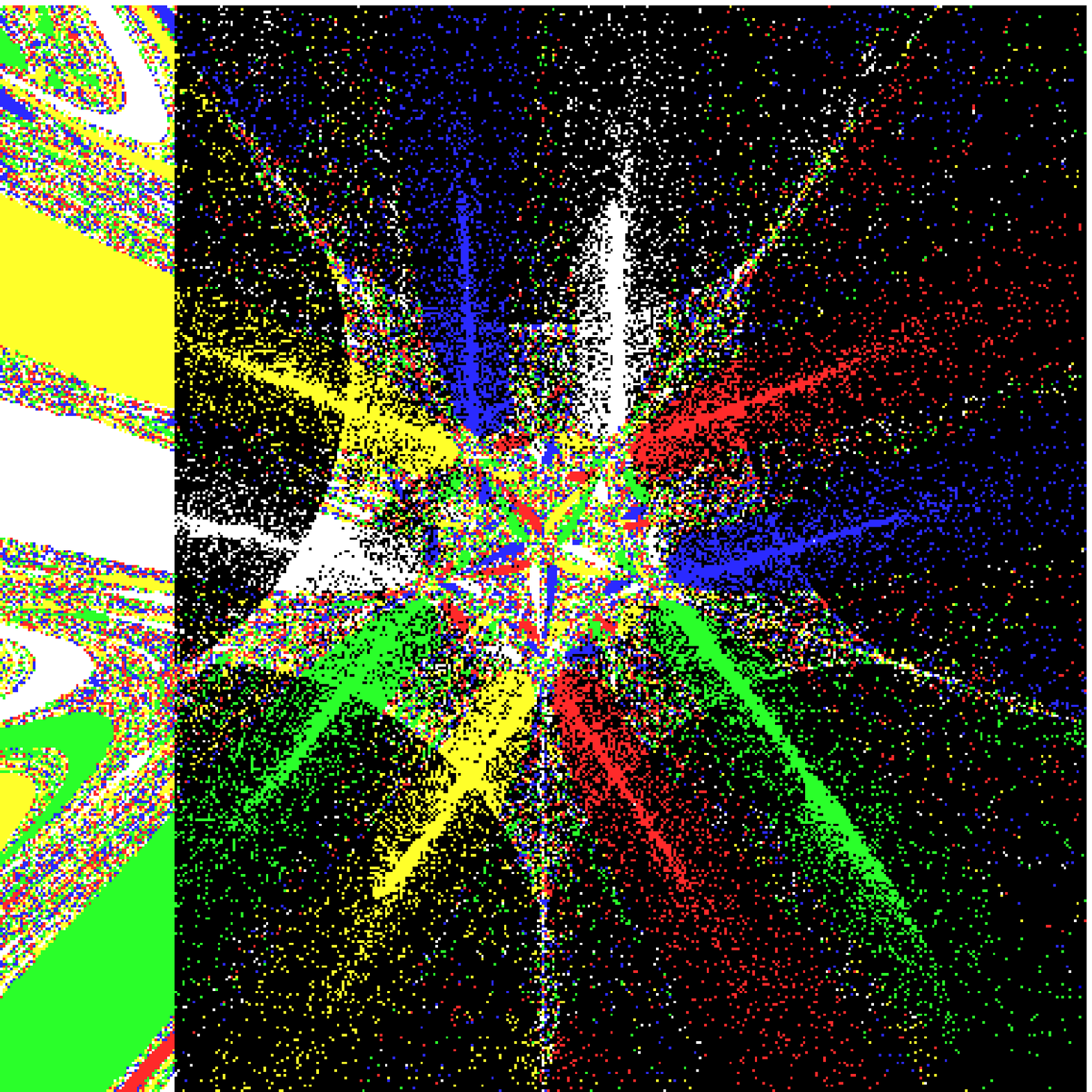}}

\caption{\RP{2} dynamics $h_{19}$}

\label{fig:h19RP2Partial}

Again, the vertical and horizontal scales are
roughly from -2 to 2. Shown here are
trajectories, colored according to their
destinations, of the points in the vertical strip
on the left.  Once more, only the five period-2
cycles of 72-points appear as attractors.  Many
of the points in the strip map inside the ``hazy
pentagon" whose vertices lie on the
45-lines---the inner ``star" is nearly filled.
``Circumscribing" this pentagon is the outer
star-like piece of the critical set shown in
Figure~\ref{fig:FOnRP2}.  Futhermore, the
pentagon seems to be the image of the inner
pentagonal oval.  Accordingly, the map folds the
plane along the pentagon's edges just outside of
which the 72-points make their presence seen in
the dense streaks at
\pP{a}{a}{1,2}.  Compare this pattern of streaks
to that of the 72-lines given in
Figure~\ref{fig:72LinesRP2}.  Since \Cb{a}\ and
\Cu{a}\ are tangent to \RR{a}{a}\ at
\pP{a}{a}{1,2}, the icosahedral 19-map opens up a
triangular angle of $\pi\over 3$ to ${4 \pi}\over
3$.  Thus, the behavior at a 72-point consists of
``fourth-powering". 
Figure~\ref{fig:h19Im36LineRP2} displays this
local ``squeezing".

\end{figure}

\begin{figure}[hp]

\resizebox{\textwidth}{!}{\includegraphics{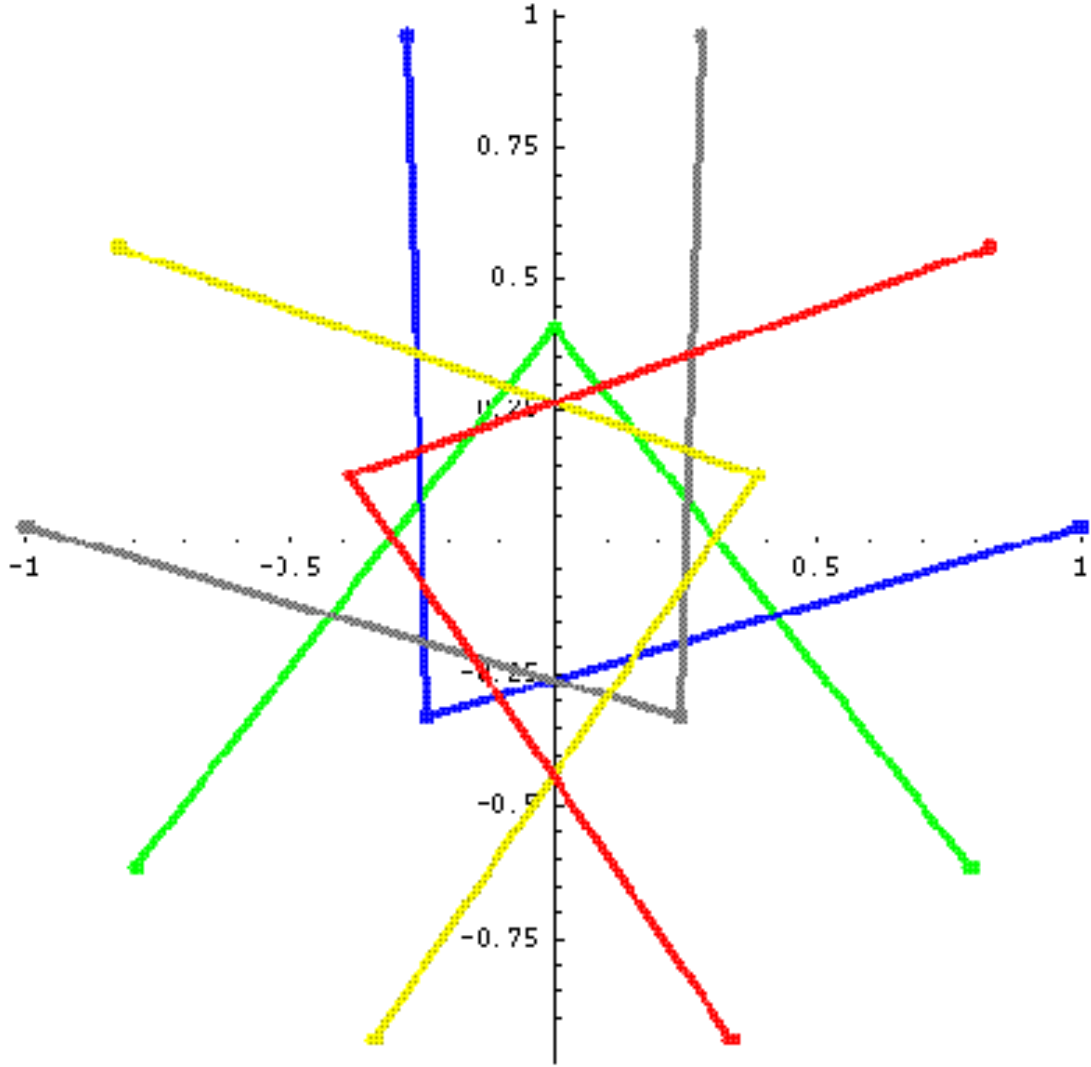}}

\caption{Configuration of 72-lines on a
\Bub-\RP{2}}

\label{fig:72LinesRP2}

Under \bub{2}{2}, five pairs of 72-lines 
$$
\{\{\LP{a}{a}{1},\LP{a}{a}{2}\}|a=1,\ldots.6\}
$$
map to themselves.  Accordingly, each line meets
\RR{2}{2}\ in an \RP{1}.  The picture shows their
configuration in the affine plane of
Figures~\ref{fig:h19RP2} and
\ref{fig:h19RP2Partial}.  Each pair receives a
single color according to the scheme of the basin
plots.  A given pair \LP{a}{a}{1}\ and
\LP{a}{a}{2}\ passes through the 72-points
\pP{a}{a}{2}\ and \pP{a}{a}{1}\ respectively;
they intersect in the corresponding repelling and
fixed 36-point \pP{a}{a}.

\end{figure}

\begin{figure}[hp]

\resizebox{\textwidth}{!}{\includegraphics{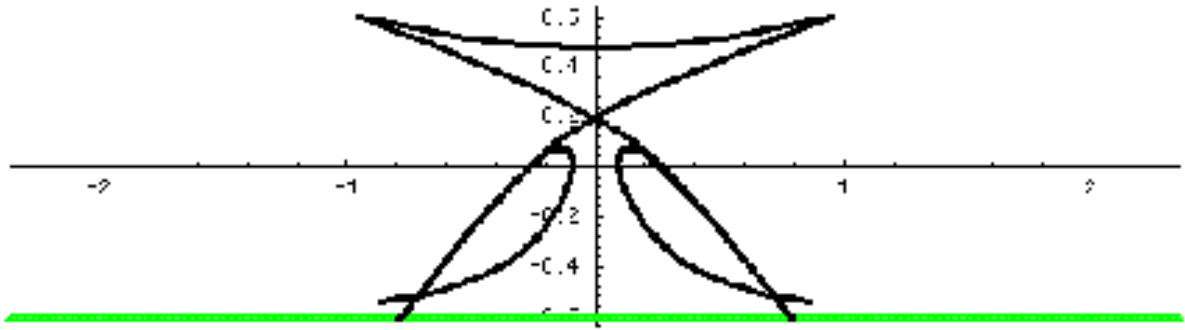}}

\caption{Image on a \Bub-\RP{2}\ of a 36-line
under $h_{19}$}

\label{fig:h19Im36LineRP2}

The green horizontal line corresponds to the
\RP{1} intersection of
\RR{2}{2} and the 36-line associated with the
pair of green 72-points from the basin plots. 
The dark curve is the image of the line under
$h_{19}$.  Sitting at the sharp cusps are the
72-points which the map exchanges.  As indicated
in the caption to Figure~\ref{fig:h19RP2Partial},
the line folds over at these critical points. 
The upper two sharp turns are not critical
values;  they occur where the line passes through
the yellow and red ``streak" that approximate
72-lines.

\end{figure}

\begin{figure}[hp]

\resizebox{\textwidth}{!}{\includegraphics{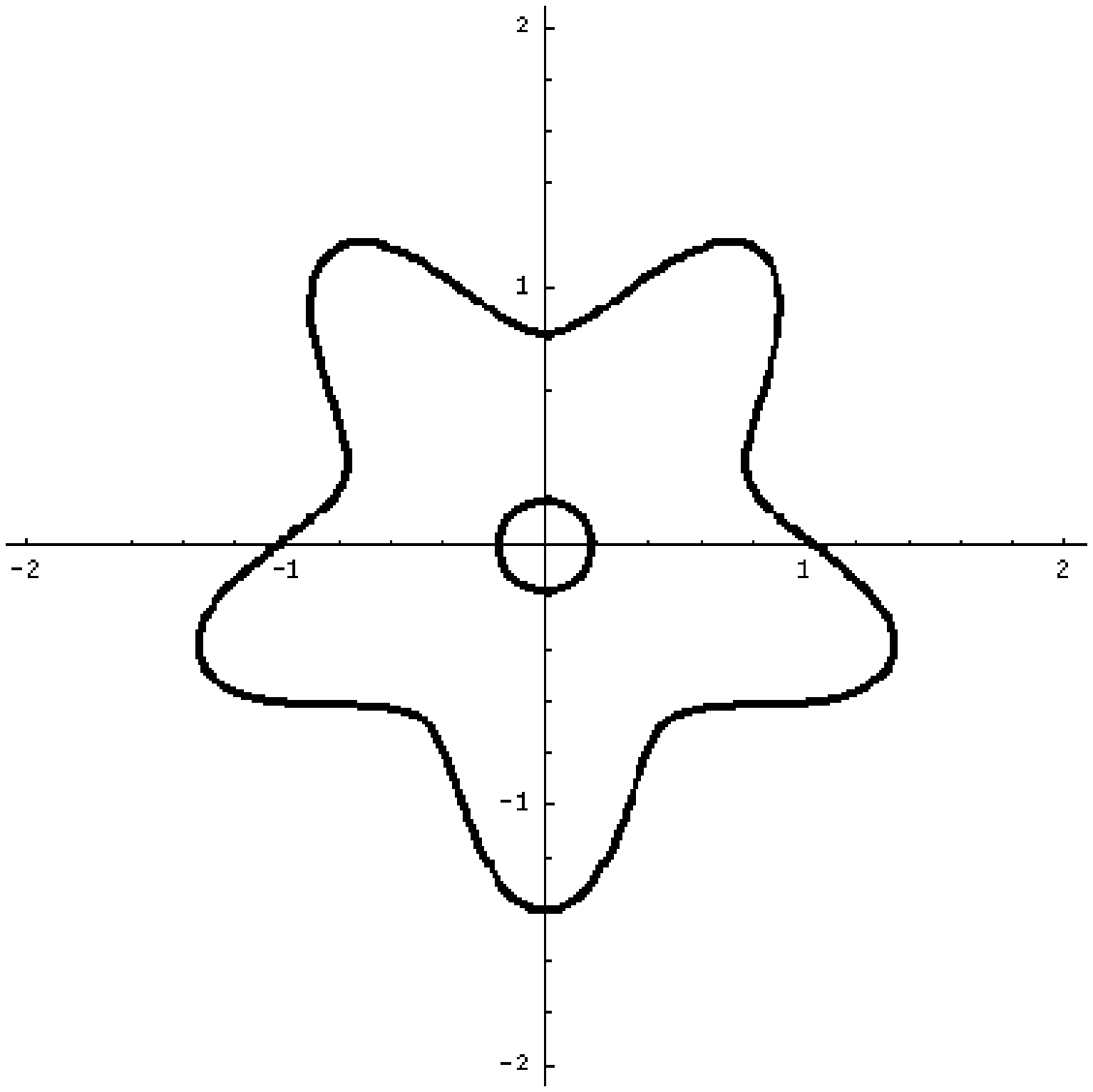}}

\caption{Critical set of $h_{19}$ on a
\Bub-\RP{2}}

\label{fig:FOnRP2}

Here is a \emph{Mathematica} contour plot of the
sixth degree curve $\{F=0\}$ on \RR{2}{2}.  At
the 10 inflection points are the superattracting
72-points.  Further computation suggests that
$\{G_{48}=0\}$ hits the \RP{2}\ in a discrete
set.  Might this set consist only of the singular
72-points?

\end{figure}

%% file: appB.tex
\section{Parametrized Expressions   
\label{sec:paramExp}}

\subsection{Sixth-degree forms}

After solving a system of linear equations for
each coefficient of the degree 6 monomials in $w$,
the
simplifying change of coordinates
$$ \begin{array}{lll}
w_1 &\rightarrow& 4\,(2\,w_1 - 23\,w_2 +
200\,w_3)\\
w_2 &\rightarrow& 2\,(w_1 - 52\,w_2 + 64\,w_3)\\
w_3 &\rightarrow& 6\,w_3
\end{array} $$
leaves the following expression for the
2-parameter family of sixth-degree Valentiner
invariants:

\small
\vspace{10pt}
\input{FInY}
\vspace{20pt}
\normalsize

In the case of the special 1-parameter family,
the coordinate change
$$ \begin{array}{rll}
3\,w_1 &\rightarrow& w_2\\
18\,w_2 &\rightarrow& (3\,w_1 + 5\,w_2 -
5\,w_3)\\
12\,w_3 &\rightarrow& (w_1 - w_3)
\end{array} $$
yields:

\small
\vspace{10pt}
\input{FInV}
\normalsize

\subsection{Root-selectors}

In the coordinates given above, the expressions
for the respective factors in the root-selecting
functions $\Bar{J}_Y(w)$ and $\Bar{K}_V(w)$ take
the forms

\small
\vspace{10pt}
\input{GammabarInY}
\vspace{20pt}
\input{ThetabarInV}
\normalsize

\subsection{Degree 19 maps}
Using (\ref{eq:f64}), the identities 
\begin{eqnarray*}
\alpha &=& \frac{1}{2^{10}\,3^{12}}\\
 |\tau_z|^2&=&F(z)^{25}\,T_Y,\\
|\sigma_z|^2&=&-2^8\,3^{17}\,\Phi(z)^2\,\Psi(z)^9\,V^2(V-1),
\end{eqnarray*}
and the results in Section~\ref{sec:ParamVals},
the following
procedures yield the conic-fixing 19-maps.  For
the 2-parameter case, substitution gives
\begin{eqnarray*} 
 f_{64}(y)
&=&\frac{(\alpha\,F(z)^{25})^{17}}{|\tau_z|^{13}}\,\tau_z
\left\{
\left[10\,\left(\frac{T_Y}{\alpha}\right)^4\,F_Y^6\,\Phi_Y
+      
100\,\left(\frac{T_Y}{\alpha}\right)^3\,F_Y^4\,\Phi_Y^2
+
       \right.\right.\\   
&&45\,\left(\frac{T_Y}{\alpha}\right)^2\,F_Y^2\,\Phi_Y^3
+      
156\,\left(\frac{T_Y}{\alpha}\right)\,\Phi_Y^4
+      
39\,\left(\frac{T_Y}{\alpha}\right)^2\,F_Y^3\,\Psi_Y
+ \\
    &&\left.\left.      
51\,\left(\frac{T_Y}{\alpha}\right)\,F_Y\,\Phi_Y\,\Psi_Y
\right]
            \cdot \psi_Y(w)  - 
       27\,\Psi_Y \cdot \phi_Y(w) + 54\,\Phi_Y^2
\cdot f_Y(w) \right\}\\
   &=&\frac{\alpha\,F(z)^{25}}{|\tau_z|^{13}}\
\;\tau_z\,(f_{64})_Y.
   \end{eqnarray*}
Thus,
\begin{eqnarray*} 
 f_{19}(y) &=& \frac{f_{64}(y)}{X(y)}\\
 &=& \frac{(\alpha\,F(z)^{25})^5}{|\tau_z|^4}\
\;\tau_z\,
     \frac{(f_{64})_Y(w)}{X_Y(w)}\\
 &=& \frac{(\alpha\,F(z)^{25})^5}{|\tau_z|^4}\
\;\tau_z\,(f_{19})_Y(w).
\end{eqnarray*}
As for the special 19-map, let
$f_Y(w)=(f_{19})_Y(w)$ so that
\begin{eqnarray*} 
 h_{19}(y) &=&1620\,F(y)^3 \cdot y + f_{19}(y)\\
      &=&1620\,(\alpha\,F(z)^{25}\,F_Y(w))^3
\cdot \tau_z(w) + 
        \frac{(\alpha\,F(z)^{25})^5}{|\tau_z|^4}\
\;\tau_z\,f_Y(w)\\
     &=&\frac{(\alpha\,F(z)^{25})^5}{|\tau_z|^4}\
\;\tau_z
       \left(1620\,
\frac{|\tau_z|^4}{(\alpha\,F(z)^{25})^2}\,
           F_Y(w)^3 \cdot w + f_Y(w)
\right)\\     
&=&\frac{(\alpha\,F(z)^{25})^5}{|\tau_z|^4}\
\;\tau_z \left(
         1620\,\left(\frac{T_Y}{\alpha}\right)^2
F_Y(w)^3 \cdot w +
          f_Y(w) \right).
\end{eqnarray*}
The family of conic-fixing maps on \CP{2}\ is
\begin{eqnarray*} 
 h_Y(w)&=&1620\,\left(\frac{T_Y}{\alpha}\right)^2
F_Y(w)^3 \cdot w + 
           f_Y(w)\\[10pt]
       &=&(2^{12}\cdot 3^{16}\cdot 5)\,T_Y^2
F_Y(w)^3 \cdot w + f_Y(w)
\end{eqnarray*} 
Similar calculations in the 1-parameter setting
yield
\begin{eqnarray*} 
h_V(w)&=&1620\,\left(-\frac{64}{9}\right)\,V^4(V-1)^2\,F_V(w)^3
\cdot w +
f_V(w)\\[10pt]
       &=&11520\,V^4(V-1)^2\,F_V(w)^3 \cdot w +
f_V(w).
\end{eqnarray*}

%% file: FInY.tex
\noindent  
$F_Y(w) = 
  \displaystyle{(2^{10}\cdot
3^{12})\,\frac{F(\tau_z(w))}{F(z)^{25}}} =$
   \vspace{10pt}
\begin{sloppypar} \raggedright \noindent
$11664\, (-{{w_1}^2} + 2\,w_1\,w_2 - {{w_2}^2} - 
     14\,w_1\,w_3 + 14\,w_2\,w_3 - 39\,{{w_3}^2}
      ) \,( -{{w_1}^2} + 2\,w_1\,w_2 - 
     {{w_2}^2} - 4\,w_1\,w_3 + 4\,w_2\,w_3 + 
     11\,{{w_3}^2} ) \,
   ( {{w_1}^2} - 2\,w_1\,w_2 + {{w_2}^2} + 
     12\,w_1\,w_3 - 12\,w_2\,w_3 + 45\,{{w_3}^2}
      )  + $\\[5pt]
$7776\,w_3\,
   ( -127\,{{w_1}^5} + 650\,{{w_1}^4}\,w_2 - 
     1285\,{{w_1}^3}\,{{w_2}^2} + 
     1225\,{{w_1}^2}\,{{w_2}^3} -
560\,w_1\,{{w_2}^4} + 
     97\,{{w_2}^5} - 3790\,{{w_1}^4}\,w_3 + 
     15760\,{{w_1}^3}\,w_2\,w_3 - 
     23820\,{{w_1}^2}\,{{w_2}^2}\,w_3 + 
     15520\,w_1\,{{w_2}^3}\,w_3 - 
     3670\,{{w_2}^4}\,w_3 -
     44385\,{{w_1}^3}\,{{w_3}^2} + 
     140715\,{{w_1}^2}\,w_2\,{{w_3}^2} - 
     144360\,w_1\,{{w_2}^2}\,{{w_3}^2} + 
     48030\,{{w_2}^3}\,{{w_3}^2} - 
     260580\,{{w_1}^2}\,{{w_3}^3} + 
     558960\,w_1\,w_2\,{{w_3}^3} - 
     290820\,{{w_2}^2}\,{{w_3}^3} - 
     756816\,w_1\,{{w_3}^4} +
     822741\,w_2\,{{w_3}^4} - 
     831054\,{{w_3}^5} ) \,Y_1 +$\\[5pt] 
$54\,( -1211\,{{w_1}^6} + 10374\,{{w_1}^5}\,w_2 - 
     31275\,{{w_1}^4}\,{{w_2}^2} + 
     42880\,{{w_1}^3}\,{{w_2}^3} - 
     26160\,{{w_1}^2}\,{{w_2}^4} +
     4176\,w_1\,{{w_2}^5} + 
     1216\,{{w_2}^6} - 79538\,{{w_1}^5}\,w_3 + 
     490810\,{{w_1}^4}\,w_2\,w_3 - 
     1044200\,{{w_1}^3}\,{{w_2}^2}\,w_3 + 
     939920\,{{w_1}^2}\,{{w_2}^3}\,w_3 - 
     307120\,w_1\,{{w_2}^4}\,w_3 + 
     128\,{{w_2}^5}\,w_3 - 
     1619325\,{{w_1}^4}\,{{w_3}^2} + 
     8179200\,{{w_1}^3}\,w_2\,{{w_3}^2} - 
     13473000\,{{w_1}^2}\,{{w_2}^2}\,{{w_3}^2} + 
     8657280\,w_1\,{{w_2}^3}\,{{w_3}^2} - 
     1728360\,{{w_2}^4}\,{{w_3}^2} - 
     14229000\,{{w_1}^3}\,{{w_3}^3} + 
     58644000\,{{w_1}^2}\,w_2\,{{w_3}^3} - 
     68856480\,w_1\,{{w_2}^2}\,{{w_3}^3} + 
     23945760\,{{w_2}^3}\,{{w_3}^3} - 
     64072440\,{{w_1}^2}\,{{w_3}^4} + 
     193407120\,w_1\,w_2\,{{w_3}^4} - 
     119268000\,{{w_2}^2}\,{{w_3}^4} - 
     116140176\,w_1\,{{w_3}^5} + 
     210347712\,w_2\,{{w_3}^5} +
     49729896\,{{w_3}^6} ) \,{{Y_1}^2} +$\\[5pt]
$(48469\,{{w_1}^6} - 
     743244\,{{w_1}^5}\,w_2 + 
     4444890\,{{w_1}^4}\,{{w_2}^2} - 
     13073120\,{{w_1}^3}\,{{w_2}^3} + 
     19723200\,{{w_1}^2}\,{{w_2}^4} - 
     14608896\,w_1\,{{w_2}^5} +
     4205056\,{{w_2}^6} + 
     1051200\,{{w_1}^5}\,w_3 - 
     21117960\,{{w_1}^4}\,w_2\,w_3 + 
     118315440\,{{w_1}^3}\,{{w_2}^2}\,w_3 - 
     259434720\,{{w_1}^2}\,{{w_2}^3}\,w_3 + 
     249145920\,w_1\,{{w_2}^4}\,w_3 - 
     88047360\,{{w_2}^5}\,w_3 + 
     21207690\,{{w_1}^4}\,{{w_3}^2} - 
     76688640\,{{w_1}^3}\,w_2\,{{w_3}^2} + 
     441644400\,{{w_1}^2}\,{{w_2}^2}\,{{w_3}^2} - 
     809291520\,w_1\,{{w_2}^3}\,{{w_3}^2} + 
     439814880\,{{w_2}^4}\,{{w_3}^2} + 
     815030640\,{{w_1}^3}\,{{w_3}^3} - 
     868838400\,{{w_1}^2}\,w_2\,{{w_3}^3} + 
     529260480\,w_1\,{{w_2}^2}\,{{w_3}^3} - 
     665167680\,{{w_2}^3}\,{{w_3}^3} + 
     8239760640\,{{w_1}^2}\,{{w_3}^4} - 
     10091148480\,w_1\,w_2\,{{w_3}^4} + 
     4723181280\,{{w_2}^2}\,{{w_3}^4} + 
     50452492032\,w_1\,{{w_3}^5} - 
     45443177280\,w_2\,{{w_3}^5} +
     145956110976\,{{w_3}^6}
) \,{{Y_1}^3} +$\\[5pt] 
$6\,( 43274\,{{w_1}^6} - 667704\,{{w_1}^5}\,w_2 + 
     4028205\,{{w_1}^4}\,{{w_2}^2} - 
     11998480\,{{w_1}^3}\,{{w_2}^3} + 
     18432480\,{{w_1}^2}\,{{w_2}^4} - 
     13967616\,w_1\,{{w_2}^5} +
     4111616\,{{w_2}^6} + 
     1805328\,{{w_1}^5}\,w_3 - 
     19051740\,{{w_1}^4}\,w_2\,w_3 + 
     75568680\,{{w_1}^3}\,{{w_2}^2}\,w_3 - 
     143512560\,{{w_1}^2}\,{{w_2}^3}\,w_3 + 
     134468640\,w_1\,{{w_2}^4}\,w_3 - 
     49479552\,{{w_2}^5}\,w_3 + 
     23063490\,{{w_1}^4}\,{{w_3}^2} - 
     97315920\,{{w_1}^3}\,w_2\,{{w_3}^2} + 
     127881720\,{{w_1}^2}\,{{w_2}^2}\,{{w_3}^2} - 
     77891040\,w_1\,{{w_2}^3}\,{{w_3}^2} + 
     44389440\,{{w_2}^4}\,{{w_3}^2} + 
     409273920\,{{w_1}^3}\,{{w_3}^3} - 
     1276901280\,{{w_1}^2}\,w_2\,{{w_3}^3} + 
     1247520960\,w_1\,{{w_2}^2}\,{{w_3}^3} - 
     650304000\,{{w_2}^3}\,{{w_3}^3} + 
     3052420200\,{{w_1}^2}\,{{w_3}^4} - 
     12009740400\,w_1\,w_2\,{{w_3}^4} + 
     10525836600\,{{w_2}^2}\,{{w_3}^4} + 
     26559204048\,w_1\,{{w_3}^5} - 
     45364278096\,w_2\,{{w_3}^5} +
     101461763160\,{{w_3}^6}
) \,{{Y_1}^4} + $\\[5pt]
$2\,( 105901\,{{w_1}^6} -
1755600\,{{w_1}^5}\,w_2 + 
     11616675\,{{w_1}^4}\,{{w_2}^2} - 
     38949440\,{{w_1}^3}\,{{w_2}^3} + 
     69116640\,{{w_1}^2}\,{{w_2}^4} - 
     60716544\,w_1\,{{w_2}^5} +
     19859200\,{{w_2}^6} - 
     6063192\,{{w_1}^5}\,w_3 + 
     82652760\,{{w_1}^4}\,w_2\,w_3 - 
     412432200\,{{w_1}^3}\,{{w_2}^2}\,w_3 + 
     898113600\,{{w_1}^2}\,{{w_2}^3}\,w_3 - 
     786384000\,w_1\,{{w_2}^4}\,w_3 + 
     220981248\,{{w_2}^5}\,w_3 - 
     146310030\,{{w_1}^4}\,{{w_3}^2} + 
     1634297040\,{{w_1}^3}\,w_2\,{{w_3}^2} - 
     6098117400\,{{w_1}^2}\,{{w_2}^2}\,{{w_3}^2}
+ 
     8499487680\,w_1\,{{w_2}^3}\,{{w_3}^2} - 
     3567715200\,{{w_2}^4}\,{{w_3}^2} + 
     880886880\,{{w_1}^3}\,{{w_3}^3} - 
     8797422960\,{{w_1}^2}\,w_2\,{{w_3}^3} + 
     21711123360\,w_1\,{{w_2}^2}\,{{w_3}^3} - 
     19644958080\,{{w_2}^3}\,{{w_3}^3} + 
     16531620480\,{{w_1}^2}\,{{w_3}^4} - 
     146923709760\,w_1\,w_2\,{{w_3}^4} + 
     177054446160\,{{w_2}^2}\,{{w_3}^4} + 
     281903296896\,w_1\,{{w_3}^5} - 
     514798256592\,w_2\,{{w_3}^5} +
     1116497698848\,{{w_3}^6}
) \,{{Y_1}^5} + $\\[5pt]
$2\,( -60989\,{{w_1}^6} +
     1013970\,{{w_1}^5}\,w_2 - 
     6509550\,{{w_1}^4}\,{{w_2}^2} + 
     19983520\,{{w_1}^3}\,{{w_2}^3} - 
     29628480\,{{w_1}^2}\,{{w_2}^4} + 
     21098496\,w_1\,{{w_2}^5} -
     11333120\,{{w_2}^6} - 
     10304730\,{{w_1}^5}\,w_3 + 
     141141510\,{{w_1}^4}\,w_2\,w_3 - 
     711217440\,{{w_1}^3}\,{{w_2}^2}\,w_3 + 
     1580051520\,{{w_1}^2}\,{{w_2}^3}\,w_3 - 
     1444147200\,w_1\,{{w_2}^4}\,w_3 + 
     433497600\,{{w_2}^5}\,w_3 + 
     7910865\,{{w_1}^4}\,{{w_3}^2} - 
     615917520\,{{w_1}^3}\,w_2\,{{w_3}^2} + 
     4948801920\,{{w_1}^2}\,{{w_2}^2}\,{{w_3}^2}
- 
     12996339840\,w_1\,{{w_2}^3}\,{{w_3}^2} + 
     10198068480\,{{w_2}^4}\,{{w_3}^2} + 
     6077720520\,{{w_1}^3}\,{{w_3}^3} - 
     63144917280\,{{w_1}^2}\,w_2\,{{w_3}^3} + 
     192593116800\,w_1\,{{w_2}^2}\,{{w_3}^3} - 
     183921926400\,{{w_2}^3}\,{{w_3}^3} + 
     58576870440\,{{w_1}^2}\,{{w_3}^4} - 
     408475787760\,w_1\,w_2\,{{w_3}^4} + 
     511734905280\,{{w_2}^2}\,{{w_3}^4} + 
     680591915568\,w_1\,{{w_3}^5} - 
     892873018080\,w_2\,{{w_3}^5} +
     1956383853480\,{{w_3}^6}
) \,{{Y_1}^6} +$\\[5pt] 
$ 9\,( 14893\,{{w_1}^6} - 291876\,{{w_1}^5}\,w_2
+ 
     2263200\,{{w_1}^4}\,{{w_2}^2} - 
     8574080\,{{w_1}^3}\,{{w_2}^3} + 
     15233280\,{{w_1}^2}\,{{w_2}^4} - 
     7590912\,w_1\,{{w_2}^5} - 6127616\,{{w_2}^6}
+ 
     1600536\,{{w_1}^5}\,w_3 - 
     28552200\,{{w_1}^4}\,w_2\,w_3 + 
     197250240\,{{w_1}^3}\,{{w_2}^2}\,w_3 - 
     650334720\,{{w_1}^2}\,{{w_2}^3}\,w_3 + 
     995343360\,w_1\,{{w_2}^4}\,w_3 - 
     529618944\,{{w_2}^5}\,w_3 + 
     150441690\,{{w_1}^4}\,{{w_3}^2} - 
     1985417760\,{{w_1}^3}\,w_2\,{{w_3}^2} + 
     9512081280\,{{w_1}^2}\,{{w_2}^2}\,{{w_3}^2}
- 
     19309670400\,w_1\,{{w_2}^3}\,{{w_3}^2} + 
     13599045120\,{{w_2}^4}\,{{w_3}^2} + 
     2795123760\,{{w_1}^3}\,{{w_3}^3} - 
     27544011840\,{{w_1}^2}\,w_2\,{{w_3}^3} + 
     87759394560\,w_1\,{{w_2}^2}\,{{w_3}^3} - 
     89221309440\,{{w_2}^3}\,{{w_3}^3} + 
     10326191040\,{{w_1}^2}\,{{w_3}^4} - 
     42916236480\,w_1\,w_2\,{{w_3}^4} + 
     17753904000\,{{w_2}^2}\,{{w_3}^4} + 
     100671672960\,w_1\,{{w_3}^5} + 
     161852587776\,w_2\,{{w_3}^5} -
     169776974592\,{{w_3}^6}
) \,{{Y_1}^7} + $\\[5pt]
$8\,( 16238\,{{w_1}^6} - 291297\,{{w_1}^5}\,w_2
+ 
     1928820\,{{w_1}^4}\,{{w_2}^2} - 
     5038240\,{{w_1}^3}\,{{w_2}^3} - 
     631680\,{{w_1}^2}\,{{w_2}^4} + 
     26204928\,w_1\,{{w_2}^5} -
     34266112\,{{w_2}^6} + 
     2811906\,{{w_1}^5}\,w_3 - 
     49346640\,{{w_1}^4}\,w_2\,w_3 + 
     339641280\,{{w_1}^3}\,{{w_2}^2}\,w_3 - 
     1138037760\,{{w_1}^2}\,{{w_2}^3}\,w_3 + 
     1835020800\,w_1\,{{w_2}^4}\,w_3 - 
     1115172864\,{{w_2}^5}\,w_3 + 
     134656560\,{{w_1}^4}\,{{w_3}^2} - 
     1881693720\,{{w_1}^3}\,w_2\,{{w_3}^2} + 
     9653294880\,{{w_1}^2}\,{{w_2}^2}\,{{w_3}^2}
- 
     21377139840\,w_1\,{{w_2}^3}\,{{w_3}^2} + 
     17012160000\,{{w_2}^4}\,{{w_3}^2} + 
     719438220\,{{w_1}^3}\,{{w_3}^3} - 
     6652711440\,{{w_1}^2}\,w_2\,{{w_3}^3} + 
     19475976960\,w_1\,{{w_2}^2}\,{{w_3}^3} - 
     17504570880\,{{w_2}^3}\,{{w_3}^3} - 
     24811199100\,{{w_1}^2}\,{{w_3}^4} + 
     218557854720\,w_1\,w_2\,{{w_3}^4} - 
     442688885280\,{{w_2}^2}\,{{w_3}^4} - 
     132608756664\,w_1\,{{w_3}^5} + 
     960112922112\,w_2\,{{w_3}^5} -
     1682221113576\,{{w_3}^6}
) \,{{Y_1}^8} + $\\[5pt]
$108\,( -171\,{{w_1}^6} + 4104\,{{w_1}^5}\,w_2 - 
     41040\,{{w_1}^4}\,{{w_2}^2} + 
     218880\,{{w_1}^3}\,{{w_2}^3} - 
     656640\,{{w_1}^2}\,{{w_2}^4} + 
     1050624\,w_1\,{{w_2}^5} - 700416\,{{w_2}^6}
+ 
     47962\,{{w_1}^5}\,w_3 - 
     959240\,{{w_1}^4}\,w_2\,w_3 + 
     7673920\,{{w_1}^3}\,{{w_2}^2}\,w_3 - 
     30695680\,{{w_1}^2}\,{{w_2}^3}\,w_3 + 
     61391360\,w_1\,{{w_2}^4}\,w_3 - 
     49113088\,{{w_2}^5}\,w_3 - 
     394570\,{{w_1}^4}\,{{w_3}^2} + 
     7285120\,{{w_1}^3}\,w_2\,{{w_3}^2} - 
     49542720\,{{w_1}^2}\,{{w_2}^2}\,{{w_3}^2} + 
     147665920\,w_1\,{{w_2}^3}\,{{w_3}^2} - 
     163217920\,{{w_2}^4}\,{{w_3}^2} - 
     153027360\,{{w_1}^3}\,{{w_3}^3} + 
     1762456320\,{{w_1}^2}\,w_2\,{{w_3}^3} - 
     6754337280\,w_1\,{{w_2}^2}\,{{w_3}^3} + 
     8611799040\,{{w_2}^3}\,{{w_3}^3} - 
     2805037200\,{{w_1}^2}\,{{w_3}^4} + 
     22007304000\,w_1\,w_2\,{{w_3}^4} - 
     43148620800\,{{w_2}^2}\,{{w_3}^4} - 
     15525586656\,w_1\,{{w_3}^5} + 
     67663897344\,w_2\,{{w_3}^5} -
     163361647584\,{{w_3}^6}
) \,{{Y_1}^9} + $\\[5pt]
$432\,w_3\,( -6075\,{{w_1}^5} + 
     121500\,{{w_1}^4}\,w_2 - 
     972000\,{{w_1}^3}\,{{w_2}^2} + 
     3888000\,{{w_1}^2}\,{{w_2}^3} - 
     7776000\,w_1\,{{w_2}^4} + 6220800\,{{w_2}^5} - 
     617700\,{{w_1}^4}\,w_3 + 
     9883200\,{{w_1}^3}\,w_2\,w_3 - 
     59299200\,{{w_1}^2}\,{{w_2}^2}\,w_3 + 
     158131200\,w_1\,{{w_2}^3}\,w_3 - 
     158131200\,{{w_2}^4}\,w_3 - 
     20973400\,{{w_1}^3}\,{{w_3}^2} + 
     251680800\,{{w_1}^2}\,w_2\,{{w_3}^2} - 
     1006723200\,w_1\,{{w_2}^2}\,{{w_3}^2} + 
     1342297600\,{{w_2}^3}\,{{w_3}^2} - 
     201064680\,{{w_1}^2}\,{{w_3}^3} + 
     1584217440\,w_1\,w_2\,{{w_3}^3} - 
     3119834880\,{{w_2}^2}\,{{w_3}^3} - 
     1374418368\,w_1\,{{w_3}^4} + 
     4253513472\,w_2\,{{w_3}^4} -
     25378959744\,{{w_3}^5}
) \,{{Y_1}^{10}} + $\\[5pt]
$23328000\,{{w_3}^3}\,( -5\,{{w_1}^3} + 
     60\,{{w_1}^2}\,w_2 - 240\,w_1\,{{w_2}^2} + 
     320\,{{w_2}^3} + 575\,{{w_1}^2}\,w_3 - 
     4600\,w_1\,w_2\,w_3 + 9200\,{{w_2}^2}\,w_3 - 
     1202\,w_1\,{{w_3}^2} + 4808\,w_2\,{{w_3}^2} - 
     124918\,{{w_3}^3} ) \,{{Y_1}^{11}} + $\\[5pt]
$2332800000\,{{w_3}^5}\,( 3\,w_1 - 12\,w_2 + 
     176\,w_3 ) \,{{Y_1}^{12}} + $\\[5pt]
 $139968\,w_3\,( 23\,{{w_1}^5} - 
     145\,{{w_1}^4}\,w_2 + 395\,{{w_1}^3}\,{{w_2}^2} - 
     545\,{{w_1}^2}\,{{w_2}^3} + 370\,w_1\,{{w_2}^4} - 
     98\,{{w_2}^5} + 620\,{{w_1}^4}\,w_3 - 
     2960\,{{w_1}^3}\,w_2\,w_3 + 
     5880\,{{w_1}^2}\,{{w_2}^2}\,w_3 - 
     5360\,w_1\,{{w_2}^3}\,w_3 + 
     1820\,{{w_2}^4}\,w_3 +
     6855\,{{w_1}^3}\,{{w_3}^2} - 
     23355\,{{w_1}^2}\,w_2\,{{w_3}^2} + 
     30060\,w_1\,{{w_2}^2}\,{{w_3}^2} - 
     13560\,{{w_2}^3}\,{{w_3}^2} + 
     39120\,{{w_1}^2}\,{{w_3}^3} - 
     84720\,w_1\,w_2\,{{w_3}^3} + 
     53160\,{{w_2}^2}\,{{w_3}^3} +  
     112674\,w_1\,{{w_3}^4} - 114654\,w_2\,{{w_3}^4} + 
     126636\,{{w_3}^5} ) \,Y_2 + $\\[5pt]
$1944\,( 115\,{{w_1}^6} - 1470\,{{w_1}^5}\,w_2 + 
     7245\,{{w_1}^4}\,{{w_2}^2} - 
     17660\,{{w_1}^3}\,{{w_2}^3} + 
     22890\,{{w_1}^2}\,{{w_2}^4} - 
     15120\,w_1\,{{w_2}^5} + 4000\,{{w_2}^6} + 
     5770\,{{w_1}^5}\,w_3 - 
     61550\,{{w_1}^4}\,w_2\,w_3 + 
     239260\,{{w_1}^3}\,{{w_2}^2}\,w_3 - 
     421840\,{{w_1}^2}\,{{w_2}^3}\,w_3 + 
     348920\,w_1\,{{w_2}^4}\,w_3 - 
     110560\,{{w_2}^5}\,w_3 + 
     48435\,{{w_1}^4}\,{{w_3}^2} - 
     503340\,{{w_1}^3}\,w_2\,{{w_3}^2} + 
     1702980\,{{w_1}^2}\,{{w_2}^2}\,{{w_3}^2} - 
     2154480\,w_1\,{{w_2}^3}\,{{w_3}^2} + 
     922200\,{{w_2}^4}\,{{w_3}^2} + 
     7020\,{{w_1}^3}\,{{w_3}^3} - 
     1034640\,{{w_1}^2}\,w_2\,{{w_3}^3} + 
     3974400\,w_1\,{{w_2}^2}\,{{w_3}^3} - 
     3058560\,{{w_2}^3}\,{{w_3}^3} - 
     540630\,{{w_1}^2}\,{{w_3}^4} + 
     1480320\,w_1\,w_2\,{{w_3}^4} + 
     2145600\,{{w_2}^2}\,{{w_3}^4} - 
     2453976\,w_1\,{{w_3}^5} + 11969856\,w_2\,{{w_3}^5} - 
     10703448\,{{w_3}^6} ) \,Y_1\,Y_2 + $\\[5pt]
  $54\,( 2599\,{{w_1}^6} - 44844\,{{w_1}^5}\,w_2 + 
     308760\,{{w_1}^4}\,{{w_2}^2} - 
     1072400\,{{w_1}^3}\,{{w_2}^3} + 
     1947840\,{{w_1}^2}\,{{w_2}^4} - 
     1718016\,w_1\,{{w_2}^5} + 572416\,{{w_2}^6} + 
     59640\,{{w_1}^5}\,w_3 - 
     638760\,{{w_1}^4}\,w_2\,w_3 + 
     2621760\,{{w_1}^3}\,{{w_2}^2}\,w_3 - 
     5860320\,{{w_1}^2}\,{{w_2}^3}\,w_3 + 
     8221440\,w_1\,{{w_2}^4}\,w_3 - 
     4462080\,{{w_2}^5}\,w_3 - 
     4497300\,{{w_1}^4}\,{{w_3}^2} + 
     25488720\,{{w_1}^3}\,w_2\,{{w_3}^2} - 
     63547200\,{{w_1}^2}\,{{w_2}^2}\,{{w_3}^2} + 
     81408960\,w_1\,{{w_2}^3}\,{{w_3}^2} - 
     34493760\,{{w_2}^4}\,{{w_3}^2} - 
     49122720\,{{w_1}^3}\,{{w_3}^3} + 
     119024640\,{{w_1}^2}\,w_2\,{{w_3}^3} - 
     177344640\,w_1\,{{w_2}^2}\,{{w_3}^3} + 
     143251200\,{{w_2}^3}\,{{w_3}^3} - 
     169406640\,{{w_1}^2}\,{{w_3}^4} + 
     94141440\,w_1\,w_2\,{{w_3}^4} - 
     51114240\,{{w_2}^2}\,{{w_3}^4} - 
     2646390528\,w_1\,{{w_3}^5} + 
     3953878272\,w_2\,{{w_3}^5} - 11168202240\,{{w_3}^6}
      ) \,{{Y_1}^2}\,Y_2 + $\\[5pt]

$216\,( 1364\,{{w_1}^6} - 24924\,{{w_1}^5}\,w_2 + 
     182055\,{{w_1}^4}\,{{w_2}^2} - 
     670960\,{{w_1}^3}\,{{w_2}^3} + 
     1287840\,{{w_1}^2}\,{{w_2}^4} - 
     1180416\,w_1\,{{w_2}^5} + 386816\,{{w_2}^6} - 
     112209\,{{w_1}^5}\,w_3 + 
     1371000\,{{w_1}^4}\,w_2\,w_3 - 
     6057780\,{{w_1}^3}\,{{w_2}^2}\,w_3 + 
     11491440\,{{w_1}^2}\,{{w_2}^3}\,w_3 - 
     8506560\,w_1\,{{w_2}^4}\,w_3 + 
     1787136\,{{w_2}^5}\,w_3 - 
     739980\,{{w_1}^4}\,{{w_3}^2} + 
     8601300\,{{w_1}^3}\,w_2\,{{w_3}^2} - 
     39277980\,{{w_1}^2}\,{{w_2}^2}\,{{w_3}^2} + 
     79248960\,w_1\,{{w_2}^3}\,{{w_3}^2} - 
     49596480\,{{w_2}^4}\,{{w_3}^2} + 
     16542000\,{{w_1}^3}\,{{w_3}^3} - 
     94577760\,{{w_1}^2}\,w_2\,{{w_3}^3} + 
     143164800\,w_1\,{{w_2}^2}\,{{w_3}^3} - 
     46719360\,{{w_2}^3}\,{{w_3}^3} - 
     150956460\,{{w_1}^2}\,{{w_3}^4} + 
     1419294960\,w_1\,w_2\,{{w_3}^4} - 
     1405900800\,{{w_2}^2}\,{{w_3}^4} - 
     5687381952\,w_1\,{{w_3}^5} + 
     12549919680\,w_2\,{{w_3}^5} - 19620243792\,{{w_3}^6}
      ) \,{{Y_1}^3}\,Y_2 + $\\[5pt]

$72\,( -2891\,{{w_1}^6} + 50565\,{{w_1}^5}\,w_2 - 
     353910\,{{w_1}^4}\,{{w_2}^2} + 
     1272640\,{{w_1}^3}\,{{w_2}^3} - 
     2556480\,{{w_1}^2}\,{{w_2}^4} + 
     3005184\,w_1\,{{w_2}^5} - 1902080\,{{w_2}^6} - 
     3213\,{{w_1}^5}\,w_3 + 
     371250\,{{w_1}^4}\,w_2\,w_3 - 
     4551120\,{{w_1}^3}\,{{w_2}^2}\,w_3 + 
     21029760\,{{w_1}^2}\,{{w_2}^3}\,w_3 - 
     40711680\,w_1\,{{w_2}^4}\,w_3 + 
     25892352\,{{w_2}^5}\,w_3 + 
     3293460\,{{w_1}^4}\,{{w_3}^2} - 
     25503660\,{{w_1}^3}\,w_2\,{{w_3}^2} + 
     30025080\,{{w_1}^2}\,{{w_2}^2}\,{{w_3}^2} + 
     140849280\,w_1\,{{w_2}^3}\,{{w_3}^2} - 
     254689920\,{{w_2}^4}\,{{w_3}^2} - 
     209584800\,{{w_1}^3}\,{{w_3}^3} + 
     2077646760\,{{w_1}^2}\,w_2\,{{w_3}^3} - 
     5395870080\,w_1\,{{w_2}^2}\,{{w_3}^3} + 
     3266213760\,{{w_2}^3}\,{{w_3}^3} - 
     7122067560\,{{w_1}^2}\,{{w_3}^4} + 
     42669011520\,w_1\,w_2\,{{w_3}^4} - 
     50244896160\,{{w_2}^2}\,{{w_3}^4} - 
     80632905408\,w_1\,{{w_3}^5} + 
     157744600848\,w_2\,{{w_3}^5} -
     205830059040\,{{w_3}^6}
      ) \,{{Y_1}^4}\,Y_2 + $\\[5pt]
$36\,( -8317\,{{w_1}^6} + 192264\,{{w_1}^5}\,w_2 - 
     1849200\,{{w_1}^4}\,{{w_2}^2} + 
     9470720\,{{w_1}^3}\,{{w_2}^3} - 
     27237120\,{{w_1}^2}\,{{w_2}^4} + 
     41699328\,w_1\,{{w_2}^5} - 26546176\,{{w_2}^6} - 
     616698\,{{w_1}^5}\,w_3 + 
     10877580\,{{w_1}^4}\,w_2\,w_3 - 
     75369600\,{{w_1}^3}\,{{w_2}^2}\,w_3 + 
     254874240\,{{w_1}^2}\,{{w_2}^3}\,w_3 - 
     416540160\,w_1\,{{w_2}^4}\,w_3 + 
     258665472\,{{w_2}^5}\,w_3 - 
     44838090\,{{w_1}^4}\,{{w_3}^2} + 
     611442000\,{{w_1}^3}\,w_2\,{{w_3}^2} - 
     3022349760\,{{w_1}^2}\,{{w_2}^2}\,{{w_3}^2} + 
     6308133120\,w_1\,{{w_2}^3}\,{{w_3}^2} - 
     4528673280\,{{w_2}^4}\,{{w_3}^2} - 
     1959857640\,{{w_1}^3}\,{{w_3}^3} + 
     17825326560\,{{w_1}^2}\,w_2\,{{w_3}^3} - 
     50015905920\,w_1\,{{w_2}^2}\,{{w_3}^3} + 
     40289287680\,{{w_2}^3}\,{{w_3}^3} - 
     29008110420\,{{w_1}^2}\,{{w_3}^4} + 
     165895011360\,w_1\,w_2\,{{w_3}^4} - 
     207940737600\,{{w_2}^2}\,{{w_3}^4} - 
     203187801456\,w_1\,{{w_3}^5} + 
     216601319808\,w_2\,{{w_3}^5} -
     165452312016\,{{w_3}^6}
      ) \,{{Y_1}^5}\,Y_2 + $\\[5pt]

$18\,( -22553\,{{w_1}^6} + 541272\,{{w_1}^5}\,w_2 - 
     5412720\,{{w_1}^4}\,{{w_2}^2} + 
     28867840\,{{w_1}^3}\,{{w_2}^3} - 
     86603520\,{{w_1}^2}\,{{w_2}^4} + 
     138565632\,w_1\,{{w_2}^5} - 92377088\,{{w_2}^6} - 
     2102004\,{{w_1}^5}\,w_3 + 
     37238400\,{{w_1}^4}\,w_2\,w_3 - 
     259493760\,{{w_1}^3}\,{{w_2}^2}\,w_3 + 
     884321280\,{{w_1}^2}\,{{w_2}^3}\,w_3 - 
     1461335040\,w_1\,{{w_2}^4}\,w_3 + 
     923222016\,{{w_2}^5}\,w_3 - 
     132943140\,{{w_1}^4}\,{{w_3}^2} + 
     1818849600\,{{w_1}^3}\,w_2\,{{w_3}^2} - 
     9063653760\,{{w_1}^2}\,{{w_2}^2}\,{{w_3}^2} + 
     19237893120\,w_1\,{{w_2}^3}\,{{w_3}^2} - 
     14306042880\,{{w_2}^4}\,{{w_3}^2} - 
     2851640640\,{{w_1}^3}\,{{w_3}^3} + 
     28046165760\,{{w_1}^2}\,w_2\,{{w_3}^3} - 
     87490575360\,w_1\,{{w_2}^2}\,{{w_3}^3} + 
     83728650240\,{{w_2}^3}\,{{w_3}^3} - 
     2171428560\,{{w_1}^2}\,{{w_3}^4} - 
     43665816960\,w_1\,w_2\,{{w_3}^4} + 
     185471596800\,{{w_2}^2}\,{{w_3}^4} + 
     257835472704\,w_1\,{{w_3}^5} - 
     1771101884160\,w_2\,{{w_3}^5} +
     2944386630144\,{{w_3}^6}
      ) \,{{Y_1}^6}\,Y_2 + $\\[5pt]

$324\,( -513\,{{w_1}^6} + 12312\,{{w_1}^5}\,w_2 - 
     123120\,{{w_1}^4}\,{{w_2}^2} + 
     656640\,{{w_1}^3}\,{{w_2}^3} - 
     1969920\,{{w_1}^2}\,{{w_2}^4} + 
     3151872\,w_1\,{{w_2}^5} - 2101248\,{{w_2}^6} - 
     3772\,{{w_1}^5}\,w_3 + 
     75440\,{{w_1}^4}\,w_2\,w_3 - 
     603520\,{{w_1}^3}\,{{w_2}^2}\,w_3 + 
     2414080\,{{w_1}^2}\,{{w_2}^3}\,w_3 - 
     4828160\,w_1\,{{w_2}^4}\,w_3 + 
     3862528\,{{w_2}^5}\,w_3 - 
     753840\,{{w_1}^4}\,{{w_3}^2} + 
     12061440\,{{w_1}^3}\,w_2\,{{w_3}^2} - 
     72368640\,{{w_1}^2}\,{{w_2}^2}\,{{w_3}^2} + 
     192983040\,w_1\,{{w_2}^3}\,{{w_3}^2} - 
     192983040\,{{w_2}^4}\,{{w_3}^2} + 
     79293120\,{{w_1}^3}\,{{w_3}^3} - 
     934669440\,{{w_1}^2}\,w_2\,{{w_3}^3} + 
     3671285760\,w_1\,{{w_2}^2}\,{{w_3}^3} - 
     4805191680\,{{w_2}^3}\,{{w_3}^3} + 
     3802053600\,{{w_1}^2}\,{{w_3}^4} - 
     29682374400\,w_1\,w_2\,{{w_3}^4} + 
     57896640000\,{{w_2}^2}\,{{w_3}^4} + 
     40244283648\,w_1\,{{w_3}^5} - 
     167285958912\,w_2\,{{w_3}^5} +
     282944332224\,{{w_3}^6}
      ) \,{{Y_1}^7}\,Y_2 + $\\[5pt]

$7776\,{{w_3}^2}\,( 59025\,{{w_1}^4} - 
     944400\,{{w_1}^3}\,w_2 + 
     5666400\,{{w_1}^2}\,{{w_2}^2} - 
     15110400\,w_1\,{{w_2}^3} + 15110400\,{{w_2}^4} + 
     4227400\,{{w_1}^3}\,w_3 - 
     50728800\,{{w_1}^2}\,w_2\,w_3 + 
     202915200\,w_1\,{{w_2}^2}\,w_3 - 
     270553600\,{{w_2}^3}\,w_3 + 
     86580180\,{{w_1}^2}\,{{w_3}^2} - 
     692641440\,w_1\,w_2\,{{w_3}^2} + 
     1385282880\,{{w_2}^2}\,{{w_3}^2} + 
     747341424\,w_1\,{{w_3}^3} - 
     2882445696\,w_2\,{{w_3}^3} + 7024285008\,{{w_3}^4} )
    \,{{Y_1}^8}\,Y_2 + $\\[5pt]
$69984000\,{{w_3}^4}\,( -185\,{{w_1}^2} + 
     1480\,w_1\,w_2 - 2960\,{{w_2}^2} - 
     1812\,w_1\,w_3 + 7248\,w_2\,w_3 + 
     171480\,{{w_3}^2} ) \,{{Y_1}^9}\,Y_2 - $\\[5pt]
$629856000000\,{{w_3}^6}\,{{Y_1}^{10}}\,Y_2 + $\\[5pt]
$17496\,( -35\,{{w_1}^6} + 510\,{{w_1}^5}\,w_2 - 
     2835\,{{w_1}^4}\,{{w_2}^2} + 
     7480\,{{w_1}^3}\,{{w_2}^3} - 
     9600\,{{w_1}^2}\,{{w_2}^4} + 5760\,w_1\,{{w_2}^5} - 
     1280\,{{w_2}^6} - 1370\,{{w_1}^5}\,w_3 + 
     16810\,{{w_1}^4}\,w_2\,w_3 - 
     75680\,{{w_1}^3}\,{{w_2}^2}\,w_3 + 
     151760\,{{w_1}^2}\,{{w_2}^3}\,w_3 - 
     131200\,w_1\,{{w_2}^4}\,w_3 + 
     39680\,{{w_2}^5}\,w_3 - 8925\,{{w_1}^4}\,{{w_3}^2} + 
     114360\,{{w_1}^3}\,w_2\,{{w_3}^2} - 
     479880\,{{w_1}^2}\,{{w_2}^2}\,{{w_3}^2} + 
     751200\,w_1\,{{w_2}^3}\,{{w_3}^2} - 
     360960\,{{w_2}^4}\,{{w_3}^2} + 
     14400\,{{w_1}^3}\,{{w_3}^3} + 
     108000\,{{w_1}^2}\,w_2\,{{w_3}^3} - 
     855360\,w_1\,{{w_2}^2}\,{{w_3}^3} + 
     1005120\,{{w_2}^3}\,{{w_3}^3} + 
     159480\,{{w_1}^2}\,{{w_3}^4} - 
     835200\,w_1\,w_2\,{{w_3}^4} + 
     147600\,{{w_2}^2}\,{{w_3}^4} + 
     363744\,w_1\,{{w_3}^5} - 2877984\,w_2\,{{w_3}^5} + 
     873072\,{{w_3}^6} ) \,{{Y_2}^2} + $\\[5pt]
$972\,( 109\,{{w_1}^6} - 2004\,{{w_1}^5}\,w_2 + 
     14730\,{{w_1}^4}\,{{w_2}^2} - 
     54560\,{{w_1}^3}\,{{w_2}^3} + 
     104640\,{{w_1}^2}\,{{w_2}^4} - 
     93696\,w_1\,{{w_2}^5} + 27136\,{{w_2}^6} + 
     840\,{{w_1}^5}\,w_3 - 
     14640\,{{w_1}^4}\,w_2\,w_3 + 
     99840\,{{w_1}^3}\,{{w_2}^2}\,w_3 - 
     330240\,{{w_1}^2}\,{{w_2}^3}\,w_3 + 
     522240\,w_1\,{{w_2}^4}\,w_3 - 
     307200\,{{w_2}^5}\,w_3 + 
     694170\,{{w_1}^4}\,{{w_3}^2} - 
     6246720\,{{w_1}^3}\,w_2\,{{w_3}^2} + 
     17943120\,{{w_1}^2}\,{{w_2}^2}\,{{w_3}^2} - 
     17210880\,w_1\,{{w_2}^3}\,{{w_3}^2} + 
     3836160\,{{w_2}^4}\,{{w_3}^2} - 
     2509920\,{{w_1}^3}\,{{w_3}^3} - 
     2617920\,{{w_1}^2}\,w_2\,{{w_3}^3} + 
     53239680\,w_1\,{{w_2}^2}\,{{w_3}^3} - 
     35631360\,{{w_2}^3}\,{{w_3}^3} - 
     80261280\,{{w_1}^2}\,{{w_3}^4} + 
     220086720\,w_1\,w_2\,{{w_3}^4} - 
     40214880\,{{w_2}^2}\,{{w_3}^4} + 
     17200512\,w_1\,{{w_3}^5} - 
     151103232\,w_2\,{{w_3}^5} + 977186592\,{{w_3}^6} ) 
    \,Y_1\,{{Y_2}^2} + $\\[5pt]$1944\,( -91\,{{w_1}^6} +
$\\[5pt]
$1986\,{{w_1}^5}\,w_2 - 
     17880\,{{w_1}^4}\,{{w_2}^2} + 
     84800\,{{w_1}^3}\,{{w_2}^3} - 
     222720\,{{w_1}^2}\,{{w_2}^4} + 
     305664\,w_1\,{{w_2}^5} - 169984\,{{w_2}^6} + 
     28158\,{{w_1}^5}\,w_3 - 
     381450\,{{w_1}^4}\,w_2\,w_3 + 
     1889520\,{{w_1}^3}\,{{w_2}^2}\,w_3 - 
     4076160\,{{w_1}^2}\,{{w_2}^3}\,w_3 + 
     3521280\,w_1\,{{w_2}^4}\,w_3 - 
     978432\,{{w_2}^5}\,w_3 + 
     328455\,{{w_1}^4}\,{{w_3}^2} - 
     3800520\,{{w_1}^3}\,w_2\,{{w_3}^2} + 
     14949360\,{{w_1}^2}\,{{w_2}^2}\,{{w_3}^2} - 
     21254400\,w_1\,{{w_2}^3}\,{{w_3}^2} + 
     4976640\,{{w_2}^4}\,{{w_3}^2} - 
     3570480\,{{w_1}^3}\,{{w_3}^3} + 
     17496000\,{{w_1}^2}\,w_2\,{{w_3}^3} - 
     14074560\,w_1\,{{w_2}^2}\,{{w_3}^3} + 
     4872960\,{{w_2}^3}\,{{w_3}^3} + 
     192588840\,{{w_1}^2}\,{{w_3}^4} - 
     980372160\,w_1\,w_2\,{{w_3}^4} + 
     793877760\,{{w_2}^2}\,{{w_3}^4} + 
     2910144672\,w_1\,{{w_3}^5} - 
     6921487584\,w_2\,{{w_3}^5} + 6774579504\,{{w_3}^6} )
    \,{{Y_1}^2}\,{{Y_2}^2} + $\\[5pt]
$ 324\,( -1777\,{{w_1}^6} + 38868\,{{w_1}^5}\,w_2 - 
     350880\,{{w_1}^4}\,{{w_2}^2} + 
     1669760\,{{w_1}^3}\,{{w_2}^3} - 
     4404480\,{{w_1}^2}\,{{w_2}^4} + 
     6079488\,w_1\,{{w_2}^5} - 3407872\,{{w_2}^6} + 
     5328\,{{w_1}^5}\,w_3 - 
     145440\,{{w_1}^4}\,w_2\,w_3 + 
     1474560\,{{w_1}^3}\,{{w_2}^2}\,w_3 - 
     7142400\,{{w_1}^2}\,{{w_2}^3}\,w_3 + 
     16773120\,w_1\,{{w_2}^4}\,w_3 - 
     15409152\,{{w_2}^5}\,w_3 - 
     1347030\,{{w_1}^4}\,{{w_3}^2} + 
     20308320\,{{w_1}^3}\,w_2\,{{w_3}^2} - 
     114384960\,{{w_1}^2}\,{{w_2}^2}\,{{w_3}^2} + 
     285120000\,w_1\,{{w_2}^3}\,{{w_3}^2} - 
     265213440\,{{w_2}^4}\,{{w_3}^2} + 
     141959520\,{{w_1}^3}\,{{w_3}^3} - 
     1184232960\,{{w_1}^2}\,w_2\,{{w_3}^3} + 
     2735389440\,w_1\,{{w_2}^2}\,{{w_3}^3} - 
     1079239680\,{{w_2}^3}\,{{w_3}^3} + 
     5862792960\,{{w_1}^2}\,{{w_3}^4} - 
     31292412480\,w_1\,w_2\,{{w_3}^4} + 
     34249703040\,{{w_2}^2}\,{{w_3}^4} + 
     54443259648\,w_1\,{{w_3}^5} - 
     118126460160\,w_2\,{{w_3}^5} + 73898741664\,{{w_3}^6}
      ) \,{{Y_1}^3}\,{{Y_2}^2} + $\\[5pt]
$17496\,( 19\,{{w_1}^6} - 456\,{{w_1}^5}\,w_2 + 
     4560\,{{w_1}^4}\,{{w_2}^2} - 
     24320\,{{w_1}^3}\,{{w_2}^3} + 
     72960\,{{w_1}^2}\,{{w_2}^4} - 
     116736\,w_1\,{{w_2}^5} + 77824\,{{w_2}^6} + 
     676\,{{w_1}^5}\,w_3 - 
     13520\,{{w_1}^4}\,w_2\,w_3 + 
     108160\,{{w_1}^3}\,{{w_2}^2}\,w_3 - 
     432640\,{{w_1}^2}\,{{w_2}^3}\,w_3 + 
     865280\,w_1\,{{w_2}^4}\,w_3 - 
     692224\,{{w_2}^5}\,w_3 + 
     47970\,{{w_1}^4}\,{{w_3}^2} - 
     601920\,{{w_1}^3}\,w_2\,{{w_3}^2} + 
     2617920\,{{w_1}^2}\,{{w_2}^2}\,{{w_3}^2} - 
     4331520\,w_1\,{{w_2}^3}\,{{w_3}^2} + 
     1681920\,{{w_2}^4}\,{{w_3}^2} + 
     4475280\,{{w_1}^3}\,{{w_3}^3} - 
     40034880\,{{w_1}^2}\,w_2\,{{w_3}^3} + 
     105465600\,w_1\,{{w_2}^2}\,{{w_3}^3} - 
     67722240\,{{w_2}^3}\,{{w_3}^3} + 
     91656720\,{{w_1}^2}\,{{w_3}^4} - 
     536624640\,w_1\,w_2\,{{w_3}^4} + 
     679991040\,{{w_2}^2}\,{{w_3}^4} + 
     118710144\,w_1\,{{w_3}^5} + 
     969224832\,w_2\,{{w_3}^5} - 3877256160\,{{w_3}^6} ) 
    \,{{Y_1}^4}\,{{Y_2}^2} + $\\[5pt]
$23328\,w_3\,( 199\,{{w_1}^5} - 
     3980\,{{w_1}^4}\,w_2 + 31840\,{{w_1}^3}\,{{w_2}^2} - 
     127360\,{{w_1}^2}\,{{w_2}^3} + 
     254720\,w_1\,{{w_2}^4} - 203776\,{{w_2}^5} - 
     5760\,{{w_1}^4}\,w_3 + 
     92160\,{{w_1}^3}\,w_2\,w_3 - 
     552960\,{{w_1}^2}\,{{w_2}^2}\,w_3 + 
     1474560\,w_1\,{{w_2}^3}\,w_3 - 
     1474560\,{{w_2}^4}\,w_3 + 
     43200\,{{w_1}^3}\,{{w_3}^2} - 
     518400\,{{w_1}^2}\,w_2\,{{w_3}^2} + 
     2073600\,w_1\,{{w_2}^2}\,{{w_3}^2} - 
     2764800\,{{w_2}^3}\,{{w_3}^2} - 
     48576240\,{{w_1}^2}\,{{w_3}^3} + 
     381028320\,w_1\,w_2\,{{w_3}^3} - 
     746893440\,{{w_2}^2}\,{{w_3}^3} - 
     1294311312\,w_1\,{{w_3}^4} + 
     5203372608\,w_2\,{{w_3}^4} - 8225713296\,{{w_3}^5} )
    \,{{Y_1}^5}\,{{Y_2}^2} + $\\[5pt]
$1119744\,{{w_3}^3}\,( -25775\,{{w_1}^3} + 
     309300\,{{w_1}^2}\,w_2 - 1237200\,w_1\,{{w_2}^2} + 
     1649600\,{{w_2}^3} - 1342845\,{{w_1}^2}\,w_3 + 
     10742760\,w_1\,w_2\,w_3 - 
     21485520\,{{w_2}^2}\,w_3 - 
     18125964\,w_1\,{{w_3}^2} + 
     72503856\,w_2\,{{w_3}^2} - 114744438\,{{w_3}^3} ) \,
   {{Y_1}^6}\,{{Y_2}^2} + $\\[5pt]
$15116544000\,( 11\,w_1 - 44\,w_2 - 972\,w_3 ) \,
   {{w_3}^5}\,{{Y_1}^7}\,{{Y_2}^2} + $\\[5pt]
$ 5832\,( -41\,{{w_1}^6} + 876\,{{w_1}^5}\,w_2 - 
     7680\,{{w_1}^4}\,{{w_2}^2} + 
     35200\,{{w_1}^3}\,{{w_2}^3} - 
     88320\,{{w_1}^2}\,{{w_2}^4} + 
     113664\,w_1\,{{w_2}^5} - 57344\,{{w_2}^6} - 
     134460\,{{w_1}^4}\,{{w_3}^2} + 
     1257120\,{{w_1}^3}\,w_2\,{{w_3}^2} - 
     3576960\,{{w_1}^2}\,{{w_2}^2}\,{{w_3}^2} + 
     2695680\,w_1\,{{w_2}^3}\,{{w_3}^2} + 
     414720\,{{w_2}^4}\,{{w_3}^2} + 
     432000\,{{w_1}^3}\,{{w_3}^3} + 
     1347840\,{{w_1}^2}\,w_2\,{{w_3}^3} - 
     14722560\,w_1\,{{w_2}^2}\,{{w_3}^3} + 
     9676800\,{{w_2}^3}\,{{w_3}^3} + 
     13433040\,{{w_1}^2}\,{{w_3}^4} - 
     40979520\,w_1\,w_2\,{{w_3}^4} - 
     622080\,{{w_2}^2}\,{{w_3}^4} - 
     5971968\,w_1\,{{w_3}^5} + 17169408\,w_2\,{{w_3}^5} - 
     148039488\,{{w_3}^6} ) \,{{Y_2}^3} + $\\[5pt]
$629856\,w_3\,( 59\,{{w_1}^5} - 
     1090\,{{w_1}^4}\,w_2 + 8000\,{{w_1}^3}\,{{w_2}^2} - 
     29120\,{{w_1}^2}\,{{w_2}^3} + 
     52480\,w_1\,{{w_2}^4} - 37376\,{{w_2}^5} + 
     920\,{{w_1}^4}\,w_3 - 
     14000\,{{w_1}^3}\,w_2\,w_3 + 
     79680\,{{w_1}^2}\,{{w_2}^2}\,w_3 - 
     200960\,w_1\,{{w_2}^3}\,w_3 + 
     189440\,{{w_2}^4}\,w_3 + 
     5000\,{{w_1}^3}\,{{w_3}^2} - 
     57840\,{{w_1}^2}\,w_2\,{{w_3}^2} + 
     222720\,w_1\,{{w_2}^2}\,{{w_3}^2} - 
     285440\,{{w_2}^3}\,{{w_3}^2} - 
     1573920\,{{w_1}^2}\,{{w_3}^3} + 
     7381440\,w_1\,w_2\,{{w_3}^3} - 
     5276160\,{{w_2}^2}\,{{w_3}^3} - 
     17284752\,w_1\,{{w_3}^4} + 
     45336672\,w_2\,{{w_3}^4} - 29611008\,{{w_3}^5} ) \,
   Y_1\,{{Y_2}^3} + $\\[5pt]
$ 52488\,( 19\,{{w_1}^6} - 456\,{{w_1}^5}\,w_2 + 
     4560\,{{w_1}^4}\,{{w_2}^2} - 
     24320\,{{w_1}^3}\,{{w_2}^3} + 
     72960\,{{w_1}^2}\,{{w_2}^4} - 
     116736\,w_1\,{{w_2}^5} + 77824\,{{w_2}^6} + 
     244\,{{w_1}^5}\,w_3 - 
     4880\,{{w_1}^4}\,w_2\,w_3 + 
     39040\,{{w_1}^3}\,{{w_2}^2}\,w_3 - 
     156160\,{{w_1}^2}\,{{w_2}^3}\,w_3 + 
     312320\,w_1\,{{w_2}^4}\,w_3 - 
     249856\,{{w_2}^5}\,w_3 + 
     18300\,{{w_1}^4}\,{{w_3}^2} - 
     292800\,{{w_1}^3}\,w_2\,{{w_3}^2} + 
     1756800\,{{w_1}^2}\,{{w_2}^2}\,{{w_3}^2} - 
     4684800\,w_1\,{{w_2}^3}\,{{w_3}^2} + 
     4684800\,{{w_2}^4}\,{{w_3}^2} - 
     708960\,{{w_1}^3}\,{{w_3}^3} + 
     5863680\,{{w_1}^2}\,w_2\,{{w_3}^3} - 
     12879360\,w_1\,{{w_2}^2}\,{{w_3}^3} + 
     3072000\,{{w_2}^3}\,{{w_3}^3} - 
     36048240\,{{w_1}^2}\,{{w_3}^4} + 
     196007040\,w_1\,w_2\,{{w_3}^4} - 
     207256320\,{{w_2}^2}\,{{w_3}^4} - 
     308650176\,w_1\,{{w_3}^5} + 
     759456000\,w_2\,{{w_3}^5} + 447550272\,{{w_3}^6} ) 
    \,{{Y_1}^2}\,{{Y_2}^3} + $\\[5pt]
$1259712\,{{w_3}^2}\,( 125\,{{w_1}^4} - 
     2000\,{{w_1}^3}\,w_2 + 12000\,{{w_1}^2}\,{{w_2}^2} - 
     32000\,w_1\,{{w_2}^3} + 32000\,{{w_2}^4} + 
     1000\,{{w_1}^3}\,w_3 - 
     12000\,{{w_1}^2}\,w_2\,w_3 + 
     48000\,w_1\,{{w_2}^2}\,w_3 - 
     64000\,{{w_2}^3}\,w_3 - 
     60840\,{{w_1}^2}\,{{w_3}^2} + 
     486720\,w_1\,w_2\,{{w_3}^2} - 
     973440\,{{w_2}^2}\,{{w_3}^2} + 
     17240256\,w_1\,{{w_3}^3} - 
     68416704\,w_2\,{{w_3}^3} + 152725392\,{{w_3}^4} ) \,
   {{Y_1}^3}\,{{Y_2}^3} + $\\[5pt]
$ 90699264\,{{w_3}^4}\,( 11215\,{{w_1}^2} - 
     89720\,w_1\,w_2 + 179440\,{{w_2}^2} + 
     329676\,w_1\,w_3 - 1318704\,w_2\,w_3 + 
     1922532\,{{w_3}^2} ) \,{{Y_1}^4}\,{{Y_2}^3}
+$\\[5pt] 
$5714053632000\,{{w_3}^6}\,{{Y_1}^5}\,{{Y_2}^3}
+$\\[5pt] 
$45349632\,w_3\,( -{{w_1}^5} + 20\,{{w_1}^4}\,w_2 - 
     160\,{{w_1}^3}\,{{w_2}^2} + 
     640\,{{w_1}^2}\,{{w_2}^3} - 1280\,w_1\,{{w_2}^4} + 
     1024\,{{w_2}^5} - 10\,{{w_1}^4}\,w_3 + 
     160\,{{w_1}^3}\,w_2\,w_3 - 
     960\,{{w_1}^2}\,{{w_2}^2}\,w_3 + 
     2560\,w_1\,{{w_2}^3}\,w_3 - 
     2560\,{{w_2}^4}\,w_3 - 40\,{{w_1}^3}\,{{w_3}^2} + 
     480\,{{w_1}^2}\,w_2\,{{w_3}^2} - 
     1920\,w_1\,{{w_2}^2}\,{{w_3}^2} + 
     2560\,{{w_2}^3}\,{{w_3}^2} + 
     16560\,{{w_1}^2}\,{{w_3}^3} - 
     80640\,w_1\,w_2\,{{w_3}^3} + 
     57600\,{{w_2}^2}\,{{w_3}^3} + 
     168048\,w_1\,{{w_3}^4} - 485568\,w_2\,{{w_3}^4} + 
     235872\,{{w_3}^5} ) \,{{Y_2}^4} - $\\[5pt]
$81498730659840\,{{w_3}^6}\,Y_1\,{{Y_2}^4} + $\\[5pt]
$78364164096\,( -201\,w_1 + 804\,w_2 - 1682\,w_3 )
    \,{{w_3}^5}\,{{Y_1}^2}\,{{Y_2}^4} +$\\[5pt] 
$ 45137758519296\,{{w_3}^6}\,{{Y_2}^5}$.
\end{sloppypar}

%% file: FInV.tex
\noindent  
$F_V(w) =
\displaystyle{
\frac{F(\sigma_z(w))|_{\{F(z)=0\}}}
     {3^{18}\,\Phi(z)^2\,\Psi(z)^9}
}
=$
\vspace{10pt}
\begin{sloppypar} \raggedright \noindent
$9\,w_2^5\,(-24\,w_1 + 19\,w_2 + 24\,w_3) +
$\\[5pt]
$w_2^4\,(365\,w_1^2 - 54\,w_1\,w_2 - 171\,w_2^2 -
640\,w_1\,w_3 -
16\,w_2\,w_3 + 275\,w_3^2)\,V + $\\[5pt] 
$(128\,w_1^6 - 328\,w_1^5\,w_2 +
165\,w_1^4\,w_2^2 -
20\,w_1^3\,w_2^3 - 365\,w_1^2\,w_2^4 +
270\,w_1\,w_2^5 -
768\,w_1^5\,w_3 + 2120\,w_1^4\,w_2\,w_3 -
2320\,w_1^3\,w_2^2\,w_3 +
2000\,w_1^2\,w_2^3\,w_3 - 360\,w_1\,w_2^4\,w_3 +
1920\,w_1^4\,w_3^2 -
5200\,w_1^3\,w_2\,w_3^2 +
5850\,w_1^2\,w_2^2\,w_3^2 -
3800\,w_1\,w_2^3\,w_3^2 + 675\,w_2^4\,w_3^2 -
2560\,w_1^3\,w_3^3 +
6160\,w_1^2\,w_2\,w_3^3 - 5400\,w_1\,w_2^2\,w_3^3
+ 1800\,w_2^3\,w_3^3
+ 1920\,w_1^2\,w_3^4 - 3560\,w_1\,w_2\,w_3^4 +
1705\,w_2^2\,w_3^4 -
768\,w_1\,w_3^5 + 808\,w_2\,w_3^5 +
128\,w_3^6)\,V^2 +$\\[5pt]
$5\,(-25\,w_1^6 + 66\,w_1^5\,w_2 -
33\,w_1^4\,w_2^2 +
4\,w_1^3\,w_2^3 + 112\,w_1^5\,w_3 -
224\,w_1^4\,w_2\,w_3 +
64\,w_1^3\,w_2^2\,w_3 - 193\,w_1^4\,w_3^2 +
280\,w_1^3\,w_2\,w_3^2 -
30\,w_1^2\,w_2^2\,w_3^2 + 152\,w_1^3\,w_3^3 -
152\,w_1^2\,w_2\,w_3^3 -
43\,w_1^2\,w_3^4 + 30\,w_1\,w_2\,w_3^4 -
8\,w_1\,w_3^5 + 5\,w_3^6)\,V^3 +
$\\[5pt]
$w_1^4\,(-3\,w_1^2 - 2\,w_1\,w_2 + 8\,w_1\,w_3 -
5\,w_3^2)\,V^4$.
 \end{sloppypar}

%% file: GammabarInY.tex
\noindent  
$\Bar{\Gamma}_Y(w) = 
  \displaystyle{\frac{5184\,F(\tau_z(w))}
               {(39-\sqrt{15}\,i)\,F(z)^{42}}} =$
   \vspace{10pt}
\begin{sloppypar} \raggedright \noindent
$((3 + \sqrt{15}\,i)\,w_1^{10} - 10\,i\,(-3\,i +
\sqrt{15})\,w_1^9\,w_2 +
     45\,(3 + \sqrt{15}\,i)\,w_1^8\,w_2^2 -
120\,i\,(-3\,i +
     \sqrt{15})\,w_1^7\,w_2^3 + 210\,(3 +
\sqrt{15}\,i)\,w_1^6\,w_2^4 -
     252\,i\,(-3\,i + \sqrt{15})\,w_1^5\,w_2^5 +
210\,(3 +
     \sqrt{15}\,i)\,w_1^4\,w_2^6 - 120\,i\,(-3\,i
+
     \sqrt{15})\,w_1^3\,w_2^7 + 45\,(3 +
\sqrt{15}\,i)\,w_1^2\,w_2^8 -
     10\,i\,(-3\,i + \sqrt{15})\,w_1\,w_2^9 + (3 +
     \sqrt{15}\,i)\,w_2^{10} + 10\,(27 +
5\,\sqrt{15}\,i)\,w_1^9\,w_3 -
     90\,i\,(-27\,i +
5\,\sqrt{15})\,w_1^8\,w_2\,w_3 + 360\,(27 +
     5\,\sqrt{15}\,i)\,w_1^7\,w_2^2\,w_3 -
840\,i\,(-27\,i +
     5\,\sqrt{15})\,w_1^6\,w_2^3\,w_3 +$\\[20pt]
$+ \cdots +$\\[20pt]
$4608\,i\,(5\,(9\,i +
   17\,\sqrt{15})\,w_1^2\,w_3^8 -
     40\,i\,(9 -
17\,\sqrt{15}\,i)\,w_1\,w_2\,w_3^8 + 80\,(9\,i +
     17\,\sqrt{15})\,w_2^2\,w_3^8 + 60\,(15\,i +
     7\,\sqrt{15})\,w_1\,w_3^9 - 240\,i\,(15 -
     7\,\sqrt{15}\,i)\,w_2\,w_3^9 + 36\,(9\,i +
     17\,\sqrt{15})\,w_3^{10})\,Y_2^8 +
663552\,(3 +
  5\,\sqrt{15}\,i)\,w_3^{10}\,Y_1\,Y_2^8$.
\end{sloppypar}

%% file: ThetabarInV.tex
\noindent  
$\Bar{\Theta}_V(w) =
\displaystyle{
\frac{2^7 3^4\,\Bar{\Theta}_z(w)|_{\{F(z)=0\}}}
{(39-\sqrt{15}\,i)\,\Phi(z)\,\Psi(z)^{16}}} =$
   \vspace{10pt}
\begin{sloppypar} \raggedright \noindent
$1944\,w_2^{10} -
 i\,(12960\,(-3\,i + \sqrt{15})\,w_1^3\,w_2^7 + 
     135\,i\,(765 +
203\,\sqrt{15}\,i)\,w_1^2\,w_2^8 + 
     90\,(-1033\,i + 228\,\sqrt{15})\,w_1\,w_2^9
+ 
     i\,(18837 + 5500\,\sqrt{15}\,i)\,w_2^{10} + 
     38880\,i\,(3 +
\sqrt{15}\,i)\,w_1^2\,w_2^7\,w_3 + 
     270\,(-693\,i +
227\,\sqrt{15})\,w_1\,w_2^8\,w_3 + 
     450\,i\,(167 + 57\,\sqrt{15}\,i)\,w_2^9\,w_3
+ 
     38880\,(-3\,i +
\sqrt{15})\,w_1\,w_2^7\,w_3^2 + 
     135\,i\,(621 +
251\,\sqrt{15}\,i)\,w_2^8\,w_3^2 + 
     12960\,i\,(3 +
\sqrt{15}\,i)\,w_2^7\,w_3^3)\,V + $\\[20pt] 
$+ \ldots +$\\[20pt]
$5\,(2\,(46 + \sqrt{15}\,i)\,w_1^{10} +
2\,i\,(38\,i +
  \sqrt{15})\,w_1^9\,w_2 +
     (8 + \sqrt{15}\,i)\,w_1^8\,w_2^2 -
2\,i\,(-171\,i +
     5\,\sqrt{15})\,w_1^9\,w_3 + 6\,(25 -
\sqrt{15}\,i)\,w_1^8\,w_2\,w_3
     + (473 + 19\,\sqrt{15}\,i)\,w_1^8\,w_3^2 +
4\,i\,(19\,i +
     \sqrt{15})\,w_1^7\,w_2\,w_3^2 -
16\,i\,(-18\,i +
     \sqrt{15})\,w_1^7\,w_3^3 + 5\,(13 +
\sqrt{15}\,i)\,w_1^6\,w_3^4)V^7
     + w_1^{10}\,V^8$.
\end{sloppypar}